\documentclass[final]{elsarticle}
\makeatletter
\def\ps@pprintTitle{%
 \let\@oddhead\@empty
 \let\@evenhead\@empty
 \def\@oddfoot{}%
 \let\@evenfoot\@oddfoot}
\makeatother
\usepackage{hyperref}
\usepackage[dvipsnames]{xcolor}
\usepackage{amsmath}
\usepackage{amsbsy}
\usepackage{amsthm}
\usepackage{amssymb}
\usepackage[tmargin=1in,bmargin=1in]{geometry}
\usepackage{verbatim}
\usepackage{graphicx}
\usepackage{latexsym}
\usepackage{rotating}
\usepackage{color}
\usepackage{caption}
\usepackage{url}
 \global\bibsep=0pt
\usepackage{multirow}
\usepackage{soul}
\usepackage{framed,multirow}
\usepackage{xcolor} 
\usepackage{booktabs}
\usepackage{array}
\usepackage{arydshln}
\usepackage[normalem]{ulem}
\usepackage[shortlabels]{enumitem}

\hypersetup{
    colorlinks=true,      		% false: boxed links; true: colored links
    linkcolor=blue,       		% color of internal links
    citecolor=blue,       		% color of links to bibliography
    filecolor=black,      		% color of file links
    urlcolor=red,       		% color of external links
    bookmarks=false,
    pdffitwindow=true,
    pdfpagelayout=SinglePage
}

\newcommand{\ie}{{\it i.e. }}
\newcommand{\eg}{{\it e.g. }}

\newcommand{\vertiii}[1]{{\left\vert\kern-0.25ex\left\vert\kern-0.25ex\left\vert #1 
    \right\vert\kern-0.25ex\right\vert\kern-0.25ex\right\vert}}

\definecolor{newcolor}{rgb}{.8,.349,.1}

\theoremstyle{plain}

\begin{document}

\begin{frontmatter}

  \title{Least SQuares Discretizations (LSQD): a robust and versatile meshless paradigm for solving elliptic PDEs}
 \journal{Journal of Computational Physics}

\address[UCMERCED]{Department of Applied Mathematics, University of California, Merced, California 95343, USA.}
\address[UCDAVIS]{Department of Mathematics, University of California, Davis, California 95616, USA.}
\author[UCMERCED]{Anna Kucherova} 
\author[UCMERCED]{Gbocho M. Terasaki} 
\author[UCDAVIS]{Selma Strango} 
\author[UCMERCED]{Maxime Theillard} 

\cortext[cor]{Corresponding author: mtheillard@ucmerced.edu}

\begin{abstract}

Searching for numerical methods that combine facility and efficiency, while remaining accurate and versatile, is critical. Often, irregular geometries challenge traditional methods that rely on structured or body-fitted meshes. Meshless methods mitigate these issues but oftentimes require the weak formulation which involves defining quadrature rules over potentially intricate geometries. To overcome these challenges, we propose the Least Squares Discretization (LSQD) method. This novel approach simplifies the application of meshless methods by eliminating the need for a weak formulation and necessitates minimal numerical analysis. It offers significant advantages in terms of ease of implementation and adaptability to complex geometries. In this paper, we demonstrate the efficacy of the LSQD method in solving elliptic partial differential equations for a variety of boundary conditions, geometries, and data layouts. We monitor h-P convergence across these parameters/settings and construct an {\it a posteriori} built-in error estimator to establish our method as a robust and accessible numerical alternative.

\end{abstract}


\end{frontmatter}

\section{Introduction}

In the ever-evolving landscape of computational science and engineering, the quest for accurate, efficient, and versatile numerical methods remains paramount. Traditional numerical approaches, such as finite element methods, finite difference methods, and finite volume methods, have long dominated the sphere of simulating complex physical phenomena. However, their reliance on elements, structured grids, or finite volumes often presents challenges when handling irregular geometries, dynamic domains, or intricate boundary conditions. In response to these limitations, a trailblazing paradigm emerged: meshless methods \cite{Gingold1977, LUCY1977}.

Gaining traction almost two decades post-emergence \cite{Atluri1998, Atluri2000, BELYTSCHKO19963, Belyt1994, KANSA1990147, Liu1995}, meshless methods have witnessed a surge in popularity across a multitude of disciplines, including computational fluid dynamics, structural mechanics, electromagnetics, and biomedical engineering \cite{MAZHAR202114, Kumar2021, HE2020225, JOLDES2019152}. Their inherent attributes (\eg{mesh generation simplicity, capacity to handle large deformations/distortions and complex geometries \cite{KUMAR2019858}, etc.}), have proven a versatile alternative to the traditional discretization techniques.

Some of the pioneering meshless methods fall into the following categories: kernel-based \cite{Gingold1977,Liu1995,Zhu1998}, Galerkin \cite{BELYTSCHKO19963,BELYTSCHKO2000385,Duarte1996,Atluri1998,Sukumar1998,Atluri2000,Gonzalez2009,Yvonnet2004,Guittet2015ASP,DU20033933,BABUSKA1997} , and Finite Difference methods \cite{LISZKA198083,ONATE1996}. Meshless Galerkin methods such as the Element-free Galerkin \cite{BELYTSCHKO19963,BELYTSCHKO2000385}, h-p clouds and h-p methods \cite{Duarte1996}, Meshless Local Petrov-Galerkin  \cite{Atluri1998,Atluri2000}, the Natural Element Method \cite{Sukumar1998,Gonzalez2009,Yvonnet2004}, the Voronoi finite volume method \cite{Guittet2015ASP,DU20033933}, and partition of unity method \cite{BABUSKA1997}, all rely on the weak formulation of the partial differential equation (PDE). Thus, although they are meshless, these methods are not truly ``geometry-free" as the weak formulation involves spatial integrals. This requires approximations of definite integrals, or quadrature rules, over - potentially - intricate geometries. These geometries, themselves, may reap extravagant computational costs (\eg Voronoi partition for the Meshless Natural Element method \cite{Sukumar1998,Gonzalez2009,Yvonnet2004} or any Voronoi Finite Volume Method \cite{Guittet2015ASP,DU20033933}). Undeniably, the theoretical foundation for the aforementioned methods allows for a comprehensive mathematical study but may serve as a deterrent for interested non-experts. Clearly, it would be a quixotic undertaking to create an encyclopedic inventory for the trajectory of all meshless methods, so for an exhaustive analysis, we direct interested readers to the following reviews \cite{NGUYEN2008763, Vivek2020, Antonio2023}.

Although commonly used in scattered data interpolation \cite{Franke1982ScatteredDI}, Kansa became one of the first to utilize radial basis functions (RBF) to solve fluid dynamics problems \cite{KANSA1990147}. Theoretical condition number estimates indicated that some meshless methods using RBFs may be unstable \cite{DUAN200866}. Fortunately, Kansa was able to successfully circumvent the ill-conditioned system by using improved multi-quadrics and adaptive techniques \cite{KANSA2000123}. Kansa's method represents the solution as a sum of global radial functions, centered at a set of points. The coefficients in this decomposition are obtained by interpolating the PDE at these points. While in general, the interpolation system can be ill-posed (Mairhuber-Curtis theorem \cite{Mairhuber1956,curtis1959}), when using appropriate radial functions, this system can be proven to be non-singular \cite{Micchelli1986,IMM2005-03600}. Not solely is it straightforward to implement, this method is particularly well-suited for problems involving complex curved surfaces (\eg{such as reaction-diffusion patterns over biological entities such as frog species \cite{PIRET20124662}, red blood cells \cite{Fuselier2013AHK}, or 3-D spherical shells \cite{Wright2010} used in geophysical calculations}). It has since been applied in the computational fluid dynamics realm \cite{cmes.2005.007.185} such as the solution of Navier–Stokes equations \cite{doi:10.1002/fld.165}, porous media flow \cite{Sarler2004}, and the solution of solid-liquid phase change problems \cite{Kovacevic}. Unfortunately, the convergence of Kansa's method is still unproven. Because the representation is global, the interpolation system is full and consequently its inversion rapidly becomes prohibitive. To address this limitation, domain decomposition techniques can be employed \cite{MAIDUY2002133}, but these abandon the meshless nature of the domain. Using local RBFs  \cite{Hardy1971}, the local Kansa method \cite{Sarler2004,SARLER20061269} offers an alternative approach to reduce the complexity. However, this approach still relies on a domain decomposition to distinguish between inner and border nodes, further requiring functional iterations on either node sets. 

One of the most pervasive methods in computational fluid dynamics is the finite point method - a meshfree method for solving PDEs on scattered distributions of points. This method is considered a generalization of the finite difference method, created for arbitrary irregular grids \cite{LISZKA198083, Chung1981}. Despite the fact that this method requires a much lower mathematical background, it entails a more tedious implementation as the discretization coefficients are obtained from high-order Taylor expansions over non-uniform neighborhoods \cite{ONATE1996}. With the Radial Basis Function Finite Difference (RBF-FD \cite{fornberg_flyer_2015}), these coefficients are computed at each point by interpolating the solution over a small neighborhood using a radial basis enriched with polynomial functions. The RBF ensures good conditioning while the polynomial enrichment grants convergence. The radial basis function-based differential quadrature approach (RBF-DQ \cite{SHU20052001}) proposes a similar strategy and has been employed to simulate compressible flows \cite{SHU20052001} and heat conduction \cite{SOLEIMANI20101411} problems.

 Despite offering significant advantages (\eg{flexibility in handling complex geometries, avoiding computational overhead associated with mesh generation}), meshless methods possess considerable limitations. One notable limitation is that most meshless methods must use background cells to integrate a weak form over the entire problem domain. The requirement of background cells for integration makes the method not ``truly” meshless. Another challenge lies in enforcing boundary conditions accurately, especially so, as the degree of freedom may not be adapted to the problem geometry. As pointed out by Atluri and Zhu, in their work on meshless local Petrov-Galerkin methods, ensuring the satisfaction of boundary conditions in meshless methods can be non-trivial, potentially leading to inaccuracies in the solution; they attempt to mitigate this by implementing a penalty function \cite{Atluri1998}. Additionally, the lack of a structured mesh results in reduced convergence rates and increased computational costs; this occurs particularly in problems with high gradients or regions of varying solution behavior \cite{BELYTSCHKO19963}. These limitations underscore the importance of carefully assessing the suitability of meshless methods for specific numerical analysis tasks and considering alternative approaches where appropriate.

In amalgamating the desired advantages of meshless methods: accessibility, simplicity, h-P convergence, and low complexity, while simultaneously relinquishing the need for rigorous theoretical validation, we present the Least Squares Discretization method (LSQD). This approach only assumes that around each point, the solution can be decomposed on a local basis. No weak formulation of the PDE, nor numerical analysis theory is needed to discretize the problem. Instead, the PDE and boundary conditions are point-wise evaluated, and the local expansions are connected through pointwise continuity conditions. All of these conditions are expressed as linear equations involving evaluations of the known local basis functions and the unknown decomposition coefficients. They are regrouped in a rectangular linear system, from which the vector of unknown coefficients is defined as the least squares solution. The method is therefore entirely meshfree, requires minimal mathematical and programming background, and scales with the number of grid points.

Our computational exploration demonstrates how this method can be easily implemented to approximate solutions of elliptic PDEs completed with arbitrary boundary conditions enforced on non-trivial geometries, and using advanced data structures. We observe the solution converging both as the spatial resolution and the maximum polynomial order increase (h-P convergence), on both uniform and arbitrary highly non-graded tree structures. Because our method is meshless in nature, the underlying set of points can chosen as virtually anything. In particular, the points do not have to conform to the problem geometry or be aligned in any specific pattern. This facilitates the construction and manipulation of the data structure and offers unbounded solution adaptivity and error control strategies.

This manuscript is organized as follows. In section \ref{sec:A}, we present the LSQD method in the general context and discuss its practical implementation. The method is illustrated with a simple one-dimensional example. In section \ref{sec:highD} we discuss the extension of the method to higher dimensions and propose an implementation on non-graded quadtree grids. Convergence results are provided in section \ref{sec:Results} for a variety of grid and boundary configurations. We conclude in section \ref{sec:Disc}.

%%%%%%%%%%%%%%%%%%%%%%%%%%%%%%%%%%%%%%%%%%%%%%%%%%
%%%%%%%          METHODS/MATH              %%%%%%%
%%%%%%%%%%%%%%%%%%%%%%%%%%%%%%%%%%%%%%%%%%%%%%%%%%

\section{The LSQD Method}\label{sec:A}

\begin{figure}[h!]
\begin{center}
\includegraphics[width=0.9\textwidth]{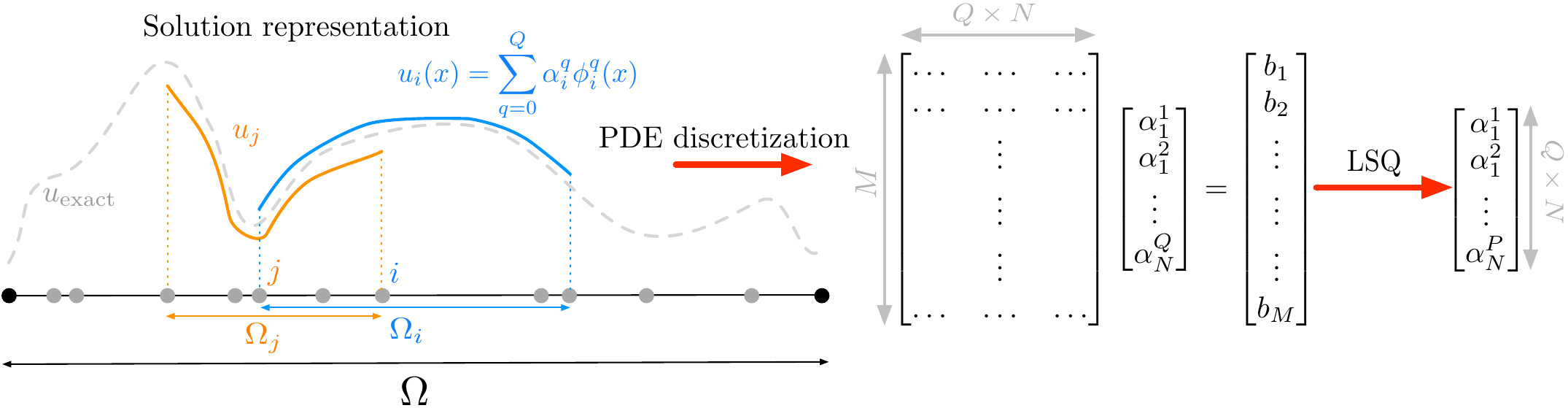}
%\vspace*{-0.5cm}
\caption{We represent the solution of a differential equation as a set of local expansions, denoted by $u_i(x)$, which are centered at arbitrary positions $x_{i=1..N}$. We build an overdetermined linear system by evaluating the differential equation around each neighborhood and enforcing continuity conditions between neighboring local representations. The coefficients $\alpha_i^q$ are the least squares solution to the overdetermined system.} 
\label{fig:1Dconvergence}
\end{center}
\end{figure}

\subsection{General Presentation}\label{sec:genrep}

For the LSQD method, we begin by considering the following elliptic problem:
\begin{eqnarray}
a u -\mu\triangle u = f(\mathbf{x}) \qquad \forall \mathbf{x} \in \Omega\subseteq\mathbb{R}^d, \quad f(\mathbf{x})\in\mathcal{C}^0(\Omega)
\label{eq:elliptic}
\end{eqnarray}
with generalized (e.g. Dirichlet, Neumann, Robin, etc.)\footnote{For Dirichlet boundary conditions: $\gamma=0$. For Neumann boundary conditions: $\beta=0$. For Robin boundary conditions: $\beta\neq 0$ and $\gamma \neq0$.} boundary conditions on the contour $\partial \Omega$, i.e.
\begin{eqnarray}
\beta u(\mathbf{x}) + \gamma \nabla u(\mathbf{x})\cdot \mathbf{n} = g(\mathbf{x}),  \qquad \forall \mathbf{x} \in \partial \Omega, \quad g(\mathbf{x})\in\mathcal{C}^0(\Omega),
\label{eq:bc}
\end{eqnarray} 
where $\mathbf{n}$ is the outward point normal to $\partial\Omega$. The parameters of the elliptic problem, $a$ and $\mu$, are positive and the solution is assumed to be continuously differentiable. We are interested in constructing a discrete approximation of the solution $u(\mathbf{x})$, using a set of points $(\mathbf{x}_1, ...,\mathbf{x}_N) \in \Omega$. 
To each point $\mathbf{x}_i$, we associate a set of neighbors $V_i$ and define the vicinity $\Omega_i$ as the convex hull of this set. For each vicinity $\Omega_i$, we approximate the solution using an expansion over a local smooth basis $(\phi^1_i(\mathbf{x}),...,\phi^Q_i(\mathbf{x}))$
\begin{eqnarray}
u(\mathbf{x})\approx \sum_{q=1}^Q\alpha_i^q\phi^q_i(\mathbf{x}), \qquad \forall \mathbf{x}\in \Omega_i.
\label{eq:exp_i}
\end{eqnarray}
Thus, by estimating the coefficients $\alpha_i^p$, at each $\Omega_i$, we hope to construct a family of local approximations to the solution $u(\mathbf{x})$.

%%%%%%%%%%%%%%%%%%%%%%%%%%%%%%%%%%%%%%%%%%%%%%%%%%
%%%%%%%       ALPHA CALCULATION            %%%%%%%
%%%%%%%%%%%%%%%%%%%%%%%%%%%%%%%%%%%%%%%%%%%%%%%%%%

To find the coefficients $\alpha_i^p$ we evaluate the PDE \eqref{eq:elliptic}, the continuity conditions  (for $\mathbf{u}$ and its first order derivatives), and the boundary condition \eqref{eq:bc} for the above local representation. For each node $\mathbf{x}_i$, we consider each neighbor pair $(\mathbf{x}_i,\mathbf{x}_j)$ with $\mathbf{x}_j\in V_i$ and construct the linear equations for $\alpha_i^p$ using the following procedure with three main steps:

%\maxime{the indices are wrong i->j} \maxime{for the boundary condition it is xj but a point on the interface between xi and xj}
\begin{enumerate} \label{LSQD_steps}
\item {\bf Partial Differential Equation:} we choose a point $\mathbf{x}^{PDE}_{ij}\in \Omega_i$ at which we evaluate the PDE for the expansion \eqref{eq:exp_i}. 

\begin{eqnarray}
a\sum_{q=1}^Q\alpha_i^q\phi^q_i(\mathbf{x}^{PDE}_{ij}) - \mu \sum_{q=1}^Q\alpha_i^q\triangle\phi^q_i(\mathbf{x}^{PDE}_{ij}) = f(\mathbf{x}^{PDE}_{ij}).
\end{eqnarray}
\item {\bf Continuity conditions:} we choose a point $\mathbf{x}^{C}_{ij}\in \Omega_i$ at which we enforce that the expansions centered at $\mathbf{x}_i$ and $\mathbf{x}_j$ are identical
\begin{eqnarray}
\sum_{q=1}^Q\alpha_i^q\phi^q_i(\mathbf{x}^{C}_{ij}) = \sum_{q=1}^Q\alpha_j^q\phi^q_j(\mathbf{x}^{C}_{ij}),
\end{eqnarray}
and that the expansions for all first-order derivatives $\partial_d$ are also identical at that point
\begin{eqnarray}
\sum_{q=1}^Q\alpha_i^q\partial_d\phi^q_i(\mathbf{x}^{C}_{ij}) = \sum_{q=1}^Q\alpha_j^q\partial_d\phi^q_j(\mathbf{x}^{C}_{ij}).
\end{eqnarray}
\item {\bf Boundary conditions:} if $\mathbf{x_i}$ is close to the contour $\partial \Omega$, we find a point $\mathbf{x}^{BC}_{i} \in \partial \Omega$, at which we enforce the boundary condition for the local expansion \eqref{eq:exp_i}. 

\begin{eqnarray}
\beta \sum_{q=1}^Q\alpha_i^q\phi^q_i(\mathbf{x}^{BC}_{i}) + \gamma \sum_{q=1}^Q\alpha_i^q\nabla \phi^q_i(\mathbf{x}^{BC}_{i})\cdot \mathbf{n} = g(\mathbf{x}^{BC}_{i}).
\end{eqnarray}
\end{enumerate}
Grouping all these equations, we form a rectangular system of width $ Q \times N$ and height $M$ (the total number of equations). 
\begin{eqnarray}
 \begin{bmatrix}
  \hdots \quad \hdots\quad \hdots  \\
  \hdots \quad \hdots\quad \hdots  \\
  \vdots \\
  \vdots \\
\vdots \\
  \hdots \quad \hdots\quad \hdots
\end{bmatrix}\begin{bmatrix}
  \alpha_1^1  \\
  \alpha_1^2  \\
  \vdots \\
  \alpha_N^Q
\end{bmatrix} &= \begin{bmatrix}
 b_{1}  \\
  b_{2}  \\
  \vdots \\
\vdots \\
\vdots \\
  b_{M}
\end{bmatrix} \qquad \Leftrightarrow \qquad A\alpha=b.
\label{eq:rect_system}
\end{eqnarray}
If the matrix $A$ is full rank, we construct the coefficient vector $\alpha$ as the least squares solution:
\begin{eqnarray}
\alpha = \underset{x}{\mathrm{argmin}}||Ax-b||_2 \qquad \Leftrightarrow \qquad A^TA\alpha=A^Tb.
\label{eq:lsqr_system1}
\end{eqnarray}
If the matrix is not full-rank, we can consider the stabilized problem
\begin{eqnarray}
 \left(A^TA+ \epsilon I\right)\alpha=A^Tb,
\label{eq:lsqr_system1}
\end{eqnarray}
which is well-known to be invertible for any $\epsilon > 0$. Nevertheless, we must be cognizant of two extremities: choosing a large $\epsilon$ will affect the accuracy of the method, but choosing its value to be too small may leave the system singular to machine precision. We address this topic in section \ref{sec:stab}.
%%%%%%%%%%%%%%%%%%%%%%%%%%%%%%%%%%%%%%%%%%%%%%%%%%
%%%%%%%          NEIGHBORS                 %%%%%%%
%%%%%%%%%%%%%%%%%%%%%%%%%%%%%%%%%%%%%%%%%%%%%%%%%%
\subsection{Practical Implementation }

\subsubsection{Basis Selection }

Formally, the method does not impose any restriction on the basis selection, yet, it is well known that a naive choice could lead to singular least squares systems, making it impossible to construct the solution (without stabilization). For example, allow us to consider the standard interpolation problem (\ie choosing $\alpha =1$ and $\mu=0 $ in \eqref{eq:elliptic}, with appropriate boundary conditions); from the Mairhuber–Curtis theorem \cite{Mairhuber1956}, we know that choosing global non-centered basis functions may result in a singular interpolation matrix. For the rest of this study, at each point $\mathbf{x}_i$ we will use polynomial functions, centered at the point $\mathbf{x}_i$ and rescaled by the size of the neighborhood $V_i$. In doing so, we aim to limit singularities in the final least squares system, as well as provide a framework for achieving high-order convergence. We will denote by $P$ the maximum polynomial order of the basis, and continue to denote by $Q$, the total number of elements in the basis.

\subsubsection{Neighborhood Construction}\label{sec:neighb}

From the method description, we see that it is important to construct a full-rank matrix A, while trying to keep it as compact and sparse as possible for computational efficiency. The matrix structure and size are primarily dictated by the neighborhoods configurations. It is therefore imperative to understand how many neighbors are needed in each neighborhood, how to construct them, and what properties they should satisfy. In particular, we identify three key properties the neighborhood must satisfy

\begin{enumerate}[label=\roman*]
    \item {\bf System Size - } The total number of equations must be greater or equal to the total number of basis functions (\ie the number of unknown coefficients). 
    \label{cond:systsize}
    \item{\bf Connectivity - } The neighbors' adjacency matrix must be invertible, or in other words the neighboring graph must be irreducible. In practice, it means that traveling along neighboring connection, you should be able to travel from any point $\mathbf{x}_i$ to any other point $\mathbf{x}_j$.   \label{cond:connectivity}
    \item{\bf Independence - } The neighborhood should contain points with at least $P + 1$ unique coordinates in each directions. 
    \label{cond:indep}
\end{enumerate}

The first condition is necessary for the matrix $A$ to be full rank. The second condition ensures that all the neighborhoods are connected, and therefore, that the computation of each basis's coefficients depends on all other coefficients. If this condition is not satisfied, there will be clusters of nodes and, consequently, basis coefficients that are separated from each other, violating the global nature of the elliptic equation \eqref{eq:elliptic}. The third condition reinforces linear independence between the coefficient equations. It is constructed from the observation that if the problem is invariant in all directions except one, the matrix A must still be full rank. In that case, the number of linearly independent basis functions is simply $P+1$ and so we need at least as many unique coordinates in the non-invariant direction.

We point out that the neighborhoods do not have to be symmetric, meaning that if $\mathbf{x}_j$ is a neighbor of $\mathbf{x}_i$ (\ie $\mathbf{x}_j \in V_i$), $\mathbf{x}_i$ does not have to be a neighbor of $\mathbf{x}_j$. They do not all have to be the same size, nor respect specific geometric configurations. In practice though, it is natural to keep them local and compact.

\subsubsection{Evaluation Points}

As mentioned in \ref{sec:genrep}, the evaluation points $\mathbf{x}^{PDE}_{ij}$ and $\mathbf{x}^{C}_{ij}$ for the PDE and continuity conditions must be chosen in the domain $\Omega_i$ where the local $i^{th}$ expansion is formally defined. Typically, the evaluation points are selected to be on the line segment $[\mathbf{x}_i,\mathbf{x}_j]$ as this ensures locality. We note that choosing different evaluation points for the PDE and continuity conditions ($\mathbf{x}^{PDE}_{ij}\neq \mathbf{x}^{C}_{ij}$) may reinforce linear independence. 

The boundary conditions are only evaluated when the point $\mathbf{x}_{i}$ is close to the interface $\partial \Omega$. As well as other method parameters, {\it proximity of the interface} must be defined in each context. For example, in 1D, we can place boundary points on the interface and evaluate the boundary condition at these points only. In higher dimension, provided that we have a smooth implicit representation (\eg a level-set representation), we can always detect whether $\Omega_i$ crosses $\partial \Omega$, and if so, find at least one point $\mathbf{x}^{BC}_{i}$ on $\Omega_i \cap \partial \Omega$.

\subsubsection{System Construction \& Inversion}
The main advantage of our method is that the matrix $A$ is straightforward to construct. Once the neighborhoods have been constructed, we need to loop over all neighboring pairs and evaluate the PDE, continuity, and boundary conditions. The matrix assembly only requires evaluation of the basis function at the evaluation points. No quadrature rule, finite difference formula, or any other numerical analysis algorithms are required.

Assuming that the Least Squares problem is non-singular, the solution vector can be computed using a preconditioned conjugate gradient method (since $A^TA$ is symmetric), or using the LSQR method \cite{10.1145/355984.355989,Saunders1995} to avoid forming $A^TA$. Forming $A^TA$ may allow us to construct a wider range of preconditioners, and so it is, {\it a priori}, unclear which of these options will be computationally more efficient. It should be noted that the matrices $A$, $A^T$, and $A^TA$ are sparse and can be constructed in $\mathcal{O}(NQ^2)$ and $\mathcal{O}(NQ^3)$ respectively. Therefore, leaving aside the neighborhood construction, the total computational cost will scale as the complexity of the linear solver, typically $\mathcal{O}(NQ^3\sqrt{\kappa})$, with $\kappa$ being the condition number of $A^TA$.

%%%%%%%%%%%%%%%%%%%%%%%%%%%%%%%%%%%%%%%%%%%%%%%%%%
%%%%%%%            1D STUFF                %%%%%%%
%%%%%%%%%%%%%%%%%%%%%%%%%%%%%%%%%%%%%%%%%%%%%%%%%%

\subsection{One-Dimensional Illustration}

We consider the following Poisson problem 
\begin{eqnarray}
-\triangle u &=& f(x) \qquad \forall x \in \Omega = [0,1],
\label{eq:1Delliptic} \\
\nonumber u(x) &=& g(x)  \qquad \forall x \in \partial \Omega .
\end{eqnarray}
In order to study the performance of the method we choose the exact solution to be
\begin{eqnarray*}
    u_{\textrm{exact}}(x)=\sin(7x)+\cos(7x)
\end{eqnarray*}
and we set $f(x)$ and $g(x)$ accordingly. All local basis function sets are chosen to be centered and rescaled polynomials, \ie for each $x_i$ and over each corresponding neighborhood $V_i$, we define
\begin{eqnarray}
\phi_i^q(x_i)= \left(\frac{x-x_i}{|V_i|}\right)^{q-1}, \qquad|V_i|= \max_{j\in V_i}|x_i-x|, \qquad x\in V_i, \quad q=\{1,..,Q = P+1\}.
\end{eqnarray}

Each neighborhood $V_i$ is first created by choosing the node $\mathbf{x}_i$ itself, its closest left and right neighbors. This ensures that the neighborhoods overlap with at least another neighborhood in each direction and therefore that the connectivity condition \ref{cond:connectivity} is satisfied. Additionally, because the cardinality of the basis set is equal to $P+1$, we need at least $P+1$ different neighbors in $V_i$ to satisfy both the system size \ref{cond:systsize} and independence \ref{cond:indep} conditions.  Therefore, if $P > 2$, we add the remaining $P-2$ closest points to the neighborhood. As a result, each neighborhood will contain $\eta = \max(3,P+1)$ elements. For simplicity, we choose the evaluation points to be $\mathbf{x}_{ij}^{PDE}=\mathbf{x}_{ij}^C=\mathbf{x}_j$.

\begin{figure}[h!]
\begin{center}
\includegraphics[width=1\textwidth]{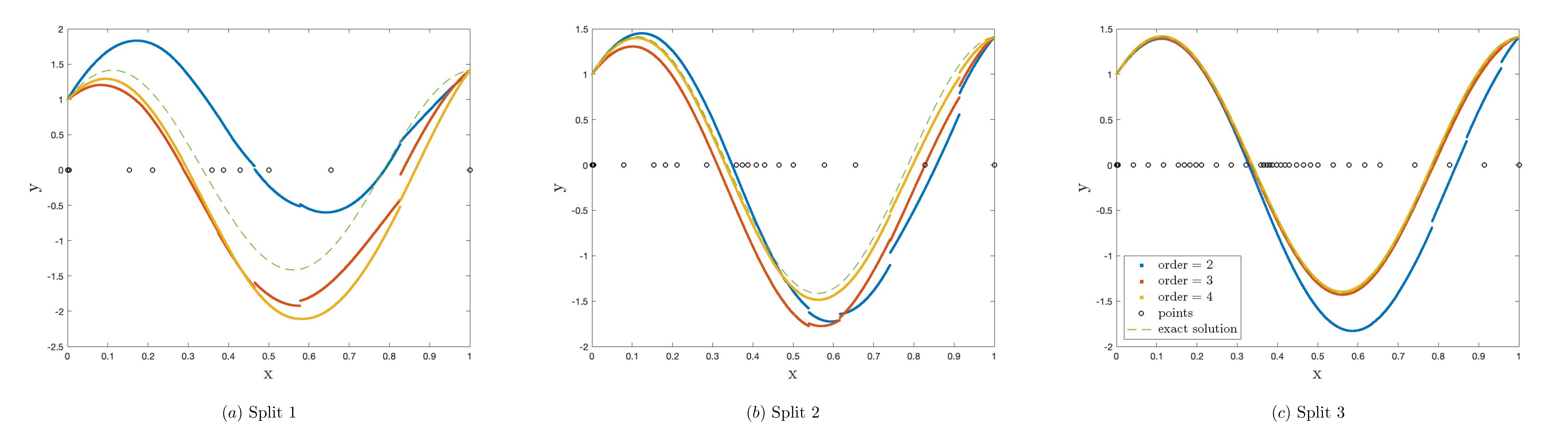}
\caption{1D Method - solution representation for increasing spatial and polynomial resolution.  We start from a randomly generated set of points $\mathbf{x}_{i=1..10}$ and then recursively split the grid by adding one grid point exactly between each existing pair of direct neighbors. The exact solution is depicted in dashed lines.}
\label{fig:1Dsolution}
\includegraphics[width=1\textwidth]{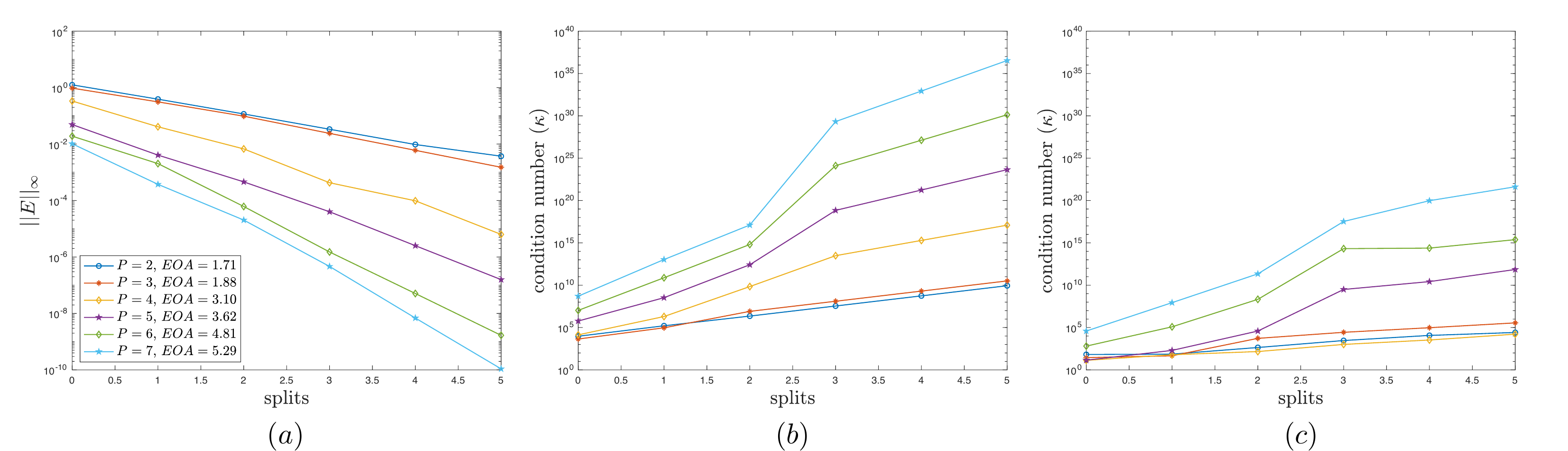}
\caption{1D Method - (a) h-P convergence analysis for Poisson's equation with Dirichlet boundary conditions. $L^\infty$ error and corresponding Estimated Order of Accuracy (EOA). (b) Condition number of the least squares matrix $A^TA$. (c) Condition number of the preconditioned least squares matrix, using an incomplete Cholesky decomposition.}
\label{fig:1Dconvergence}
\end{center}
\end{figure}

Once again, we emphasize the fact that advanced numerical analysis techniques are not required to construct our system, nor to obtain our solution (e.g. no Taylor expansions, convoluted integrals, or quadrature rules). Our implementation contains around 100 lines of MATLAB code and the least squares system is solved using the simple backslash solve command (i.e. $\alpha = A^TA \backslash A^Tb$). The backslash command for general square matrices uses LU decomposition; however, for symmetric matrices (which is our case) the Cholesky decomposition is used instead. The results can be seen in Figure \ref{fig:1Dsolution}. 

To study the h-P convergence of the method, we initialize the points $x_{i=1..10}$ at random and then recursively add midpoints between adjacent neighbors to double the local spatial resolution. We will refer to splits as the number of subdivisions performed from the original points set, and point out that the local resolution scales like $h^{-\text{splits}}$, where $h$ is the original local spatial resolution. In parallel, we vary the polynomial order using $P=2, 3, 4$. The error is computed by comparing the exact solution and the local approximation at each point. We will focus on the maximum error only. 

Convergence results are compiled in Figure \ref{fig:1Dconvergence}, where we see the error diminishing both as the resolution and the polynomial order increase. On the other hand, the condition number inflates with both the resolution and maximum polynomial order; this phenomenon should be seen as a proxy for the computational time needed to solve the least squares problem. This is as one might expect and is the exchange we are willing to make in the name of simplicity for this method.

With the system being symmetric, we can precondition it using a standard incomplete Cholesky decomposition, which, for this example, drastically reduces the condition number. Typically, the condition number of the preconditioned system is the square root of the unpreconditioned one. As we will see in section \ref{sec:errorest} the discontinuities between adjacent expansions can be used to estimate.

%%%%%%%%%%%%%%%%%%%%%%%%%%%%%%%%%%%%%%%%%%%%%%%%%%
%%%%%%%            2D STUFF                %%%%%%%
%%%%%%%%%%%%%%%%%%%%%%%%%%%%%%%%%%%%%%%%%%%%%%%%%%

\section{Extension to Higher Dimensions \& Implementation on Quadtree Grids}\label{sec:highD}
\begin{figure}[h!]
\begin{center}
\includegraphics[width=0.4\textwidth]{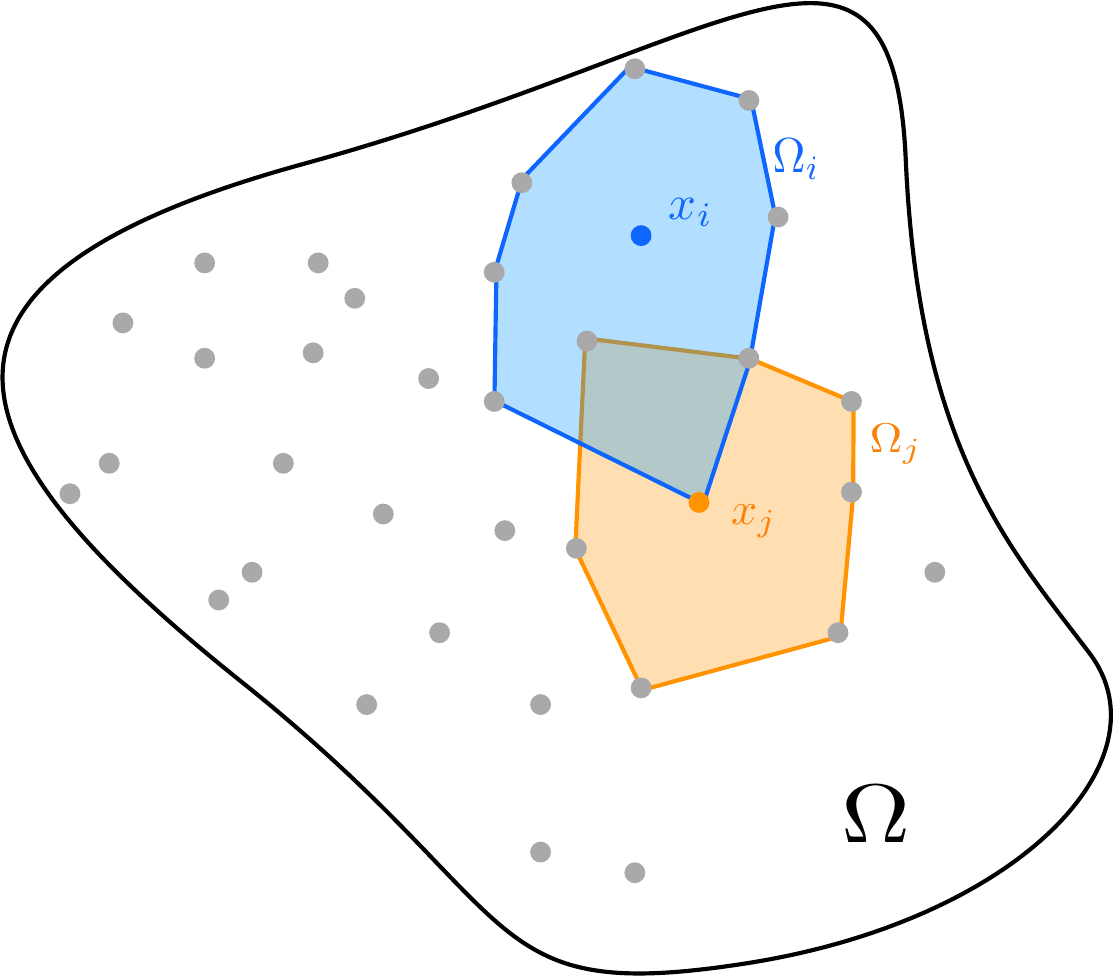}
\caption{Extension to higher dimensions: domain of interest $\Omega$, and local vicinities $\Omega_i, \Omega_j$. The neighborhoods $V_i, V_j$ are defined as the set of points contained in the corresponding vicinities.} 
\label{fig:2drep}
\end{center}
\end{figure}

Despite the one-dimensional example serving as an obligatory proof of concept in terms of convergence, the full capabilities of our method, such as dealing with highly adaptive data structure and arbitrary geometries and boundary conditions, only manifest themselves in higher dimensions. 

For our method, the extension to higher dimensions is straightforward from both the theoretical and practical perspectives. As illustrated in Figure \ref{fig:2drep}, the domain of interest $\Omega$ is still discretized by a set of points $\mathbf{x}_{i=1..N}$. The solution is then decomposed over a local basis $\phi_i^1(\mathbf{x})..\phi_i^Q(\mathbf{x})$ around the set of points, and the coefficients for these expansions are obtained by evaluating the partial differential equations, boundary, and continuity conditions over local neighborhoods, as underscored in section \ref{sec:genrep}. The construction of the linear system remains unambiguous and the properties of the final system are unchanged ($A^TA$ is symmetric, positive definite).

The main difference lies in the neighborhood construction: both in their shape and size. In 1D, finding close neighbors from a pre-sorted list is trivial. In higher dimensions, for an efficient construction, the pre-sorting must be done using more complex data structures, such as sorting trees or hash tables. 

Also, because the derivative continuity condition will be evaluated in each direction, the total number of equation written per neighboring pair will differ, and thus the minimal neighborhood size. Namely, in $d$-dimensions, for a neighborhood $V_i$ composed of $\eta$ neighbors (including $\mathbf{x}_i$), we will evaluate
\begin{itemize}
    \item[-] $\eta$ PDE equations,
    \item[-] $\eta - 1$ continuity conditions,
    \item[-] $d(\eta - 1)$ first-order derivative continuity conditions ($\eta - 1$ for each partial derivative),
\end{itemize}
so a total of $\eta\left(d+2\right)-d-1$ equations. Since we want the matrix $A$ to be full rank, we will systematically pick the number of basis of local basis function $Q$ such that
\begin{eqnarray}
    Q \leq \eta\left(d+2\right)-d-1. \label{eq:ngbcondition}
\end{eqnarray}
Note that this strategy does not account for the boundary condition evaluations. This is an intentional decision to derive a robust and simple strategy: the total number of boundary condition evaluations is hard to predict {\it a priori}, and these conditions will only increase the system size and rank.

Since sorting structures will be needed to construct neighborhoods, we will locate the points at the cell centers of a Quadtree to limit the number of data structures and facilitate the presentation and implementation of our method.
 Quadtree grids \cite{Samet1988AnOO}, a paradigmatic data structure for high-performance computations (see \eg \cite{KUCHEROVA2021110591,BLOMQUIST2024112695,D3SM01196H, HEYDARI2022110755}). The points $\mathbf{x}_i$ will be located at the tree's cell centers. We will represent the domain $\Omega$ through a level-set function and only approximate the solution at the points located inside $\Omega$ (see Figure \ref{fig:neighb_cons}).

\subsection{2D Basis Functions} 
\begin{figure}[h]
\begin{center}
\includegraphics[width=0.95\textwidth]{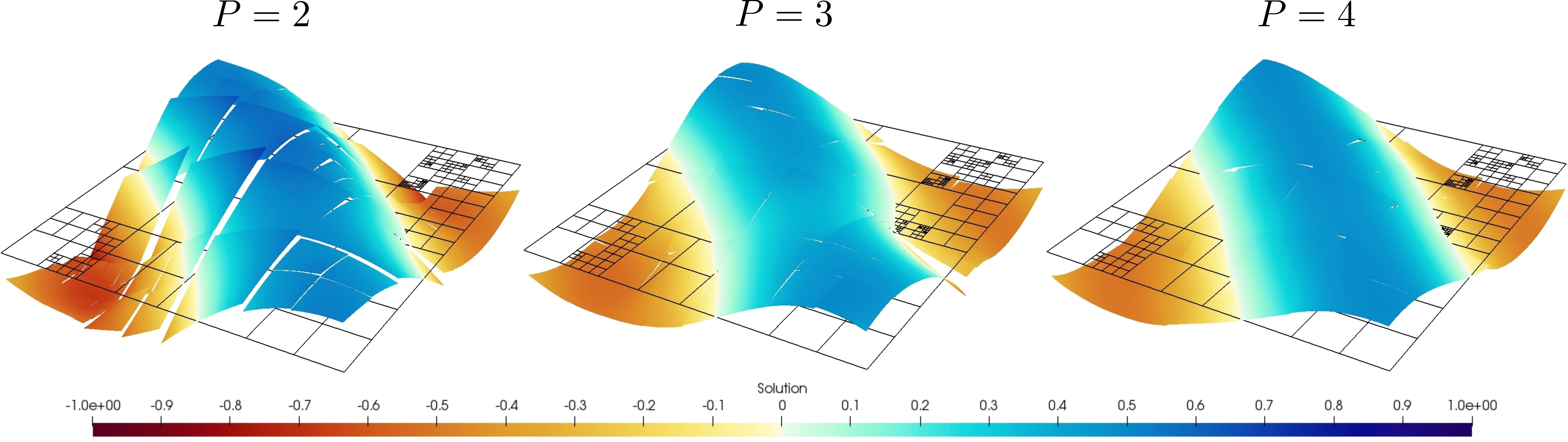}
\caption{LSQD Approximation of the function $u_{\textrm{exact}}=\cos(x+y)$ for increasing polynomial order on an arbitrary Quadtree grid.   }
\label{fig:quad_ex}
\end{center}
\end{figure}
We chose our basis functions as products of one-dimensional rescaled polynomials. Namely around $\mathbf{x}_i=(x_i,y_i)$:
\begin{eqnarray}
\phi_i^q\left(\mathbf{x}\right) = \left(\frac{x-x_i}
{|V^x_i|}\right)^{p_x}\left(\frac{y-y_i}{|V^y_i|}\right)^{p_y} 
\qquad \forall \mathbf{x} = (x,y) \in\Omega_i,
 \end{eqnarray}
where the rescaling coefficients are
\begin{eqnarray}
\qquad |V_i^x| = \max_{\mathbf{x}\in V_i}|x_i-x|, \qquad |V_i^x| = \max_{\mathbf{y}\in V_i}|y_i-y|.   
\end{eqnarray}
The rescaling helps keep the coefficients in $A$ -- therefore of $A^TA$ -- of similar magnitude, consequently keeping the condition number of $A^TA$ low. In this context, this is essential as we are dealing with very disparate length scales and potentially high-order polynomial functions.

We consider the set of polynomial functions of degree $P$, meaning that $p_x+p_y\leq P$ (see Figure \ref{fig:quad_ex}). The total number of basis functions will be
\begin{eqnarray}
    Q = \frac{(P+1)(P+2)}{2},
\end{eqnarray}
and so conditions \eqref{eq:ngbcondition} dictates that we will need to choose the neighborhood size $\eta$ such that
\begin{eqnarray}
      \eta \geq \frac{(P+1)(P+2)+6}{8}. \label{eq:ineqbeta}
\end{eqnarray}
Consider, for example, using a polynomial order $P=7$, from the above condition we see that $\eta\geq 9$. The neighborhood will grow with the number of basis functions and, thus, local unknown coefficients. However, this growth will occur slowly, implying that the neighborhood will remain very compact; this attribute is highly desirable, especially in higher dimensions, as keeping ``things” as local as possible implies keeping $A$ and $A^TA$ as sparse as possible. The above inequality \eqref{eq:ineqbeta} is not enough to specify $\eta$ as the connectivity \eqref{cond:connectivity} and independence \eqref{cond:indep} conditions must also be satisfied.

\subsection{Tree-Based Neighborhood Construction} 
% \label{2D_neigh}

\begin{figure}[h!]
\begin{center}
    \includegraphics[width=1\textwidth]{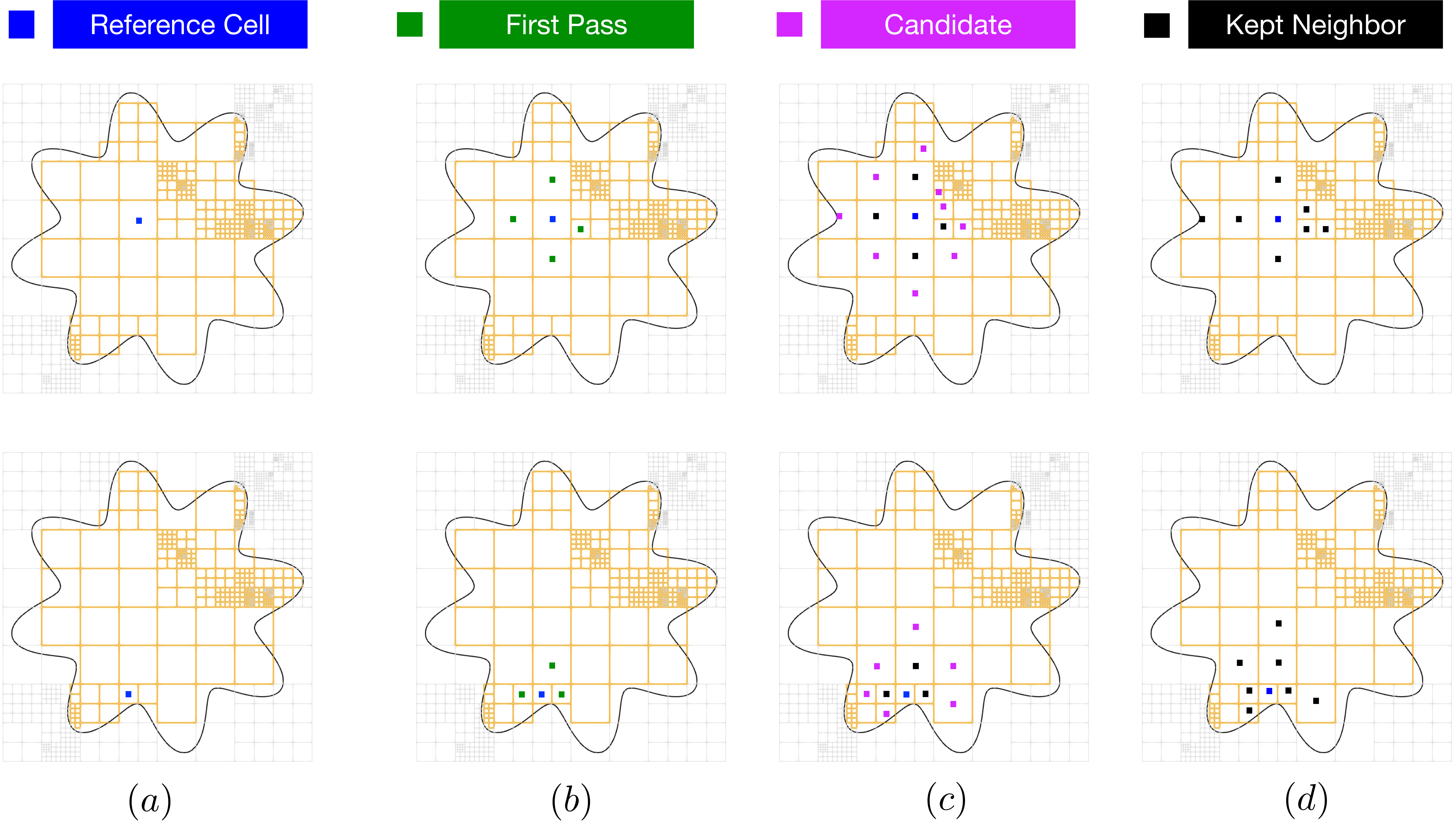}
   \caption{Illustration of neighborhood construction around a reference cell (in blue) for polynomial basis order $P=4$.  The process is carried out one layer at a time until enough neighbors are identified to satisfy the system size and independence conditions. When the reference cell lays next to the interface, the neighbors are chosen within $\Omega$ (bottom row). The solution is only constructed at the cell centers of the orange cells located inside $\Omega$. The outside grey cells are only used to locate the interface while evaluating the boundary condition.}
    \label{fig:neighb_cons}
\end{center}
\end{figure}

Figure \ref{fig:neighb_cons} illustrates the construction of the neighborhood $V_i$ for any given point $\mathbf{x}_i$. We start by adding the point itself (blue) to the neighborhood. Next, the algorithm recursively
picks a direct neighbor in each direction for each point in the neighborhood and store them in a candidate set (pink). In this set, the candidates are sorted in counterclockwise order starting to the right. In case one of the candidates lies outside $\Omega$ (see bottom row of Figure \ref{fig:neighb_cons}), we move it to the ghost set $G_i$ that will be used to find evaluation points for the boundary condition. 

On the first pass, we add the entire candidate set to the neighborhood to ensure that the connectivity condition \eqref{cond:connectivity} is met. On the subsequent passes, while the other two conditions are not met, we go through the candidates set and add a candidate $\mathbf{c}$ to the neighborhood if either

\begin{enumerate}[(a)] 
    \item The system size condition \eqref{cond:systsize}, reformulated as inequality \eqref{eq:ineqbeta}, is not met.
    \item The independence condition \eqref{cond:indep} is not met in at least a direction $d$ and the coordinate $c_d$ of the candidate in that direction is not already present in the neighborhood, \ie $\forall \mathbf{v}\in V_i \quad v_d\neq c_d$.
\end{enumerate}

We point out that using non-uniform point distributions reduces repetition in the coordinates, which eases the fulfillment of condition \eqref{cond:indep} and ultimately generates smaller and more compact neighborhoods. If the grid is refined enough, the algorithm mentioned above will always converge, which allows us to systematically construct neighborhoods. However, in case the grid is under-resolved, there may not be enough neighbors available, which can lead to infinite loops. To prevent this, we set a maximum number of passes (usually $P+4$). In these degenerated situations, the grid is inconsistent with the polynomial order and geometry, so the problem is completely under-resolved, and the method will produce an inaccurate solution. In the next section, we will explain how we use the neighborhood to evaluate the problem and form a linear system.

\subsection{Evaluation points}
For simplicity, we keep the continuity and pde evaluation points identical ($
    \mathbf{x}^{PDE}_{ij} =  \mathbf{x}^{C}_{ij} = \mathbf{x}_{j}$). For the evaluation of the boundary condition, we go through the set of ghost neighbors $G_i$ obtained through the neighborhood construction, and for each ghost $\mathbf{g}$ neighbor we find a point $\mathbf{x}^{BC}_{ig} \in [\mathbf{x}_i,\mathbf{g}]$ lying on the interface and evaluate the boundary condition there. Since we are using an implicit representation for the interface (level-set), finding $\mathbf{x}^{BC}_{ig}$ boils down to a zero-finding problem on the segment $[\mathbf{x}_i,\mathbf{g}]$ which is easily solved using, for example, the bisection method.

If $\mathbf{x}_i$ is located at the center of a cell at the edge of the computational domain, the boundary condition is evaluated on the corresponding edges.

\subsection{Matrix Rescaling}
In addition to the basis function rescaling, we rescale the rows of the matrix $A$ by their respective norms. This operation should be seen as a re-weighting of each linear equation, \ie multiplying the entire system by a diagonal matrix with non-zero diagonal entries. In doing so, we keep the magnitude of the coefficients in $A$ as close to each other as possible to ensure that the final condition number is as low as possible. Specifically, we start by computing the weighting coefficients $w_k$ as the inverse of the row norms of $A = (a_{kl})$
\begin{eqnarray}
     w_k = \left(\sum_{l=1}^{Q\times N}|a_{kl}|\right)^{-1}\qquad \forall k=1,2,...M 
\end{eqnarray}
and then rescale the entire system
\begin{eqnarray}
     a_{kl} &=& {a_{kl}}{w_k} \qquad \forall k=1, 2,...M, \quad l=1,2,..Q\times N\\
     b_{k} &=& {b_{k}}{w_k} \qquad \forall k=1, 2, ...M.
\end{eqnarray}
Recall $M$ is the total number of equations in the overdetermined system. 

\subsection{Stabilization of the LSQ System} \label{sec:stab}
As discussed earlier, the final least squares system is guaranteed to be invertible if a small diagonal perturbation $\epsilon$ is added. This addition may not be necessary, especially as it essentially penalizes large coefficients, thus affecting the accuracy of the method. In practice, this stabilization should be seen as a means to keep the condition number of the system bounded. 

Adding the diagonal perturbation effectively shifts the spectrum of the system. We know that the eigenvalues $\lambda_k$ of $A^TA$ are positive, and so the condition number of $A^TA$ is the ratio of the largest and smallest eigenvalues $\lambda_{\max}$ and $\lambda_{\min}$ 
\begin{eqnarray*}
    \kappa(A^TA) = \frac{\lambda_{\max}}{\lambda_{\min}}.
\end{eqnarray*}
If $\lambda_k$ is an eigenvalue of $A^TA$,  $\lambda_k+\epsilon$ is an eigenvalue of the stabilized problem $A^TA+\epsilon I$, and so the condition number of the stabilized problem is
\begin{eqnarray*}
    \kappa(A^TA+\epsilon I) = \frac{\lambda_{\max}+\epsilon }{\lambda_{\min} + \epsilon }.
\end{eqnarray*}
If we set a practical maximum condition number $\kappa_{M}$ for the stabilized system, and seek the minimal perturbation needed $\epsilon_{\min}$ to remain below this limit
\begin{eqnarray*}
    \kappa_{\max} = \frac{\lambda_{\max}+\epsilon_{\min}}{\lambda_{\min} + \epsilon_{\min}},
\end{eqnarray*}
we find that
\begin{eqnarray*}
    \epsilon_{\min} = \frac{\lambda_{\max}-\kappa_{\max}\lambda_{\min}}{\kappa_{\max} -1}.
\end{eqnarray*}
For the current implementation we set $\kappa_{\max}=10^{40}$, and estimate the largest and smallest eigenvalues of $A^TA$ using Gershgorin circles theorem
\begin{eqnarray*}
    \lambda_{\max} &=& \max_{k=1..M} c_{kk} + \sum_{l\neq k}|c_{kl}|, \\
     \lambda_{\min} &=& \max\left(0,\min_{k=1..M} c_{kk} + \sum_{l\neq k}|c_{kl}| \right),
\end{eqnarray*}
where the $c_{kl}$ are the coefficients of the matrix $A^TA$. In practice, we measured $\lambda_{\max}\approx10^2$, $\lambda_{\min}=0$, leading to $\epsilon\approx 10^{-38}$.
We note that $\kappa_{\max}$ is the maximum allowed condition number for the preconditioned system. In practice, we preconditioned the system with an incomplete Cholesky precondition, and so the final condition number will be much lower. Based on our one-dimensional observations, we expect the preconditioned condition number to be bounded above by $\approx 10^{20}$. Again, this is an upper limit that should be seen as the absolute worst-case scenario. For all cases presented in the following section, we find the least squares system to be invertible and construct the numerical solution using an iterative Conjugate Gradient method, typically taking $10^2-10^4$ iterations to converge to residuals $||(A^TA+\epsilon_{\min} I)\alpha-A^Tb||_{\infty}$ smaller than $10^{-26}$.

\begin{figure}[t!]
\begin{center}
        \includegraphics[width=0.8\textwidth]{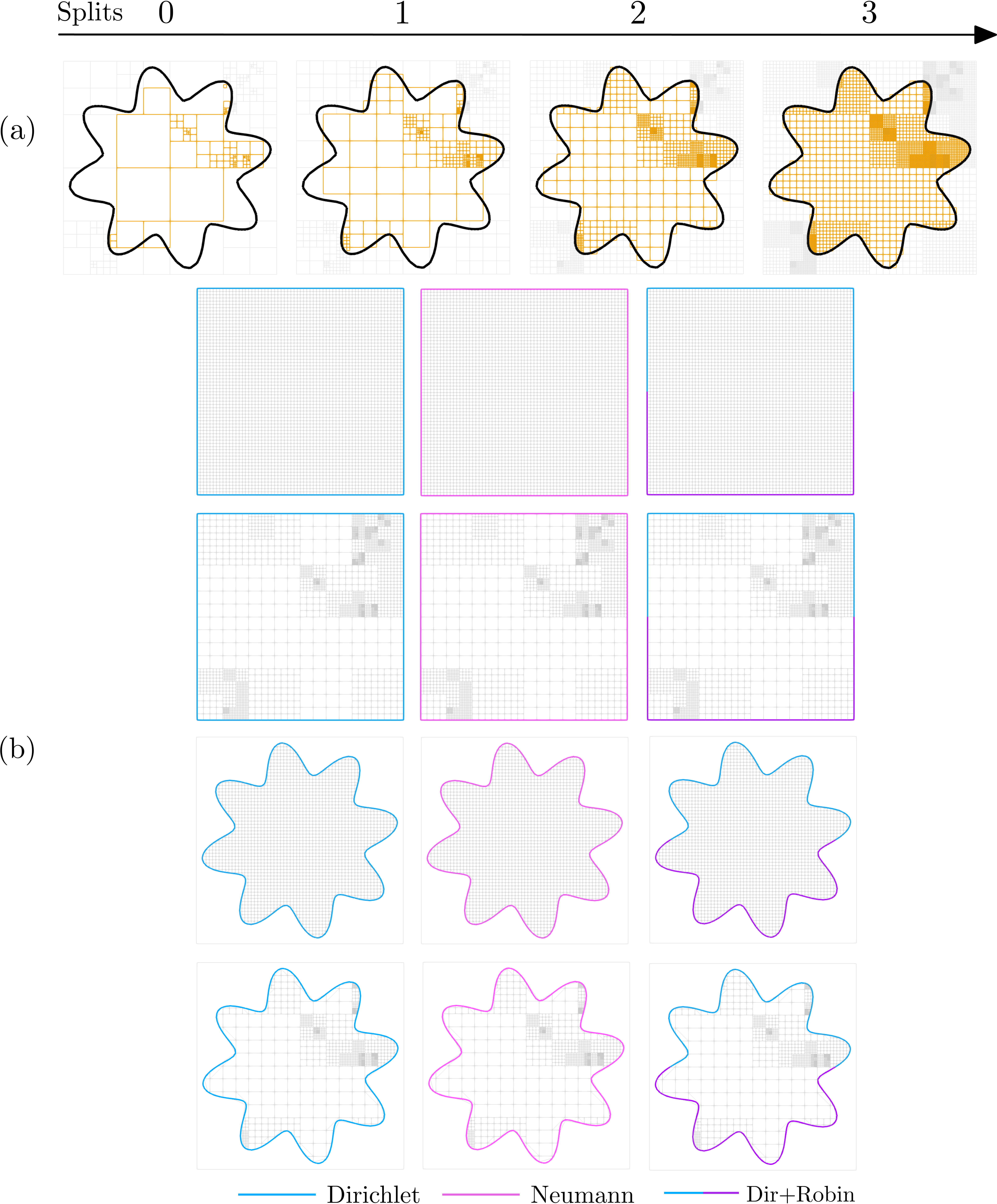}
    \caption{%\selma{if we feel like it, change the colors to (current colors A51EE9, 3FBEFF, FF55FF,) and change all other figures ..}
    Configurations for the convergence analysis: (a) Illustration of arbitrary tree structures. The grid is randomly generated (first grid) which we call Split 0, then recursively refined (left to right - Split 1 to Split 3). 
    (b) Grids and domains used -  For each boundary condition type (Dirichlet, Neumann, Mixed Dirichlet Robin), enforced on a square (top rows) or the octofoil domain (bottom rows), we will consider both uniform (row 1 and 3) and adaptive grids (rows 2 and 4). This representation will be used to label all convergence figures hereafter}.
    \label{fig:BC}
\end{center}
\end{figure}

%%%%%%%%%%%%%%%%%%%%%%%%%%%%%%%%%%%%%%%%%%%%%%%%%%
%%%%%%%           2D -- RESULTS              %%%%%%%
%%%%%%%%%%%%%%%%%%%%%%%%%%%%%%%%%%%%%%%%%%%%%%%%%%
\setlength\dashlinedash{0.2pt}
\setlength\dashlinegap{1.5pt}
\setlength\arrayrulewidth{0.3pt}
\begin{table*}[] \centering
\caption{Convergence summary - Average estimated order of convergence for various boundary conditions, configurations, and polynomial orders $P=2,..,5$.}
\label{tab:EOC}
%\ra{1.3}
\begin{small}
\begin{tabular}{@{}lrrrrrr@{}}\toprule
 & \textbf{Dirichlet} & \textbf{Neumann} & \textbf{Dir.+Robin} & \textbf{Vary $\frac{\alpha}{\mu}$} & \textbf{Degenerate} & \textbf{Overall Avg} \\ \midrule
$\mathbf{P=2}$         & 1.86 & 2.17 & 2.02 & 1.63 & 1.46 & 1.83   \\ \hdashline
$\mathbf{P=3}$         & 1.96 & 2.20 & 2.07 & 2.16 & 1.91 & 2.06   \\ \hdashline
$\mathbf{P=4}$         & 3.45 & 3.30 & 3.46 & 3.51 & 3.45 & 3.43   \\ \hdashline
$\mathbf{P=5}$         & 3.66 & 4.53 & 3.89 & 4.46 & 3.60 & 4.03   \\ %\hdashline
% \textbf{Total Avg}   &  & 2.94 & 2.86 & 2.61 &   \\
\bottomrule
\end{tabular}
\end{small}
\end{table*}

\section{Results}\label{sec:Results}
Presented here are convergence results for various domain types, levels of mesh refinement, boundary condition types, and numerous formulations of equation \eqref{eq:elliptic} (by changing the coefficients). Furthermore, we present built-in error estimate. Finally, we discuss the implications of Neumann boundary conditions with a computational time analysis.

\subsection{Problem Setup}
We consider the proposed elliptic problem (\ref{eq:elliptic}) with varying $\alpha , \mu$ values over two domains. The first domain is the unit square $\Omega = [-1,1] \times [-1,1]$. The second domain is the interior of an octofoil centered at $(0,0)$ and defined by the zero contour of the following level set equation (see Figure \ref{fig:BC})
\begin{eqnarray}
    \phi(x,y) = -\Biggl[ \frac{4}{5} +\frac{4}{25}\sin\biggl(8\arctan\left(\frac{y}{x}\right)\biggr)\Biggr] + \sqrt{x^2+y^2}.
\end{eqnarray}
For both of these domains, we use the boundary condition formulation \eqref{eq:bc} and adjust the weights $\beta$ and $\gamma$ to achieve the desired condition types. We also consider two types of mesh refinement: uniform and randomly generated (see Figure \ref{fig:BC}). Each of the representations will be used for the convergence analysis in section \ref{sec:Results}.

We choose the exact solution to be $u_{\textrm{exact}}(x,y) = e^{x+y}$ and we set the forcing term $f(\mathbf{x})$ and the boundary condition $g(\mathbf{x})$ accordingly. For each of the representations, we vary the polynomial order of our basis from 2 to 5, and recursively split all grid cells to study the h-P convergence of the $L^{\infty}$-errors (see Figure \ref{fig:BC}(a)). \\
The local $L^1$ and global $L^{\infty}$-error are defined respectively as
\begin{eqnarray}
    ||e_u(\mathbf{x_i})||_1 = \left\lvert u_{\textrm{exact}}(\mathbf{x}_i) - \sum_{q=1}^Q\alpha_i^q\phi_i^q(\mathbf{x}_i)\right\rvert .
\end{eqnarray}
\begin{eqnarray}
    ||e_u||_{\infty} = \max_{i=1..N}\left\lvert u_{\textrm{exact}}(\mathbf{x}_i) - \sum_{q=1}^Q\alpha_i^q\phi_i^q(\mathbf{x}_i)\right\rvert .
\end{eqnarray}
The Estimated Order of Convergence (EOC) will be computed from the global error, using power fitting curves.

\begin{figure}[t]
\begin{center}
    \includegraphics[width=0.8\textwidth]{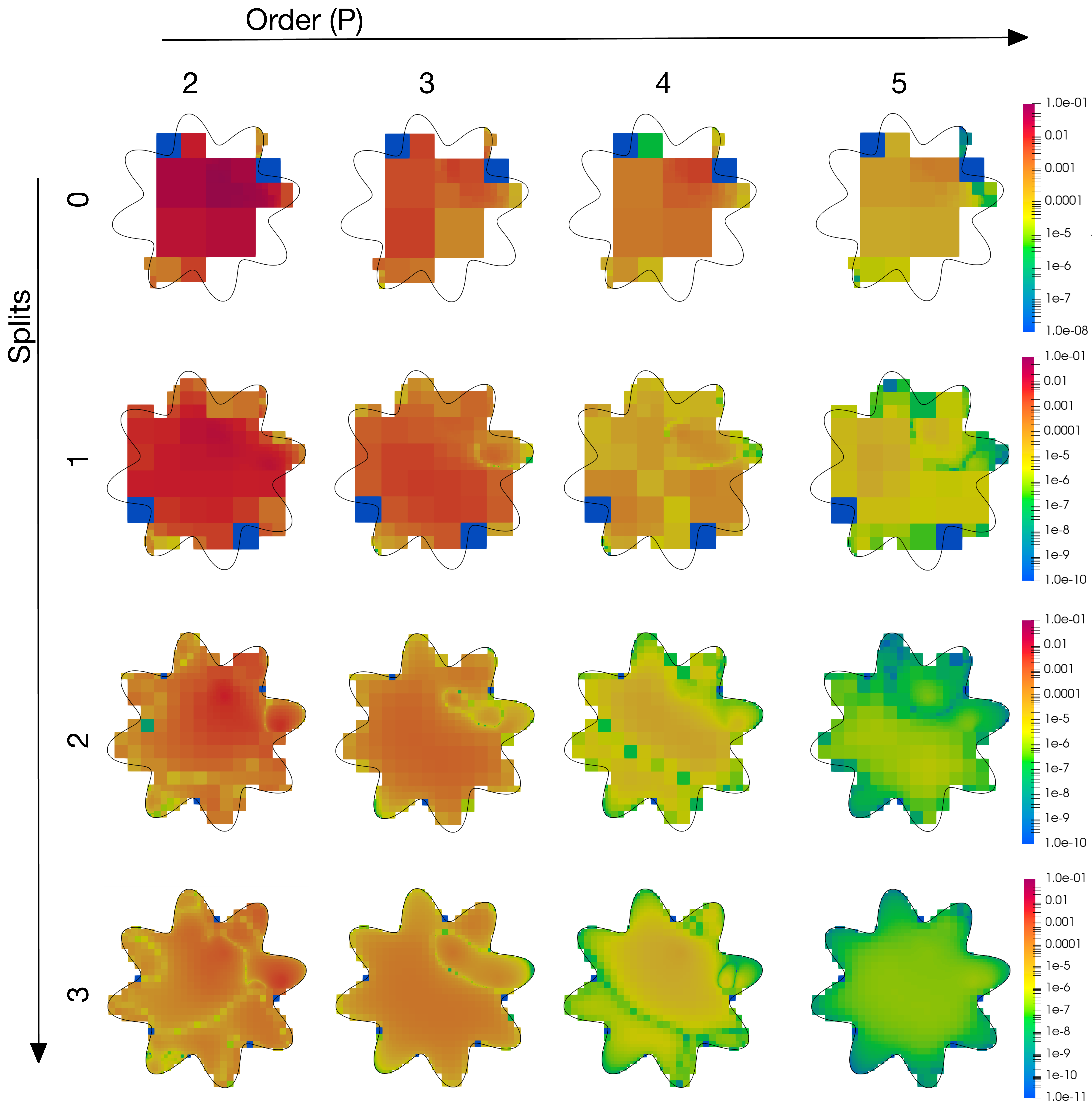}
    \caption{Convergence of the local error: Octofoil domain with Dirichlet boundary conditions and an adaptive grid. We see the error decaying to zero both as the polynomial order increases or the grid is refined.}
    \label{fig:err_con}
\end{center}
\end{figure}

\subsection{Classical Boundary Conditions: Dirichlet, Neumann, and Robin}
\begin{figure}[t!]
\begin{center}
    \includegraphics[width=0.85\textwidth]{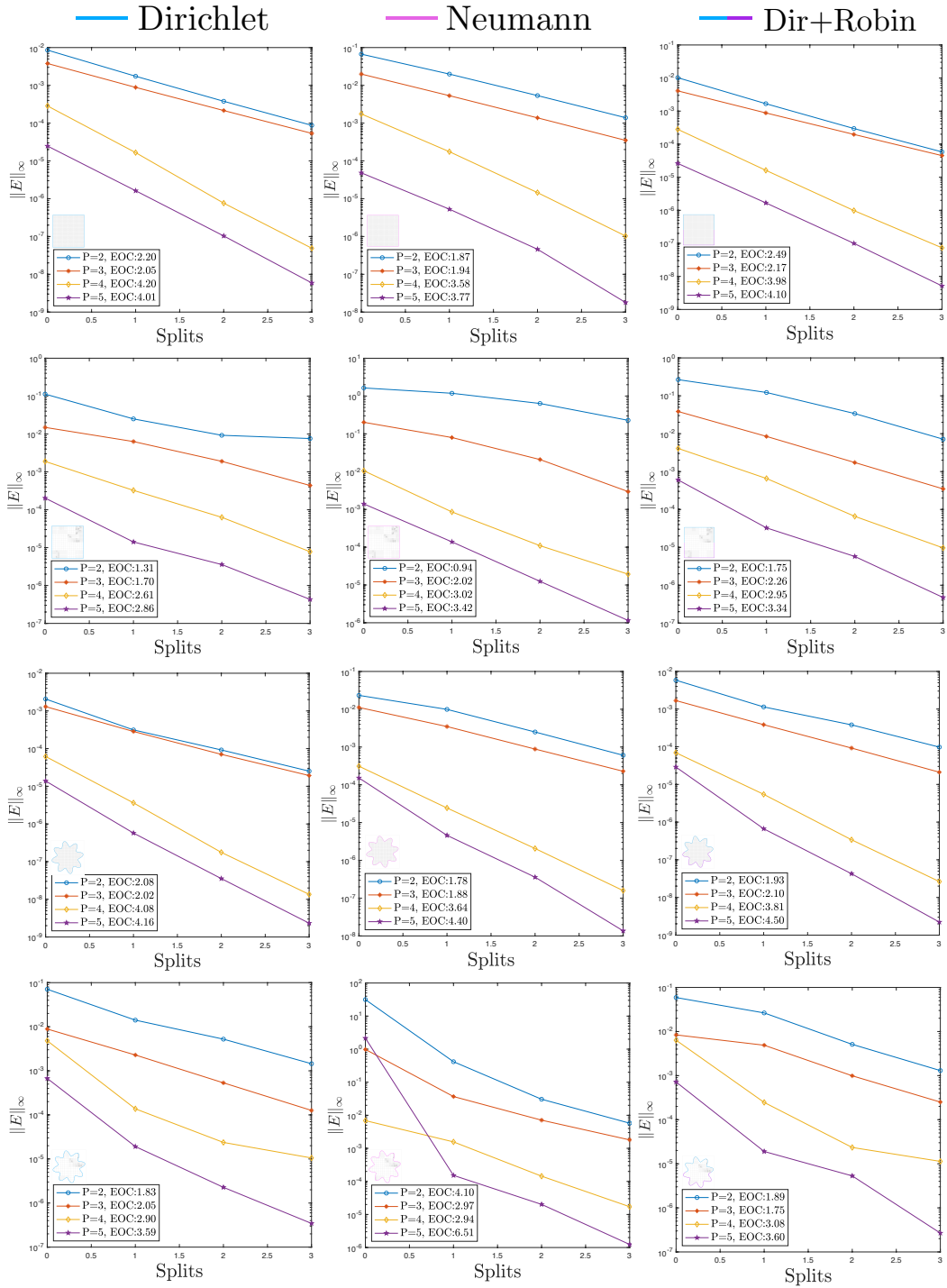}   
    \caption{Convergence of $L^\infty$ global error for several configurations. Columns denote boundary condition type. Rows correspond to (1) a uniform grid on a square domain, (2) adaptive square, (3) uniform octofoil, and (4) adaptive octofoil (top to bottom). See Figure \ref{fig:BC} for mesh examples.}
    \label{fig:error_all_cases}
\end{center}
\end{figure}

Consider $\alpha=\mu=1$ for all twelve configurations presented in Figure \ref{fig:BC}. The local error convergence for the octofoil domain with a Dirichlet Boundary condition ($\beta = 1$), using an adaptive refinement is depicted on Figure \ref{fig:err_con}. The h-P convergence of the global error for each set of boundary conditions is represented in Figure \ref{fig:error_all_cases}. For all twelve configurations, we see the error decaying to zero as the grid resolution and the polynomial order increase. The estimated order of convergence (see Table \ref{tab:EOC} and Figure \ref{fig:error_all_cases}) are fairly constant across all examples, ranging from approximately $2$ (P=2) to $4$ (P=5). 

We point out that the octofoil examples on the adaptive grid, and even more so with the mixed boundary condition case (labeled "Dir. + Robin" in Figure \ref{fig:BC}), are particularly challenging to deal with traditional discretization techniques. Again, non-trivial boundary conditions are enforced on complex geometries, and the grid is randomly generated and not conforming to the interface. In fact, for these last two reasons, the traditional Finite Element Method cannot be employed here, as they require conforming mesh, and when implemented on quadtree grids, requires graded tree (\ie the size ratio between neighboring cells cannot exceed 2). There are in theory no obstructions from using the Finite Differences of Volumes methods, but the constructions of the discrete differential operators and fluxes rapidly becomes a practical nightmare, especially in the high-order context. The modal Discontinuous Galerkin method allows for a systematic construction of the discrete system, where volume and contour integrals must be performed on the fractions of the computational cells and faces contained in the domain of interest, which can become very involved when the mesh is not adapted to the interface. Again, LSQD requires no convoluted discretization nor quadrature formula. Even in the presence of complex data structure or geometries, simple basis functions evaluations at pre-specified locations are enough to yield above second-order convergence.

\subsection{Degenerate Cases} \label{sec:degenerate}
Here, we investigate degenerate cases where the PDE or boundary conditions may lead to a non-coercive weak formulation and thus potentially ill-posed problems, non-positive or even non-invertible systems. The LSQD method obtains the polynomial coefficient from an SPD system and does not rely on the weak formulation. 
Similarly to the Least Square Finite Element Method (LSFEM \cite{JIANG199013, Jiang1998}), the LSQD method constructs the discrete system by evaluating the PDE and boundary conditions everywhere (formation of the matrix A), and then defining the Galerkin matrix taking the scalar product of these evaluations ($A^TA$). As a result, it is not affected by the potential ill-posedness of the weak form, and the resulting linear system is guaranteed to be invertible.

\begin{figure}[t]
\begin{center}
        \includegraphics[width=1\textwidth]{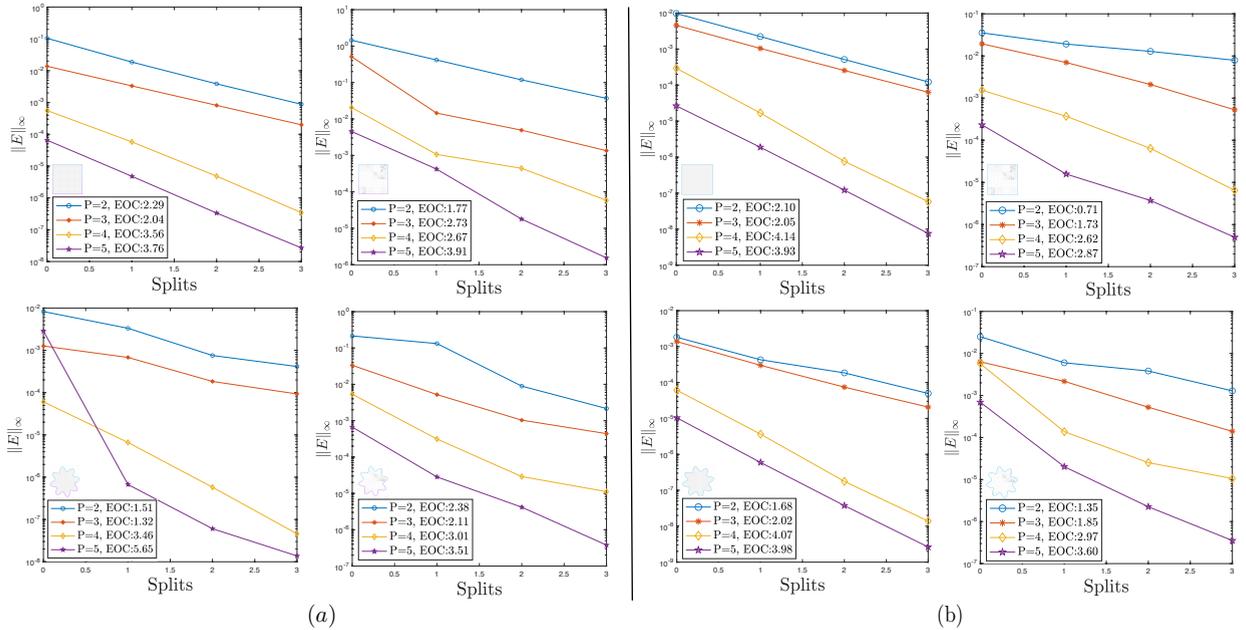}   
    \caption{Two instances of a degenerate case: (a) Mixed Robin and Dirichlet boundary conditions, (b) Dirichlet boundary condition for a highly negative value of $\mu$.}
        \label{fig:degenerate}
    \end{center}
\end{figure}

For the first case, we considered equation \ref{eq:elliptic} with $a=1, \mu = 1$ and imposed a split boundary conditions of Dirichlet and Robin. The degeneration comes from the Robin boundary condition, for which we set $\beta= 1, \gamma=-1$. For the second case, we keep $a = 1$, but set $\mu =-100$ (\ie Helmholtz's equation), and enforce a Dirichlet boundary condition. The convergence results are depicted in Figure \ref{fig:degenerate}. The method achieves h-P convergence for both of the degenerate cases and all grid configurations. The numerical error dips rapidly below reasonable error tolerances, and the EOC behaves normally, increasing with the polynomial order and split number. We note that results using exclusively traditional Robin boundary conditions ($\beta = 1$, $\gamma = 1$) on the entire domain achieved the similar results.

%%%%%%% NEUMANN BC'S %%%%%%%%%%%%%%

\subsection{Sensitivity to Problem Parameters} \label{sec:sensit}
Despite consistent patterns of convergence for various boundary conditions tested here (\ie{dirichlet, dirichlet + robin, degenerate, \it{etc.}}), Neumann boundary conditions seem to perform poorly -- especially on adaptive grids. In fact, in Figure \ref{fig:error_all_cases}, we observe errors of magnitude upwards of $10^{1} \text{ and even } 10^{2}$.

From (\ref{eq:elliptic}), we note that when $\mu$ is significantly smaller than $a$, the problem resembles an interpolation problem, which is relatively straightforward to solve. On the other hand, if $a$ is considerably smaller than $\mu$, the equation resembles a -- potentially -- ill-posed Poisson problem, \eg failing the solution uniqueness criterion. In fact, the well-posedness of the Poisson problem under Neumann boundary conditions hinges on factors such as the smoothness of the domain and boundary, the characteristics of the Neumann boundary data $g$, and how well the boundary conditions align with the given problem. When the Neumann boundary condition fails to offer sufficient information to determine $u$ uniquely, it can result in non-uniqueness or instability within the solution, rendering the problem ill-posed. Therefore, by setting $a=\mu = 1$, the well-posedness of the problem could have impacted the convergence/accuracy in the Neumann boundary condition setting.
To investigate this, we select four different values of $\frac{a}{\mu}$. The expectation is that as we increase the ratio, the problem should become easier to solve.

\begin{figure}[t!]
\begin{center}
        \includegraphics[width=0.9\textwidth]{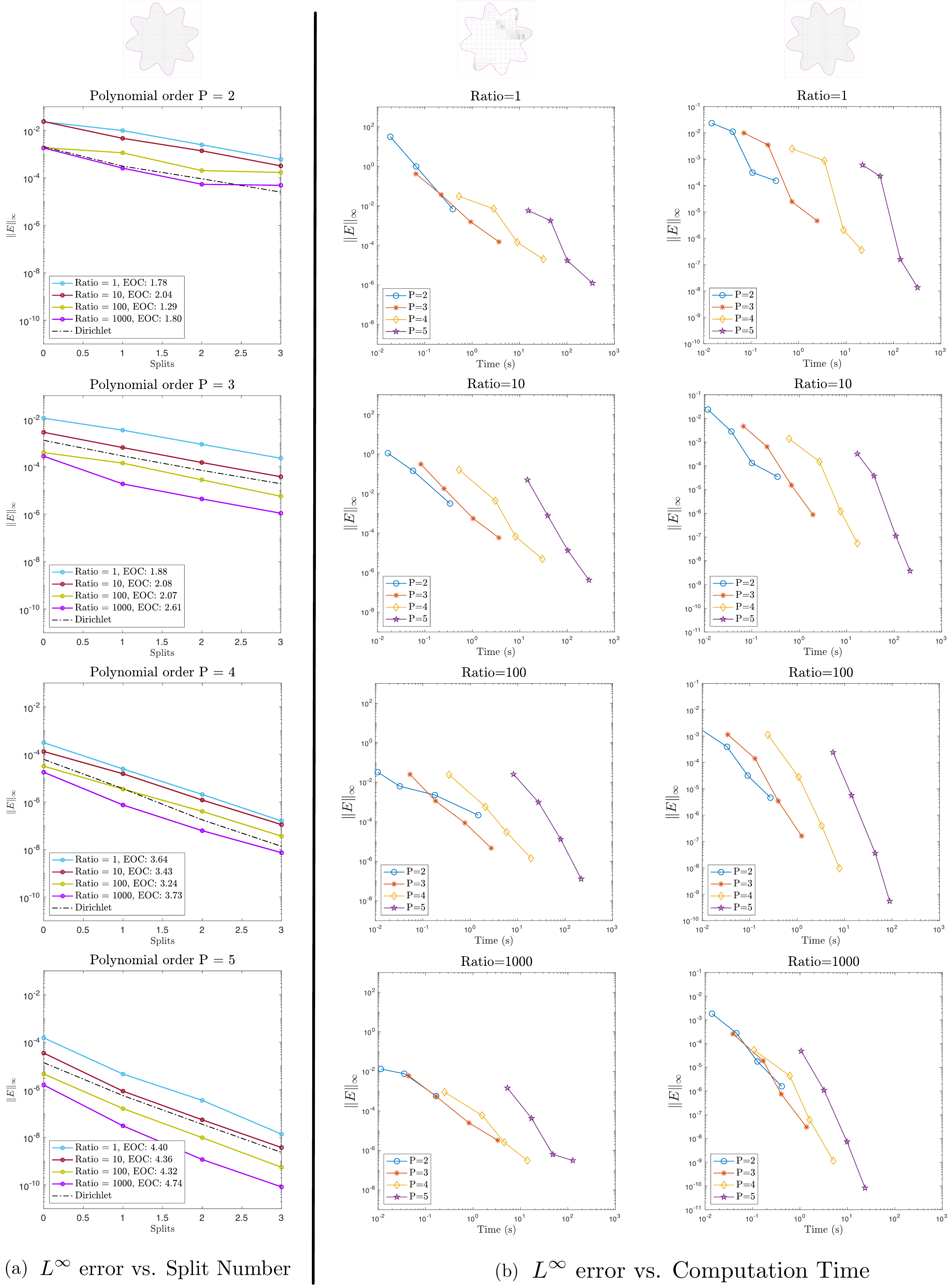}
        \caption{Neumann boundary condition on the octofoil domain: impact of $\frac{\alpha}{\mu}$ ratios. (a) Global error trend for increasing split number for a uniform grid for each polynomial basis. (b) Impact of computational time on global error for both adaptive and uniform grid types.  }
    \label{fig:ratio}
\end{center}
\end{figure}

Presented in Figure \ref{fig:ratio} (a) are convergence results for Neumann boundary conditions imposed on the octofoil interface with a uniform grid for $\frac{a}{\mu} = 1,10,100,1000$. We observe that for each polynomial basis order, $P$, the global error maintains a similar trend among the four tested ratios $\frac{a}{\mu}$, respectively. This behavior corroborates the congruity of the Estimated Order of Convergence (EOC) between various $\frac{a}{\mu}$ given order $P$. Clearly, the method performs as expected with the Neumann constraint across the various ratios; nonetheless, it is noteworthy that larger values of $\frac{a}{\mu}$ allow for a significant decrease in the global error ($\approx 2$  orders of magnitude). In fact, when $\frac{a}{\mu} \geq 100$, the method demonstrates accuracy comparable to, or surpassing, that which was achieved with Dirichlet boundary conditions. Hence, with Neumann boundary conditions, smaller values of $\frac{a}{\mu}$ can compromise the effectiveness of the method.

\subsection{Computational Time}\label{sec:cputime}

Previously, it was established that lower $\frac{\alpha}{\mu}$ ratio values associated with a Neumann boundary condition create a more challenging computational problem (\ref{sec:sensit}). For this reason, we showcase the performance and efficiency of the LSQD method by exploring a variety of $\frac{\alpha}{\mu}$ ratios. Figure \ref{fig:ratio} (b) illustrates the maximum error obtained by our solver in comparison to the time it took to run in serial on a Dell precision 3240. It is apparent that, in general, achieving a lower maximum error requires a longer runtime and higher polynomial order. However, we also see that as the ratio $\frac{\alpha}{\mu}$ decreases, not only does it take more time to solve the system, but the order of accuracy declines. 

We also note that the slope of the error versus computational times increases with the polynomial order, suggesting that in the asymptotic limit of highly resolved simulations (small error or large CPU time), high-polynomial computations will always be the most efficient. Though, for the presented results this asymptotic regime is never reached. Therefore, in practice, it is imperative to decide which variable (minimizing error or computational cost) takes precedence in a given setting. Presuming that the modeled system has a tight tolerance for error, while computational resources are not a primary concern, the conspicuous objective is achieving a low global error. Figure \ref{fig:ratio}(b) hints at the necessary setup for attaining specific accuracy levels; the higher the polynomial basis order, the higher the accuracy. However, if the priority shifts towards saving computational time and/or memory or simply capturing the general behavior of the solution, opting for a moderately compact setup and a lower-order polynomial basis might suffice.

%%%%%%%%%%%%%% Error Estimator %%%%%%%%%%%%%%

\subsection{Error Estimator}\label{sec:errorest}

We propose an alternative way to compute local error estimate as: 
\begin{eqnarray}
    ||e_u(\mathbf{\widetilde{x}})||_{Est} = \max \left\lvert \sum_{q=1}^Q\alpha_{l_1}^q\phi_{l_1}^q(\mathbf{\widetilde{x}}) - \sum_{q=1}^Q\alpha_{l_2}^q\phi_{l_2}^q(\mathbf{\widetilde{x}})\right\rvert, \quad l_1,l_2 \in L = (\text{neighboring cells}).
\end{eqnarray}

Obtaining analytical solutions to PDEs is only possible in the simplest cases, so it should be expected that the exact solution to the PDE \eqref{eq:elliptic} may not be available. This error estimator is built into the algorithm, as evaluating accuracy predicated on traditional methods tends to be infeasible. For a comparison between the theoretical prediction and practical results, we restrict our analysis to the Dirichlet case. 

\begin{figure}[t]
\begin{center}
    \includegraphics[width=0.8\textwidth]{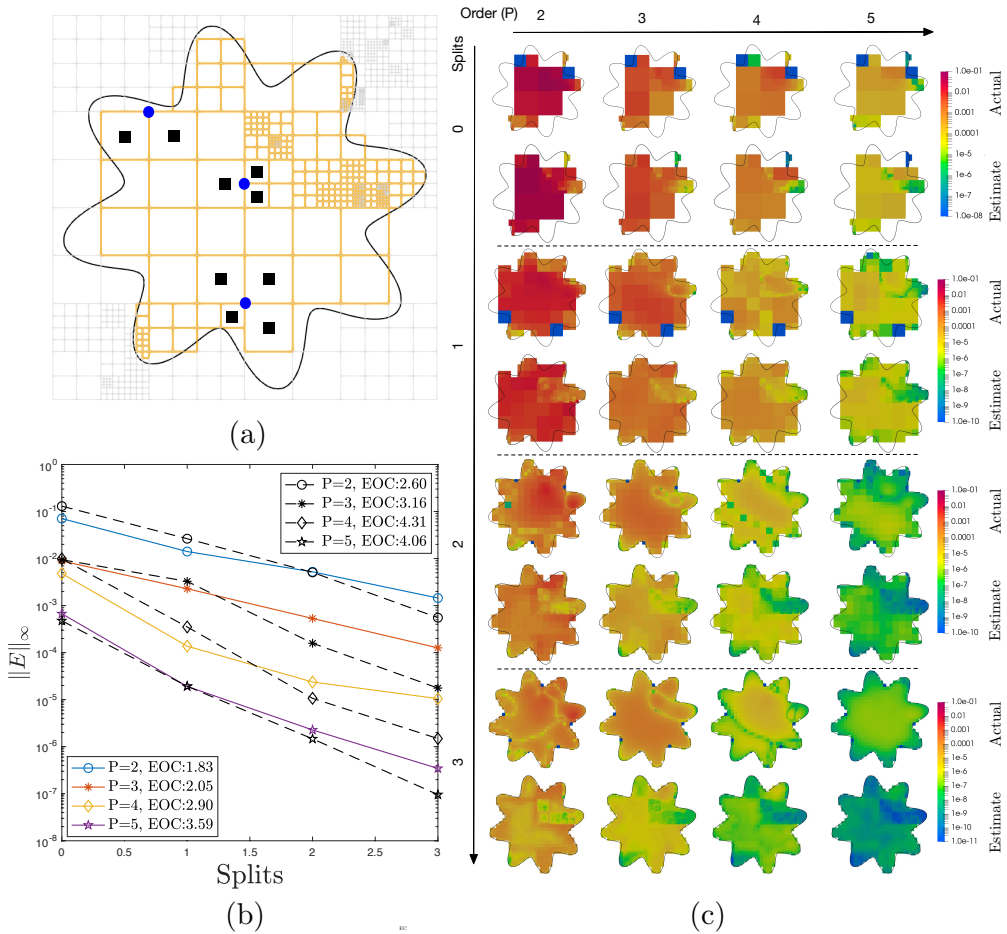}
    \caption{Error Estimator Analysis on adaptive grids: (a) Three distinct node locations and what is considered their respective surrounding/neighboring cells. (b) Comparison between computational order of accuracy (colored lines) and predicted order of accuracy (dashed lines). (c) Illustration of computational and predicted error for increasing polynomial order and grid refinement.}
    \label{diri_both_err}    
\end{center}
\end{figure}

The function of our node-based error estimator is to compare local approximations from neighboring cells. All corners of the cells located inside the domain become node locations $\mathbf{\widetilde{x}_i}$. For each node, the touching cells are used to calculate the error. In one dimension, we can picture the error as the discontinuity between two local expansions around that node (see Figure \ref{fig:1Dsolution}). To visualize this process in two dimensional space, refer once again to Figure \ref{diri_both_err}(a). If we want to compute the estimate error around the reference node (in blue), we use the local approximation functions constructed for the closest neighbor cells (black), and evaluate at the considered reference node $\mathbf{\widetilde{x}_i}$ (note that the maximum amount of neighboring cells is four: left, right, top, and bottom). Finally, the error estimate is computed as the maximum difference between the local values obtained.

Focusing, in particular, on Figure \ref{diri_both_err}(b-c), we are able to see a comparison between simulation error and estimated error for multiple polynomial orders and for increasing grid refinement. These results are encapsulated, numerically, in \ref{diri_both_err}(b), where calculations (on legend) generalize that the predicted EOC is consistently higher than the actual EOC. However, when taking a closer look at the distribution of the colored (simulated) and dashed lines (predicted), the predicted global error is not free of variance, meaning that neither the simulated, nor the predicted global error is consistently superior. It becomes evident (see fig. \ref{diri_both_err}(c)) that the main differences between our actual and estimated error are transpiring at the boundary of our octofoil domain. Nonetheless, the most stark differences dissipate as the grid is refined (this is referred to as `Splits'). Overall, we find that this cheap built-in estimator offers meaningful global error and convergence estimates.

\section{Discussion / Conclusion}\label{sec:Disc}

In this paper, we proposed a new paradigm for solving elliptic PDEs, one that would retain all of the advantages of meshless methods, while remaining accessible in theory and in practice. Least Squares Discretization Method (LSQD) does not require any quadarature rules, weak formulations, body-fitted grids, nor does it necessitate advanced mathematical theory. Nevertheless, it constructs h-P converging numerical solutions on non-trivial data structures with reasonable computational complexity. LSQD method easily handles complex geometries and is virtually impervious to ill-posed problems.

Although we began with a one-dimensional proof of concept, we demonstrated the feasibility in extending our method to a higher dimensional space, and its applicability to a broad range of problems, boundary conditions, and computational parameters (\eg grid resolution, polynomial orders). For all considered cases, we consistently observe h-P convergence, with Estimated Order of Convergence (EOC) ranging from 1.83 to 4.53 (see Table \ref{tab:EOC}), and in practice the global error falls rapidly below any reasonable error tolerance. Moreover, when investigating degenerate cases (\ie where the problem is potentially ill-posed as a consequence of the PDE coefficients or boundary conditions), the method observes kindred results (see \ref{sec:degenerate}), exemplifying the robustness of the LSQD method. 

Finally, we proposed an {\it a posteriori} node-based error estimator to unveil the local solution agreement. From Fig \ref{diri_both_err} (c), we see that as the polynomial order increases and the grid is refined, the estimated and measured/experimental errors display the same convergence pattern and similar spatial localization. This estimator cannot be used as a lower bound for actual error as seen in \ref{diri_both_err}(c), but provides a general blueprint for where error is to be expected and, thus, can be used for optimal mesh refinement. Coupled with a spatial adaptivity capability, this error estimator could be used as an error controller. 

From our computational time analysis (see \ref{sec:cputime}), we observed steeper efficiency curves for higher polynomial orders, suggesting that using high-order representations might be advantageous, at least for very large systems. In practice, for small test problems like the ones considered here, these benefits may not manifest, and low polynomial orders seem to be the best strategy to get an accurate solution rapidly. In the examples presented here, computational times peaked around $10^3$ seconds on a single core, or $\approx 17$ minutes for the largest systems ($\approx 1.3\times 10^6$ degrees of freedom for 4 splits and P=5). Even though we are willing to sacrifice computational time for ease of implementation, these values are significant and suggest that the method must be optimized to be of full practical interest.

In this case, we have two avenues for optimization, either we can prioritize reducing the CPU time or the numerical error. Since most of the computational time is dedicated to solving the least square system, accelerating the method would entail speeding up the linear solver and ultimately reducing the condition number. The method offers a number of free parameters and tuning knobs, such as the positions of the evaluation points, the number of neighbors, the basis functions selection, or the possibility to weight equations differently in the overdetermined system. We believe the method can be optimized by investigating the impact of these parameters on the accuracy of the solution and conditioning of the system.

Despite our presentation and analysis being restricted to low-dimensional examples and an un-optimized approach, we see that the LSQD method is a robust and versatile tool for computational scientists. Its minimal theoretical and programming barriers render it accessible to a broad scientific community, and we anticipate that further research endeavors could propel it into a formidable simulation paradigm.

%%%%%%%%%%%%%%%%%%%%%%%%%%%%%%%%%%%%%%%%%%%%%%%%%%
%%%%%%%        ACKNOWLEDGEMENTS            %%%%%%%
%%%%%%%%%%%%%%%%%%%%%%%%%%%%%%%%%%%%%%%%%%%%%%%%%%
\section{Acknowledgments}
The authors thank Arnold Kim for valuable discussions. This material is based upon work supported by the National Science Foundation under Grant No. DMS-1840265.

\section*{Credit author statement }  
{\bf Anna Kucherova}: 
Investigation, Formal Analysis, Software, Visualization, Writing - Original Draft, Writing - Review $\&$ Editing. 
 {\bf Gbocho Terasaki}: 
Formal Analysis, Investigation, Software, Visualization, Writing - Original Draft, Writing - Review $\&$ Editing.
{\bf Selma Strango}: Preliminary Investigation, Formal Analysis, Software, Writing - Review $\&$ Editing.
{\bf Maxime Theillard}: Conceptualization, Supervision, Formal Analysis, Methodology, Software, Project Administration, Writing - Original Draft, Writing - Review $\&$ Editing.

\newpage

\bibliographystyle{abbrv}
\bibliography{refs.bib}
@article{NGAI2011559,
title = "The application of data mining techniques in financial fraud detection: A classification framework and an academic review of literature",
journal = "Decision Support Systems",
volume = "50",
number = "3",
pages = "559 - 569",
year = "2011",
note = "On quantitative methods for detection of financial fraud",
issn = "0167-9236",
doi = "https://doi.org/10.1016/j.dss.2010.08.006",
url = "http://www.sciencedirect.com/science/article/pii/S0167923610001302",
author = "E.W.T. Ngai and Yong Hu and Y.H. Wong and Yijun Chen and Xin Sun",
keywords = "Financial fraud, Fraud detection, Literature review, Data mining, Business intelligence",
abstract = "This paper presents a review of — and classification scheme for — the literature on the application of data mining techniques for the detection of financial fraud. Although financial fraud detection (FFD) is an emerging topic of great importance, a comprehensive literature review of the subject has yet to be carried out. This paper thus represents the first systematic, identifiable and comprehensive academic literature review of the data mining techniques that have been applied to FFD. 49 journal articles on the subject published between 1997 and 2008 was analyzed and classified into four categories of financial fraud (bank fraud, insurance fraud, securities and commodities fraud, and other related financial fraud) and six classes of data mining techniques (classification, regression, clustering, prediction, outlier detection, and visualization). The findings of this review clearly show that data mining techniques have been applied most extensively to the detection of insurance fraud, although corporate fraud and credit card fraud have also attracted a great deal of attention in recent years. In contrast, we find a distinct lack of research on mortgage fraud, money laundering, and securities and commodities fraud. The main data mining techniques used for FFD are logistic models, neural networks, the Bayesian belief network, and decision trees, all of which provide primary solutions to the problems inherent in the detection and classification of fraudulent data. This paper also addresses the gaps between FFD and the needs of the industry to encourage additional research on neglected topics, and concludes with several suggestions for further FFD research."
}

@article{BLOMQUIST2024112695,
title = {Stable nodal projection method on octree grids},
journal = {Journal of Computational Physics},
volume = {499},
pages = {112695},
year = {2024},
issn = {0021-9991},
doi = {https://doi.org/10.1016/j.jcp.2023.112695},
url = {https://www.sciencedirect.com/science/article/pii/S0021999123007908},
author = {Matthew Blomquist and Scott R. West and Adam L. Binswanger and Maxime Theillard},
keywords = {Incompressible Navier-Stokes, Collocated node-based, Projection, Stability, Octree grids, Sharp interface},
abstract = {We propose a novel collocated projection method for solving the incompressible Navier-Stokes equations with arbitrary boundaries. Our approach employs non-graded octree grids, where all variables are stored at the nodes. To discretize the viscosity and projection steps, we utilize supra-convergent finite difference approximations with sharp boundary treatments. We demonstrate the stability of our projection on uniform grids, identify a sufficient stability condition on adaptive grids, and validate these findings numerically. We further demonstrate the accuracy and capabilities of our solver with several canonical two- and three-dimensional simulations of incompressible fluid flows. Overall, our method is second-order accurate, allows for dynamic grid adaptivity with arbitrary geometries, and reduces the overhead in code development through data collocation.}
}

@article{NGUYEN2008763,
title = {Meshless methods: A review and computer implementation aspects},
journal = {Mathematics and Computers in Simulation},
volume = {79},
number = {3},
pages = {763-813},
year = {2008},
issn = {0378-4754},
doi = {https://doi.org/10.1016/j.matcom.2008.01.003},
url = {https://www.sciencedirect.com/science/article/pii/S0378475408000062},
author = {Vinh Phu Nguyen and Timon Rabczuk and Stéphane Bordas and Marc Duflot},
keywords = {Meshless methods, Intrinsic enrichment, Extrinsic discontinuities, Computer implementation, MATLAB},
abstract = {The aim of this manuscript is to give a practical overview of meshless methods (for solid mechanics) based on global weak forms through a simple and well-structured MATLAB code, to illustrate our discourse. The source code is available for download on our website and should help students and researchers get started with some of the basic meshless methods; it includes intrinsic and extrinsic enrichment, point collocation methods, several boundary condition enforcement schemes and corresponding test cases. Several one and two-dimensional examples in elastostatics are given including weak and strong discontinuities and testing different ways of enforcing essential boundary conditions.}
}

@inproceedings{Samet1988AnOO,
  title={An Overview of Quadtrees, Octrees, and Related Hierarchical Data Structures},
  author={Hanan Samet},
  year={1988}
}

@article{SARSCOV2nucleocapsid,
author = {Chang, C., Michalska, K., Jedrzejczak, R., Maltseva, N., Endres, M., Godzik, A., Kim, Y., Joachimiak, A.},
title = {Crystal structure of RNA binding domain of nucleocapsid phosphoprotein from SARS coronavirus 2},
journal = {(to be published)},
year = {2020},
doi = {10.2210/pdb6vyo/pdb},
}

@misc{Foldinghome,
  title = {{Folding @ home} },
  howpublished = {\url{http:/foldingathome.org}},
 
}

@article{Politzer1985,
author = {P Politzer  and P R Laurence  and K Jayasuriya },
title = {Molecular electrostatic potentials: an effective tool for the elucidation of biochemical phenomena.},
journal = {Environmental Health Perspectives},
volume = {61},
number = {},
pages = {191-202},
year = {1985},
doi = {10.1289/ehp.8561191},
URL = {https://ehp.niehs.nih.gov/doi/abs/10.1289/ehp.8561191},
eprint = {https://ehp.niehs.nih.gov/doi/pdf/10.1289/ehp.8561191},
    abstract = { The electrostatic potential V(r) that is created in the space around a molecule by its nuclei and electrons (treated as static distributions of charge) is a very useful property for analyzing and predicting molecular reactive behavior. It is rigorously defined and can be determined experimentally as well as computationally. The potential has been particularly useful as an indicator of the sites or regions of a molecule to which an approaching electrophile is initially attracted, and it has also been applied successfully to the study of interactions that involve a certain optimum relative orientation of the reactants, such as between a drug and its cellular receptor. A variety of methods for calculating V(r) is available, at different levels of rigor. For large biologically active molecules, multipole expansions and superposition of potentials computed for subunits have been found to be effective. A large number of chemical and biochemical systems and processes have now been studied in terms of electrostatic potentials. Three examples of such applications are surveyed in this paper. These deal with: (a) reactive properties of nucleic acids, including their component bases; (b) biological recognition processes, including drug-receptors and enzyme-substrate interactions; and (c) chemical carcinogenesis, referring specifically to the polycyclic aromatic hydrocarbons and halogenated olefins and their epoxides. For each of these areas, examples of the use of electrostatic potentials in elucidating structure-activity patterns are given. }
}

@ARTICLE{Shashikala2019,
  
AUTHOR={Shashikala, H. B. Mihiri and Chakravorty, Arghya and Alexov, Emil},   
	 
TITLE={Modeling Electrostatic Force in Protein-Protein Recognition},      
	
JOURNAL={Frontiers in Molecular Biosciences},      
	
VOLUME={6},      

PAGES={94},     
	
YEAR={2019},      
	  
URL={https://www.frontiersin.org/article/10.3389/fmolb.2019.00094},       
	
DOI={10.3389/fmolb.2019.00094},      
	
ISSN={2296-889X},   
   
ABSTRACT={Electrostatic interactions are important for understanding molecular interactions, since they are long-range interactions and can guide binding partners to their correct binding positions. To investigate the role of electrostatic forces in molecular recognition, we calculated electrostatic forces between binding partners separated at various distances. The investigation was done on a large set of 275 protein complexes using recently developed DelPhiForce tool and in parallel, evaluating the total electrostatic force via electrostatic association energy. To accomplish the goal, we developed a method to find an appropriate direction to move one chain of protein complex away from its bound position and then calculate the corresponding electrostatic force as a function of separation distance. It is demonstrated that at large distances between the partners, the electrostatic force (magnitude and direction) is consistent among the protocols used and the main factors contributing to it are the net charge of the partners and their interfaces. However, at short distances, where partners form specific pair-wise interactions or de-solvation penalty becomes significant, the outcome depends on the precise balance of these factors. Based on the electrostatic force profile (force as a function of distance), we group the cases into four distinctive categories, among which the most intriguing is the case termed “soft landing.” In this case, the electrostatic force at large distances is favorable assisting the partners to come together, while at short distance it opposes binding, and thus slows down the approach of the partners toward their physical binding.}
}

@article{AZN_LANCET,
journal = {The Lancet},
year = {2020},
author = {Pedro M Folegatti and Katie J Ewer and Parvinder K Aley, 
Brian Angus and 
Prof Stephan Becker and
Sandra Belij-Rammerstorfer and 
Duncan Bellamy and
Sagida Bibi and
Mustapha Bittaye and
Elizabeth A Clutterbuck and 
Christina Dold and
Prof Saul N Faust and 
Prof Adam Finn and 
Amy L Flaxman and 
Bassam Hallis and
Prof Paul Heath and 
Daniel Jenkin and 
Rajeka Lazarus and 
Rebecca Makinson and 
Angela M Minassian and 
Katrina M Pollock and
Maheshi Ramasamy and 
Hannah Robinson and
Matthew Snape and 
Richard Tarrant and 
Merryn Voysey and
Catherine Green and 
Alexander D Douglas and
Prof Adrian V S Hill and
Teresa Lambe and 
Sarah C Gilbert and 
Andrew J Pollard},
title = {Safety and immunogenicity of the ChAdOx1 nCoV-19 vaccine against SARS-CoV-2: a preliminary report of a phase 1/2, single-blind, randomised controlled trial},
}

@article{AZN_LANCET1,
journal = {The Lancet},
year = {2020},
author = {Pedro M Folegatti and Katie J Ewer and Parvinder K Aley},
title = {Safety and immunogenicity of the ChAdOx1 nCoV-19 vaccine against SARS-CoV-2: a preliminary report of a phase 1/2, single-blind, randomised controlled trial},
}

@INPROCEEDINGS{6406271,

  author={Y. {Chen} and Y. {Zhou} and S. {Zhu} and H. {Xu}},

  booktitle={2012 International Conference on Privacy, Security, Risk and Trust and 2012 International Confernece on Social Computing}, 

  title={Detecting Offensive Language in Social Media to Protect Adolescent Online Safety}, 

  year={2012},

  volume={},

  number={},

  pages={71-80},}
  
@Article{Stilgoe2019,
author={Stilgoe, Jack},
title={Self-driving cars will take a while to get right},
journal={Nature Machine Intelligence},
year={2019},
month={May},
day={01},
volume={1},
number={5},
pages={202-203},
abstract={If we are to realize the potential of self-driving cars, we need to recognize the limits of machine learning. We should not pretend self-driving cars are around the corner: it will still take substantial time and effort to integrate the technology safely and fairly into our societies.},
issn={2522-5839},
doi={10.1038/s42256-019-0046-z},
url={https://doi.org/10.1038/s42256-019-0046-z}
}

@article{HEYDARI2022110755,
title = {Conservative finite volume method on deforming geometries: The case of protein aggregation in dividing yeast cells},
journal = {Journal of Computational Physics},
volume = {448},
pages = {110755},
year = {2022},
issn = {0021-9991},
doi = {https://doi.org/10.1016/j.jcp.2021.110755},
url = {https://www.sciencedirect.com/science/article/pii/S0021999121006501},
author = {A. Ali Heydari and Suzanne S. Sindi and Maxime Theillard},
keywords = {Conservative, Finite volume, Level set, Octree, Protein aggregation, Dividing cell},
abstract = {Motivated by experimental observations of asymmetric protein aggregates distributions in dividing yeast cells, we present a conservative finite volume approach for reaction-diffusion systems defined over deforming geometries. The key idea of our approach is to use spatio-temporal control volumes instead of integrating the time-discretized equations in space, as it is common practice. Both our theoretical and computational results demonstrate the convergence of our method and highlight how traditional approaches can lead to inaccurate solutions. We employ this novel approach to investigate the partitioning of protein aggregates in dividing yeast cells, leveraging the flexibility of the level set method to construct realistic biological geometries. Using a simple reaction-diffusion model, we find that spatial heterogeneity in yeast cells during division can alone create asymmetries in the concentration of protein aggregates. Moreover, we find that obstructing intracellular entities, such as nuclei or insoluble protein compartments, amplify these asymmetries, suggesting that they may play an essential role in regulating molecular partitioning. Beyond these findings, our results illustrate the flexibility of our approach and its potential to design realistic predictive tools to explore intracellular bio-mechanisms.}
}

@Article{D3SM01196H,
author ="Fylling, Cayce and Tamayo, Joshua and Gopinath, Arvind and Theillard, Maxime",
title  ="Multi-population dissolution in confined active fluids",
journal  ="Soft Matter",
year  ="2024",
volume  ="20",
issue  ="7",
pages  ="1392-1409",
publisher  ="The Royal Society of Chemistry",
doi  ="10.1039/D3SM01196H",
url  ="http://dx.doi.org/10.1039/D3SM01196H",
abstract  ="Autonomous out-of-equilibrium agents or cells in suspension are ubiquitous in biology and engineering. Turning chemical energy into mechanical stress{,} they generate activity in their environment{,} which may trigger spontaneous large-scale dynamics. Often{,} these systems are composed of multiple populations that may reflect the coexistence of multiple species{,} differing phenotypes{,} or chemically varying agents in engineered settings. Here{,} we present a new method for modeling such multi-population active fluids subject to confinement. We use a continuum multi-scale mean-field approach to represent each phase by its first three orientational moments and couple their evolution with those of the suspending fluid. The resulting coupled system is solved using a parallel adaptive level-set-based solver for high computational efficiency and maximal flexibility in the confinement geometry. Motivated by recent experimental work{,} we employ our method to study the spatiotemporal dynamics of confined bacterial suspensions and swarms dominated by fluid hydrodynamic effects. Our in silico explorations reproduce observed emergent collective patterns{,} including features of active dissolution in two-population active-passive swarms{,} with results clearly suggesting that hydrodynamic effects dominate dissolution dynamics. Our work lays the foundation for a systematic characterization and study of collective phenomena in natural or synthetic multi-population systems such as bacteria colonies{,} bird flocks{,} fish schools{,} colloid swimmers{,} or programmable active matter."}

@article{JIANG199013,
title = {Least-squares finite element method for fluid dynamics},
journal = {Computer Methods in Applied Mechanics and Engineering},
volume = {81},
number = {1},
pages = {13-37},
year = {1990},
issn = {0045-7825},
doi = {https://doi.org/10.1016/0045-7825(90)90139-D},
url = {https://www.sciencedirect.com/science/article/pii/004578259090139D},
author = {Bo-Nan Jiang and Louis A. Povinelli},
abstract = {This paper gives an overview of new developments of the least-squares finite element method (LSFEM) in fluid dynamics. Special emphasis is placed on the universality of LSFEM, the symmetry and positiveness of the algebraic systems obtained from LSFEM, the accommodation of LSFEM to equal-order interpolations for incompressible viscous flows, and the natural numerical dissipation of LSFEM for convective transport problems and high-speed compressible flows. The performance of LSFEM is illustrated by numerical examples.}
}

@Inbook{Jiang1998,
author="Jiang, Bo-nan",
title="Basis of LSFEM",
bookTitle="The Least-Squares Finite Element Method: Theory and Applications in Computational Fluid Dynamics and Electromagnetics",
year="1998",
publisher="Springer Berlin Heidelberg",
address="Berlin, Heidelberg",
pages="47--79",
abstract="In this chapter an attempt is made to present a unified theory and formulation of the least-squares finite element method for multi-dimensional problems so that in principle, it is not necessary to repeat the same argument in different instances. In the following chapters, these theorems, inequalities and formulations will be extensively utilized. The mathematics is kept as simple as possible, since our purpose is to become familiar with certain modern concepts that will help readers understand the basic theory, formulation and properties of the least-squares method.",
isbn="978-3-662-03740-9",
doi="10.1007/978-3-662-03740-9_4",
url="https://doi.org/10.1007/978-3-662-03740-9_4"
}

@article{KUCHEROVA2021110591,
title = {Computational modeling of protein conformational changes - Application to the opening SARS-CoV-2 spike},
journal = {Journal of Computational Physics},
volume = {444},
pages = {110591},
year = {2021},
issn = {0021-9991},
doi = {https://doi.org/10.1016/j.jcp.2021.110591},
url = {https://www.sciencedirect.com/science/article/pii/S0021999121004861},
author = {Anna Kucherova and Selma Strango and Shahar Sukenik and Maxime Theillard},
keywords = {SARS-CoV-2, Covid-19, Spike protein, Molecular trajectory, Poisson-Boltzmann, Multiscale modeling},
abstract = {We present a new approach to compute and analyze the dynamical electro-geometric properties of proteins undergoing conformational changes. The molecular trajectory is obtained from Markov state models, and the electrostatic potential is calculated using the continuum Poisson-Boltzmann equation. The numerical electric potential is constructed using a parallel sharp numerical solver implemented on adaptive Octree grids. We introduce novel a posteriori error estimates to quantify the solution's accuracy on the molecular surface. To illustrate the approach, we consider the opening of the SARS-CoV-2 spike protein using the recent molecular trajectory simulated through the Folding@home initiative. We analyze our results, focusing on the characteristics of the receptor-binding domain and its vicinity. This work lays the foundation for a new class of hybrid computational approaches, producing high-fidelity dynamical computational measurements serving as a basis for protein bio-mechanism investigations.}
}

@article{googlelanguage2016,
  author    = {Yonghui Wu and
               Mike Schuster and
               Zhifeng Chen and
               Quoc V. Le and
               Mohammad Norouzi and
               Wolfgang Macherey and
               Maxim Krikun and
               Yuan Cao and
               Qin Gao and
               Klaus Macherey and
               Jeff Klingner and
               Apurva Shah and
               Melvin Johnson and
               Xiaobing Liu and
               Lukasz Kaiser and
               Stephan Gouws and
               Yoshikiyo Kato and
               Taku Kudo and
               Hideto Kazawa and
               Keith Stevens and
               George Kurian and
               Nishant Patil and
               Wei Wang and
               Cliff Young and
               Jason Smith and
               Jason Riesa and
               Alex Rudnick and
               Oriol Vinyals and
               Greg Corrado and
               Macduff Hughes and
               Jeffrey Dean},
  title     = {Google's Neural Machine Translation System: Bridging the Gap between
               Human and Machine Translation},
  journal   = {CoRR},
  volume    = {abs/1609.08144},
  year      = {2016},
  url       = {http://arxiv.org/abs/1609.08144},
  archivePrefix = {arXiv},
  eprint    = {1609.08144},
  timestamp = {Thu, 14 Mar 2019 09:34:18 +0100},
  biburl    = {https://dblp.org/rec/journals/corr/WuSCLNMKCGMKSJL16.bib},
  bibsource = {dblp computer science bibliography, https://dblp.org}
}

@article{Weiss2018,
  author    = {Ron J. Weiss and
               Jan Chorowski and
               Navdeep Jaitly and
               Yonghui Wu and
               Zhifeng Chen},
  title     = {Sequence-to-Sequence Models Can Directly Transcribe Foreign Speech},
  journal   = {CoRR},
  volume    = {abs/1703.08581},
  year      = {2017},
  url       = {http://arxiv.org/abs/1703.08581},
  archivePrefix = {arXiv},
  eprint    = {1703.08581},
  timestamp = {Mon, 13 Aug 2018 16:46:35 +0200},
  biburl    = {https://dblp.org/rec/journals/corr/WeissCJWC17.bib},
  bibsource = {dblp computer science bibliography, https://dblp.org}
}

@inproceedings{Nobata2016,
author = {Nobata, Chikashi and Tetreault, Joel and Thomas, Achint and Mehdad, Yashar and Chang, Yi},
title = {Abusive Language Detection in Online User Content},
year = {2016},
isbn = {9781450341431},
publisher = {International World Wide Web Conferences Steering Committee},
address = {Republic and Canton of Geneva, CHE},
url = {https://doi.org/10.1145/2872427.2883062},
doi = {10.1145/2872427.2883062},
booktitle = {Proceedings of the 25th International Conference on World Wide Web},
pages = {145–153},
numpages = {9},
keywords = {hate speech, stylistic classification, natural language processing, abusive language, nlp, discourse classification},
location = {Montr\'{e}al, Qu\'{e}bec, Canada},
series = {WWW ’16}
}

@article{Baskin2016,
author = {Igor I. Baskin and David Winkler and Igor V. Tetko},
title = {A renaissance of neural networks in drug discovery},
journal = {Expert Opinion on Drug Discovery},
volume = {11},
number = {8},
pages = {785-795},
year  = {2016},
publisher = {Taylor \& Francis},
doi = {10.1080/17460441.2016.1201262},
    note ={PMID: 27295548},

URL = { 
        https://doi.org/10.1080/17460441.2016.1201262
    
},
eprint = { 
        https://doi.org/10.1080/17460441.2016.1201262
    
}

}

@Article{Ekins2016,
author={Ekins, Sean},
title={The Next Era: Deep Learning in Pharmaceutical Research},
journal={Pharmaceutical Research},
year={2016},
month={Nov},
day={01},
volume={33},
number={11},
pages={2594-2603},
abstract={Over the past decade we have witnessed the increasing sophistication of machine learning algorithms applied in daily use from internet searches, voice recognition, social network software to machine vision software in cameras, phones, robots and self-driving cars. Pharmaceutical research has also seen its fair share of machine learning developments. For example, applying such methods to mine the growing datasets that are created in drug discovery not only enables us to learn from the past but to predict a molecule's properties and behavior in future. The latest machine learning algorithm garnering significant attention is deep learning, which is an artificial neural network with multiple hidden layers. Publications over the last 3 years suggest that this algorithm may have advantages over previous machine learning methods and offer a slight but discernable edge in predictive performance. The time has come for a balanced review of this technique but also to apply machine learning methods such as deep learning across a wider array of endpoints relevant to pharmaceutical research for which the datasets are growing such as physicochemical property prediction, formulation prediction, absorption, distribution, metabolism, excretion and toxicity (ADME/Tox), target prediction and skin permeation, etc. We also show that there are many potential applications of deep learning beyond cheminformatics. It will be important to perform prospective testing (which has been carried out rarely to date) in order to convince skeptics that there will be benefits from investing in this technique.},
issn={1573-904X},
doi={10.1007/s11095-016-2029-7},
url={https://doi.org/10.1007/s11095-016-2029-7}
}

@article{Gawehn2016,
author = {Gawehn, Erik and Hiss, Jan A. and Schneider, Gisbert},
title = {Deep Learning in Drug Discovery},
journal = {Molecular Informatics},
volume = {35},
number = {1},
pages = {3-14},
keywords = {bioinformatics, cheminformatics, drug design, machine-learning, neural network, virtual screening},
doi = {10.1002/minf.201501008},
url = {https://onlinelibrary.wiley.com/doi/abs/10.1002/minf.201501008},
eprint = {https://onlinelibrary.wiley.com/doi/pdf/10.1002/minf.201501008},
abstract = {Abstract Artificial neural networks had their first heyday in molecular informatics and drug discovery approximately two decades ago. Currently, we are witnessing renewed interest in adapting advanced neural network architectures for pharmaceutical research by borrowing from the field of “deep learning”. Compared with some of the other life sciences, their application in drug discovery is still limited. Here, we provide an overview of this emerging field of molecular informatics, present the basic concepts of prominent deep learning methods and offer motivation to explore these techniques for their usefulness in computer-assisted drug discovery and design. We specifically emphasize deep neural networks, restricted Boltzmann machine networks and convolutional networks.},
year = {2016}
}

@article{LAVECCHIA20192017,
title = "Deep learning in drug discovery: opportunities, challenges and future prospects",
journal = "Drug Discovery Today",
volume = "24",
number = "10",
pages = "2017 - 2032",
year = "2019",
issn = "1359-6446",
doi = "https://doi.org/10.1016/j.drudis.2019.07.006",
url = "http://www.sciencedirect.com/science/article/pii/S135964461930282X",
author = "Antonio Lavecchia",
abstract = "Artificial Intelligence (AI) is an area of computer science that simulates the structures and operating principles of the human brain. Machine learning (ML) belongs to the area of AI and endeavors to develop models from exposure to training data. Deep Learning (DL) is another subset of AI, where models represent geometric transformations over many different layers. This technology has shown tremendous potential in areas such as computer vision, speech recognition and natural language processing. More recently, DL has also been successfully applied in drug discovery. Here, I analyze several relevant DL applications and case studies, providing a detailed view of the current state-of-the-art in drug discovery and highlighting not only the problematic issues, but also the successes and opportunities for further advances."
}

@article{CHEN20181241,
title = "The rise of deep learning in drug discovery",
journal = "Drug Discovery Today",
volume = "23",
number = "6",
pages = "1241 - 1250",
year = "2018",
issn = "1359-6446",
doi = "https://doi.org/10.1016/j.drudis.2018.01.039",
url = "http://www.sciencedirect.com/science/article/pii/S1359644617303598",
author = "Hongming Chen and Ola Engkvist and Yinhai Wang and Marcus Olivecrona and Thomas Blaschke",
abstract = "Over the past decade, deep learning has achieved remarkable success in various artificial intelligence research areas. Evolved from the previous research on artificial neural networks, this technology has shown superior performance to other machine learning algorithms in areas such as image and voice recognition, natural language processing, among others. The first wave of applications of deep learning in pharmaceutical research has emerged in recent years, and its utility has gone beyond bioactivity predictions and has shown promise in addressing diverse problems in drug discovery. Examples will be discussed covering bioactivity prediction, de novo molecular design, synthesis prediction and biological image analysis."
}

@ARTICLE{8661629,

  author={Y. {Zhao} and Q. {Chen} and W. {Cao} and J. {Yang} and J. {Xiong} and G. {Gui}},

  journal={IEEE Access}, 

  title={Deep Learning for Risk Detection and Trajectory Tracking at Construction Sites}, 

  year={2019},

  volume={7},

  number={},

  pages={30905-30912},}

@article{PALTRINIERI2019475,
title = "Learning about risk: Machine learning for risk assessment",
journal = "Safety Science",
volume = "118",
pages = "475 - 486",
year = "2019",
issn = "0925-7535",
doi = "https://doi.org/10.1016/j.ssci.2019.06.001",
url = "http://www.sciencedirect.com/science/article/pii/S0925753518311184",
author = "Nicola Paltrinieri and Louise Comfort and Genserik Reniers",
keywords = "Risk assessment, Dynamic risk analysis, Machine learning, Deep learning",
abstract = "Risk assessment has a primary role in safety-critical industries. However, it faces a series of overall challenges, partially related to technology advancements and increasing needs. There is currently a call for continuous risk assessment, improvement in learning past lessons and definition of techniques to process relevant data, which are to be coupled with adequate capability to deal with unexpected events and provide the right support to enable risk management. Through this work, we suggest a risk assessment approach based on machine learning. In particular, a deep neural network (DNN) model is developed and tested for a drive-off scenario involving an Oil \& Gas drilling rig. Results show reasonable accuracy for DNN predictions and general suitability to (partially) overcome risk assessment challenges. Nevertheless, intrinsic model limitations should be taken into account and appropriate model selection and customization should be carefully carried out to deliver appropriate support for safety-related decision-making."
}

@ARTICLE{8316850,

  author={C. {Jiang} and J. {Song} and G. {Liu} and L. {Zheng} and W. {Luan}},

  journal={IEEE Internet of Things Journal}, 

  title={Credit Card Fraud Detection: A Novel Approach Using Aggregation Strategy and Feedback Mechanism}, 

  year={2018},

  volume={5},

  number={5},

  pages={3637-3647},}

@article{DORNADULA2019631,
title = "Credit Card Fraud Detection using Machine Learning Algorithms",
journal = "Procedia Computer Science",
volume = "165",
pages = "631 - 641",
year = "2019",
note = "2nd International Conference on Recent Trends in Advanced Computing ICRTAC -DISRUP - TIV INNOVATION , 2019 November 11-12, 2019",
issn = "1877-0509",
doi = "https://doi.org/10.1016/j.procs.2020.01.057",
url = "http://www.sciencedirect.com/science/article/pii/S187705092030065X",
author = "Vaishnavi Nath Dornadula and S Geetha",
keywords = "Card-Not-Present frauds, Card-Present-Frauds, Concept Drift",
abstract = "Credit card frauds are easy and friendly targets. E-commerce and many other online sites have increased the online payment modes, increasing the risk for online frauds. Increase in fraud rates, researchers started using different machine learning methods to detect and analyse frauds in online transactions. The main aim of the paper is to design and develop a novel fraud detection method for Streaming Transaction Data, with an objective, to analyse the past transaction details of the customers and extract the behavioural patterns. Where cardholders are clustered into different groups based on their transaction amount. Then using sliding window strategy [1], to aggregate the transaction made by the cardholders from different groups so that the behavioural pattern of the groups can be extracted respectively. Later different classifiers [3],[5],[6],[8] are trained over the groups separately. And then the classifier with better rating score can be chosen to be one of the best methods to predict frauds. Thus, followed by a feedback mechanism to solve the problem of concept drift [1]. In this paper, we worked with European credit card fraud dataset."
}

@article{RAISSI2018125,
title = "Hidden physics models: Machine learning of nonlinear partial differential equations",
journal = "Journal of Computational Physics",
volume = "357",
pages = "125 - 141",
year = "2018",
issn = "0021-9991",
doi = "https://doi.org/10.1016/j.jcp.2017.11.039",
url = "http://www.sciencedirect.com/science/article/pii/S0021999117309014",
author = "Maziar Raissi and George Em Karniadakis",
keywords = "Probabilistic machine learning, System identification, Bayesian modeling, Uncertainty quantification, Fractional equations, Small data",
abstract = "While there is currently a lot of enthusiasm about “big data”, useful data is usually “small” and expensive to acquire. In this paper, we present a new paradigm of learning partial differential equations from small data. In particular, we introduce hidden physics models, which are essentially data-efficient learning machines capable of leveraging the underlying laws of physics, expressed by time dependent and nonlinear partial differential equations, to extract patterns from high-dimensional data generated from experiments. The proposed methodology may be applied to the problem of learning, system identification, or data-driven discovery of partial differential equations. Our framework relies on Gaussian processes, a powerful tool for probabilistic inference over functions, that enables us to strike a balance between model complexity and data fitting. The effectiveness of the proposed approach is demonstrated through a variety of canonical problems, spanning a number of scientific domains, including the Navier–Stokes, Schrödinger, Kuramoto–Sivashinsky, and time dependent linear fractional equations. The methodology provides a promising new direction for harnessing the long-standing developments of classical methods in applied mathematics and mathematical physics to design learning machines with the ability to operate in complex domains without requiring large quantities of data."
}
@article{RAISSI2019686,
title = "Physics-informed neural networks: A deep learning framework for solving forward and inverse problems involving nonlinear partial differential equations",
journal = "Journal of Computational Physics",
volume = "378",
pages = "686 - 707",
year = "2019",
issn = "0021-9991",
doi = "https://doi.org/10.1016/j.jcp.2018.10.045",
url = "http://www.sciencedirect.com/science/article/pii/S0021999118307125",
author = "M. Raissi and P. Perdikaris and G.E. Karniadakis",
keywords = "Data-driven scientific computing, Machine learning, Predictive modeling, Runge–Kutta methods, Nonlinear dynamics",
abstract = "We introduce physics-informed neural networks – neural networks that are trained to solve supervised learning tasks while respecting any given laws of physics described by general nonlinear partial differential equations. In this work, we present our developments in the context of solving two main classes of problems: data-driven solution and data-driven discovery of partial differential equations. Depending on the nature and arrangement of the available data, we devise two distinct types of algorithms, namely continuous time and discrete time models. The first type of models forms a new family of data-efficient spatio-temporal function approximators, while the latter type allows the use of arbitrarily accurate implicit Runge–Kutta time stepping schemes with unlimited number of stages. The effectiveness of the proposed framework is demonstrated through a collection of classical problems in fluids, quantum mechanics, reaction–diffusion systems, and the propagation of nonlinear shallow-water waves."
}

@article{Botu2015,
author = {Botu, Venkatesh and Ramprasad, Rampi},
title = {Adaptive machine learning framework to accelerate ab initio molecular dynamics},
journal = {International Journal of Quantum Chemistry},
volume = {115},
number = {16},
pages = {1074-1083},
keywords = {ab initio molecular dynamics, accelerate, machine learning, fingerprint, adaptive},
doi = {10.1002/qua.24836},
url = {https://onlinelibrary.wiley.com/doi/abs/10.1002/qua.24836},
eprint = {https://onlinelibrary.wiley.com/doi/pdf/10.1002/qua.24836},
abstract = {Quantum mechanics-based ab initio molecular dynamics (MD) simulation schemes offer an accurate and direct means to monitor the time evolution of materials. Nevertheless, the expensive and repetitive energy and force computations required in such simulations lead to significant bottlenecks. Here, we lay the foundations for an accelerated ab initio MD approach integrated with a machine learning framework. The proposed algorithm learns from previously visited configurations in a continuous and adaptive manner on-the-fly, and predicts (with chemical accuracy) the energies and atomic forces of a new configuration at a minuscule fraction of the time taken by conventional ab initio methods. Key elements of this new accelerated ab initio MD paradigm include representations of atomic configurations by numerical fingerprints, a learning algorithm to map the fingerprints to the properties, a decision engine that guides the choice of the prediction scheme, and requisite amount of ab initio data. The performance of each aspect of the proposed scheme is critically evaluated for Al in several different chemical environments. This work has enormous implications beyond ab initio MD acceleration. It can also lead to accelerated structure and property prediction schemes, and accurate force fields. © 2014 Wiley Periodicals, Inc.},
year = {2015}
}

@Article{Smith2019,
author={Smith, Justin S.
and Nebgen, Benjamin T.
and Zubatyuk, Roman
and Lubbers, Nicholas
and Devereux, Christian
and Barros, Kipton
and Tretiak, Sergei
and Isayev, Olexandr
and Roitberg, Adrian E.},
title={Approaching coupled cluster accuracy with a general-purpose neural network potential through transfer learning},
journal={Nature Communications},
year={2019},
month={Jul},
day={01},
volume={10},
number={1},
pages={2903},
abstract={Computational modeling of chemical and biological systems at atomic resolution is a crucial tool in the chemist's toolset. The use of computer simulations requires a balance between cost and accuracy: quantum-mechanical methods provide high accuracy but are computationally expensive and scale poorly to large systems, while classical force fields are cheap and scalable, but lack transferability to new systems. Machine learning can be used to achieve the best of both approaches. Here we train a general-purpose neural network potential (ANI-1ccx) that approaches CCSD(T)/CBS accuracy on benchmarks for reaction thermochemistry, isomerization, and drug-like molecular torsions. This is achieved by training a network to DFT data then using transfer learning techniques to retrain on a dataset of gold standard QM calculations (CCSD(T)/CBS) that optimally spans chemical space. The resulting potential is broadly applicable to materials science, biology, and chemistry, and billions of times faster than CCSD(T)/CBS calculations.},
issn={2041-1723},
doi={10.1038/s41467-019-10827-4},
url={https://doi.org/10.1038/s41467-019-10827-4}
}

@article{Moreau2006,
author = {Moreau,A.  and Teytaud,O.  and Bertoglio,J. P. },
title = {Optimal estimation for large-eddy simulation of turbulence and application to the analysis of subgrid models},
journal = {Physics of Fluids},
volume = {18},
number = {10},
pages = {105101},
year = {2006},
doi = {10.1063/1.2357974},

URL = { 
        https://doi.org/10.1063/1.2357974
    
},
eprint = { 
        https://doi.org/10.1063/1.2357974
    
}

}
@article{fukami_fukagata_taira_2019, title={Super-resolution reconstruction of turbulent flows with machine learning}, volume={870}, DOI={10.1017/jfm.2019.238}, journal={Journal of Fluid Mechanics}, publisher={Cambridge University Press}, author={Fukami, Kai and Fukagata, Koji and Taira, Kunihiko}, year={2019}, pages={106–120}}

@article{UNGER20081994,
title = "Coupling of scales in a multiscale simulation using neural networks",
journal = "Computers \& Structures",
volume = "86",
number = "21",
pages = "1994 - 2003",
year = "2008",
issn = "0045-7949",
doi = "https://doi.org/10.1016/j.compstruc.2008.05.004",
url = "http://www.sciencedirect.com/science/article/pii/S0045794908001430",
author = "Jörg F. Unger and Carsten Könke",
keywords = "Neural network, Homogenization, Multilayer perceptron, Support vector machines, Constitutive modelling",
abstract = "Multiscale approaches require the coupling of models on different spatial scales. In this paper, a coupling using neural networks is proposed. Based on a set of mesoscale simulations of concrete a system of neural networks is trained to approximate the response. A macroscale constitutive model is obtained by homogenizing the mesoscale response. Special focus is put on the mesh sensitivity, since the mesoscale model includes softening and consequently the width of the localization zone compared to the dimension of the mesoscale model is an additional parameter in the model."
}

@article{MingTryggvason2015,
author = {Ma,Ming  and Lu,Jiacai  and Tryggvason,Gretar },
title = {Using statistical learning to close two-fluid multiphase flow equations for a simple bubbly system},
journal = {Physics of Fluids},
volume = {27},
number = {9},
pages = {092101},
year = {2015},
doi = {10.1063/1.4930004},

URL = { 
        https://aip.scitation.org/doi/abs/10.1063/1.4930004
    
},
eprint = { 
        https://aip.scitation.org/doi/pdf/10.1063/1.4930004
    
}

}

@article{GIBOU2019442,
title = "Sharp interface approaches and deep learning techniques for multiphase flows",
journal = "Journal of Computational Physics",
volume = "380",
pages = "442 - 463",
year = "2019",
issn = "0021-9991",
doi = "https://doi.org/10.1016/j.jcp.2018.05.031",
url = "http://www.sciencedirect.com/science/article/pii/S0021999118303371",
author = "Frederic Gibou and David Hyde and Ron Fedkiw",
keywords = "Multiphase flows, Ghost-fluid method, Voronoi interface method, Surface tension, Quadtree and octree, Parallel",
abstract = "We present a review on numerical methods for simulating multiphase and free surface flows. We focus in particular on numerical methods that seek to preserve the discontinuous nature of the solutions across the interface between phases. We provide a discussion on the Ghost-Fluid and Voronoi Interface methods, on the treatment of surface tension forces that avoid stringent time step restrictions, on adaptive grid refinement techniques for improved efficiency and on parallel computing approaches. We present the results of some simulations obtained with these treatments in two and three spatial dimensions. We also provide a discussion of Machine Learning and Deep Learning techniques in the context of multiphase flows and propose several future potential research thrusts for using deep learning to enhance the study and simulation of multiphase flows."
}

@article{Chang6591,
author = {Chang, Fi-John and Yang, Han-Chung and Lu, Jau-Yau and Hong, Jian-Hao},
title = {Neural network modelling for mean velocity and turbulence intensities of steep channel flows},
journal = {Hydrological Processes},
volume = {22},
number = {2},
pages = {265-274},
keywords = {artificial neural network (ANN), velocity profile, turbulent open channel flow, log-law model, Reynolds stress model},
doi = {10.1002/hyp.6591},
url = {https://onlinelibrary.wiley.com/doi/abs/10.1002/hyp.6591},
eprint = {https://onlinelibrary.wiley.com/doi/pdf/10.1002/hyp.6591},
abstract = {Abstract The main purpose of this study is to evaluate the potential of simulating the profiles of the mean velocity and turbulence intensities for the steep open channel flows over a smooth boundary using artificial neural networks. In a laboratory flume, turbulent flow conditions were measured using a fibre-optic laser doppler velocimeter (FLDV). One thousand and sixty-four data sets were collected for different slopes and aspect ratios at different locations. These data sets were randomly split into two subsets, i.e. training and validation sets. The multi-layer functional link network (MFLN) was used to construct the simulation model based on the training data. The constructed MFLN models can almost perfectly simulate the velocity profile and turbulence intensity. The values of correlation coefficient (γ) are close to one and the values of root mean square error (RMSE) are close to zero in all conditions. The results demonstrate that the MFLN can precisely simulate the velocity profiles, while the log law and Reynolds stress model (RSM) are less effective when used to simulate the velocity profiles close to the side wall. The simulated longitudinal turbulence intensities yielded by the MFLN were also fairly consistent with the measured data, while the simulated vertical turbulence intensities by the RSM were not consistent with the measured data. Copyright © 2007 John Wiley \& Sons, Ltd.},
year = {2008}
}

@article{MILANO20021,
title = "Neural Network Modeling for Near Wall Turbulent Flow",
journal = "Journal of Computational Physics",
volume = "182",
number = "1",
pages = "1 - 26",
year = "2002",
issn = "0021-9991",
doi = "https://doi.org/10.1006/jcph.2002.7146",
url = "http://www.sciencedirect.com/science/article/pii/S0021999102971469",
author = "Michele Milano and Petros Koumoutsakos",
abstract = "A neural network methodology is developed in order to reconstruct the near wall field in a turbulent flow by exploiting flow fields provided by direct numerical simulations. The results obtained from the neural network methodology are compared with the results obtained from prediction and reconstruction using proper orthogonal decomposition (POD). Using the property that the POD is equivalent to a specific linear neural network, a nonlinear neural network extension is presented. It is shown that for a relatively small additional computational cost nonlinear neural networks provide us with improved reconstruction and prediction capabilities for the near wall velocity fields. Based on these results advantages and drawbacks of both approaches are discussed with an outlook toward the development of near wall models for turbulence modeling and control."
}

@article{Giralt2000,
author = {Giralt,Francesc  and Arenas,A.  and Ferre-Giné,J.  and Rallo,R.  and Kopp,G. A. },
title = {The simulation and interpretation of free turbulence with a cognitive neural system},
journal = {Physics of Fluids},
volume = {12},
number = {7},
pages = {1826-1835},
year = {2000},
doi = {10.1063/1.870430},

URL = { 
        https://doi.org/10.1063/1.870430
    
},
eprint = { 
        https://doi.org/10.1063/1.870430
    
}

}

@Article{Asproulis2013,
author={Asproulis, Nikolaos
and Drikakis, Dimitris},
title={An artificial neural network-based multiscale method for hybrid atomistic-continuum simulations},
journal={Microfluidics and Nanofluidics},
year={2013},
month={Oct},
day={01},
volume={15},
number={4},
pages={559-574},
abstract={This paper presents an artificial neural network-based multiscale method for coupling continuum and molecular simulations. Molecular dynamics modelling is employed as a local ``high resolution'' refinement of computational data required by the continuum computational fluid dynamics solver. The coupling between atomistic and continuum simulations is obtained by an artificial neural network (ANN) methodology. The ANN aims to optimise the transfer of information through minimisation of (1) the computational cost by avoiding repetitive atomistic simulations of nearly identical states, and (2) the fluctuation strength of the atomistic outputs that are fed back to the continuum solver. Results are presented for prototype flows such as the isothermal Couette flow with slip boundary conditions and the slip Couette flow with heat transfer.},
issn={1613-4990},
doi={10.1007/s10404-013-1154-4},
url={https://doi.org/10.1007/s10404-013-1154-4}
}

@article{Zhenwei2015,
  title = {Molecular Dynamics with On-the-Fly Machine Learning of Quantum-Mechanical Forces},
  author = {Li, Zhenwei and Kermode, James R. and De Vita, Alessandro},
  journal = {Phys. Rev. Lett.},
  volume = {114},
  issue = {9},
  pages = {096405},
  numpages = {5},
  year = {2015},
  month = {Mar},
  publisher = {American Physical Society},
  doi = {10.1103/PhysRevLett.114.096405},
  url = {https://link.aps.org/doi/10.1103/PhysRevLett.114.096405}
}

@article {Chmielae1603015,
	author = {Chmiela, Stefan and Tkatchenko, Alexandre and Sauceda, Huziel E. and Poltavsky, Igor and Sch{\"u}tt, Kristof T. and M{\"u}ller, Klaus-Robert},
	title = {Machine learning of accurate energy-conserving molecular force fields},
	volume = {3},
	number = {5},
	elocation-id = {e1603015},
	year = {2017},
	doi = {10.1126/sciadv.1603015},
	publisher = {American Association for the Advancement of Science},
	abstract = {Using conservation of energy{\textemdash}a fundamental property of closed classical and quantum mechanical systems{\textemdash}we develop an efficient gradient-domain machine learning (GDML) approach to construct accurate molecular force fields using a restricted number of samples from ab initio molecular dynamics (AIMD) trajectories. The GDML implementation is able to reproduce global potential energy surfaces of intermediate-sized molecules with an accuracy of 0.3 kcal mol-1 for energies and 1 kcal mol-1 {\r A}̊-1 for atomic forces using only 1000 conformational geometries for training. We demonstrate this accuracy for AIMD trajectories of molecules, including benzene, toluene, naphthalene, ethanol, uracil, and aspirin. The challenge of constructing conservative force fields is accomplished in our work by learning in a Hilbert space of vector-valued functions that obey the law of energy conservation. The GDML approach enables quantitative molecular dynamics simulations for molecules at a fraction of cost of explicit AIMD calculations, thereby allowing the construction of efficient force fields with the accuracy and transferability of high-level ab initio methods.},
	URL = {https://advances.sciencemag.org/content/3/5/e1603015},
	eprint = {https://advances.sciencemag.org/content/3/5/e1603015.full.pdf},
	journal = {Science Advances}
}

@Article{computation8010015,
AUTHOR = {Frank, Michael and Drikakis, Dimitris and Charissis, Vassilis},
TITLE = {Machine-Learning Methods for Computational Science and Engineering},
JOURNAL = {Computation},
VOLUME = {8},
YEAR = {2020},
NUMBER = {1},
ARTICLE-NUMBER = {15},
URL = {https://www.mdpi.com/2079-3197/8/1/15},
ISSN = {2079-3197},
ABSTRACT = {The re-kindled fascination in machine learning (ML), observed over the last few decades, has also percolated into natural sciences and engineering. ML algorithms are now used in scientific computing, as well as in data-mining and processing. In this paper, we provide a review of the state-of-the-art in ML for computational science and engineering. We discuss ways of using ML to speed up or improve the quality of simulation techniques such as computational fluid dynamics, molecular dynamics, and structural analysis. We explore the ability of ML to produce computationally efficient surrogate models of physical applications that circumvent the need for the more expensive simulation techniques entirely. We also discuss how ML can be used to process large amounts of data, using as examples many different scientific fields, such as engineering, medicine, astronomy and computing. Finally, we review how ML has been used to create more realistic and responsive virtual reality applications.},
DOI = {10.3390/computation8010015}
}

@article{DESPRES2020109275,
title = "Machine Learning design of Volume of Fluid schemes for compressible flows",
journal = "Journal of Computational Physics",
volume = "408",
pages = "109275",
year = "2020",
issn = "0021-9991",
doi = "https://doi.org/10.1016/j.jcp.2020.109275",
url = "http://www.sciencedirect.com/science/article/pii/S0021999120300498",
author = "Bruno Despr\'{e}s and Herv\'{e} Jourdren",
keywords = "VOF, CFD, ML",
abstract = "Our aim is to establish the feasibility of Machine-Learning-designed Volume of Fluid algorithms for compressible flows. We detail the incremental steps of the construction of a new family of Volume of Fluid-Machine Learning (VOF-ML) schemes adapted to bi-material compressible Euler calculations on Cartesian grids. An additivity principle is formulated for the Machine Learning datasets. We explain a key feature of this approach which is how to adapt the compressible solver to the preservation of natural symmetries. The VOF-ML schemes show good accuracy for advection of a variety of interfaces, including regular interfaces (straight lines and arcs of circle), Lipschitz interfaces (corners) and non Lipschitz triple point (the Trifolium test problem). Basic comparisons with a SLIC/Downwind scheme are presented together with elementary bi-material calculations with shocks."
}

@article {Brunton3932,
	author = {Brunton, Steven L. and Proctor, Joshua L. and Kutz, J. Nathan},
	title = {Discovering governing equations from data by sparse identification of nonlinear dynamical systems},
	volume = {113},
	number = {15},
	pages = {3932--3937},
	year = {2016},
	doi = {10.1073/pnas.1517384113},
	publisher = {National Academy of Sciences},
	abstract = {Understanding dynamic constraints and balances in nature has facilitated rapid development of knowledge and enabled technology, including aircraft, combustion engines, satellites, and electrical power. This work develops a novel framework to discover governing equations underlying a dynamical system simply from data measurements, leveraging advances in sparsity techniques and machine learning. The resulting models are parsimonious, balancing model complexity with descriptive ability while avoiding overfitting. There are many critical data-driven problems, such as understanding cognition from neural recordings, inferring climate patterns, determining stability of financial markets, predicting and suppressing the spread of disease, and controlling turbulence for greener transportation and energy. With abundant data and elusive laws, data-driven discovery of dynamics will continue to play an important role in these efforts.Extracting governing equations from data is a central challenge in many diverse areas of science and engineering. Data are abundant whereas models often remain elusive, as in climate science, neuroscience, ecology, finance, and epidemiology, to name only a few examples. In this work, we combine sparsity-promoting techniques and machine learning with nonlinear dynamical systems to discover governing equations from noisy measurement data. The only assumption about the structure of the model is that there are only a few important terms that govern the dynamics, so that the equations are sparse in the space of possible functions; this assumption holds for many physical systems in an appropriate basis. In particular, we use sparse regression to determine the fewest terms in the dynamic governing equations required to accurately represent the data. This results in parsimonious models that balance accuracy with model complexity to avoid overfitting. We demonstrate the algorithm on a wide range of problems, from simple canonical systems, including linear and nonlinear oscillators and the chaotic Lorenz system, to the fluid vortex shedding behind an obstacle. The fluid example illustrates the ability of this method to discover the underlying dynamics of a system that took experts in the community nearly 30 years to resolve. We also show that this method generalizes to parameterized systems and systems that are time-varying or have external forcing.},
	issn = {0027-8424},
	URL = {https://www.pnas.org/content/113/15/3932},
	eprint = {https://www.pnas.org/content/113/15/3932.full.pdf},
	journal = {Proceedings of the National Academy of Sciences}
}

@Article{Yan2020,
author={Yan, Li
and Zhang, Hai-Tao
and Goncalves, Jorge
and Xiao, Yang
and Wang, Maolin
and Guo, Yuqi
and Sun, Chuan
and Tang, Xiuchuan
and Jing, Liang
and Zhang, Mingyang
and Huang, Xiang
and Xiao, Ying
and Cao, Haosen
and Chen, Yanyan
and Ren, Tongxin
and Wang, Fang
and Xiao, Yaru
and Huang, Sufang
and Tan, Xi
and Huang, Niannian
and Jiao, Bo
and Cheng, Cheng
and Zhang, Yong
and Luo, Ailin
and Mombaerts, Laurent
and Jin, Junyang
and Cao, Zhiguo
and Li, Shusheng
and Xu, Hui
and Yuan, Ye},
title={An interpretable mortality prediction model for COVID-19 patients},
journal={Nature Machine Intelligence},
year={2020},
month={May},
day={01},
volume={2},
number={5},
pages={283-288},
abstract={The sudden increase in COVID-19 cases is putting high pressure on healthcare services worldwide. At this stage, fast, accurate and early clinical assessment of the disease severity is vital. To support decision making and logistical planning in healthcare systems, this study leverages a database of blood samples from 485 infected patients in the region of Wuhan, China, to identify crucial predictive biomarkers of disease mortality. For this purpose, machine learning tools selected three biomarkers that predict the mortality of individual patients more than 10 days in advance with more than 90{\%} accuracy: lactic dehydrogenase (LDH), lymphocyte and high-sensitivity C-reactive protein (hs-CRP). In particular, relatively high levels of LDH alone seem to play a crucial role in distinguishing the vast majority of cases that require immediate medical attention. This finding is consistent with current medical knowledge that high LDH levels are associated with tissue breakdown occurring in various diseases, including pulmonary disorders such as pneumonia. Overall, this Article suggests a simple and operable decision rule to quickly predict patients at the highest risk, allowing them to be prioritized and potentially reducing the mortality rate.},
issn={2522-5839},
doi={10.1038/s42256-020-0180-7},
url={https://doi.org/10.1038/s42256-020-0180-7}
}

@Article{Vinyals2019,
author={Vinyals, Oriol
and Babuschkin, Igor
and Czarnecki, Wojciech M.
and Mathieu, Micha{\"e}l
and Dudzik, Andrew
and Chung, Junyoung
and Choi, David H.
and Powell, Richard
and Ewalds, Timo
and Georgiev, Petko
and Oh, Junhyuk
and Horgan, Dan
and Kroiss, Manuel
and Danihelka, Ivo
and Huang, Aja
and Sifre, Laurent
and Cai, Trevor
and Agapiou, John P.
and Jaderberg, Max
and Vezhnevets, Alexander S.
and Leblond, R{\'e}mi
and Pohlen, Tobias
and Dalibard, Valentin
and Budden, David
and Sulsky, Yury
and Molloy, James
and Paine, Tom L.
and Gulcehre, Caglar
and Wang, Ziyu
and Pfaff, Tobias
and Wu, Yuhuai
and Ring, Roman
and Yogatama, Dani
and W{\"u}nsch, Dario
and McKinney, Katrina
and Smith, Oliver
and Schaul, Tom
and Lillicrap, Timothy
and Kavukcuoglu, Koray
and Hassabis, Demis
and Apps, Chris
and Silver, David},
title={Grandmaster level in StarCraft II using multi-agent reinforcement learning},
journal={Nature},
year={2019},
month={Nov},
day={01},
volume={575},
number={7782},
pages={350-354},
abstract={Many real-world applications require artificial agents to compete and coordinate with other agents in complex environments. As a stepping stone to this goal, the domain of StarCraft has emerged as an important challenge for artificial intelligence research, owing to its iconic and enduring status among the most difficult professional esports and its relevance to the real world in terms of its raw complexity and multi-agent challenges. Over the course of a decade and numerous competitions1--3, the strongest agents have simplified important aspects of the game, utilized superhuman capabilities, or employed hand-crafted sub-systems4. Despite these advantages, no previous agent has come close to matching the overall skill of top StarCraft players. We chose to address the challenge of StarCraft using general-purpose learning methods that are in principle applicable to other complex domains: a multi-agent reinforcement learning algorithm that uses data from both human and agent games within a diverse league of continually adapting strategies and counter-strategies, each represented by deep neural networks5,6. We evaluated our agent, AlphaStar, in the full game of StarCraft II, through a series of online games against human players. AlphaStar was rated at Grandmaster level for all three StarCraft races and above 99.8{\%} of officially ranked human players.},
issn={1476-4687},
doi={10.1038/s41586-019-1724-z},
url={https://doi.org/10.1038/s41586-019-1724-z}
}

@Article{Silver2016,
author={Silver, David
and Huang, Aja
and Maddison, Chris J.
and Guez, Arthur
and Sifre, Laurent
and van den Driessche, George
and Schrittwieser, Julian
and Antonoglou, Ioannis
and Panneershelvam, Veda
and Lanctot, Marc
and Dieleman, Sander
and Grewe, Dominik
and Nham, John
and Kalchbrenner, Nal
and Sutskever, Ilya
and Lillicrap, Timothy
and Leach, Madeleine
and Kavukcuoglu, Koray
and Graepel, Thore
and Hassabis, Demis},
title={Mastering the game of Go with deep neural networks and tree search},
journal={Nature},
year={2016},
month={Jan},
day={01},
volume={529},
number={7587},
pages={484-489},
abstract={The game of Go has long been viewed as the most challenging of classic games for artificial intelligence owing to its enormous search space and the difficulty of evaluating board positions and moves. Here we introduce a new approach to computer Go that uses `value networks' to evaluate board positions and `policy networks' to select moves. These deep neural networks are trained by a novel combination of supervised learning from human expert games, and reinforcement learning from games of self-play. Without any lookahead search, the neural networks play Go at the level of state-of-the-art Monte Carlo tree search programs that simulate thousands of random games of self-play. We also introduce a new search algorithm that combines Monte Carlo simulation with value and policy networks. Using this search algorithm, our program AlphaGo achieved a 99.8{\%} winning rate against other Go programs, and defeated the human European Go champion by 5 games to 0. This is the first time that a computer program has defeated a human professional player in the full-sized game of Go, a feat previously thought to be at least a decade away.},
issn={1476-4687},
doi={10.1038/nature16961},
url={https://doi.org/10.1038/nature16961}
}

@manual{Haydin2018,
      title  = "How Machine Learning Algorithms Make Self-Driving Cars a Reality",
      author = "Victor Haydin",
      note   = "https://www.intellias.com/how-machine-learning-algorithms-make-self-driving-cars-a-reality/",
      year   = "2018 (accessed 2020)"
    }

@article{DBLP:journals/corr/cs-NA-0111063,
  author    = {W. Chen},
  title     = {New {RBF} collocation methods and kernel {RBF} with applications},
  journal   = {CoRR},
  volume    = {cs.NA/0111063},
  year      = {2001},
  url       = {https://arxiv.org/abs/cs/0111063},
  timestamp = {Fri, 10 Jan 2020 12:58:18 +0100},
  biburl    = {https://dblp.org/rec/journals/corr/cs-NA-0111063.bib},
  bibsource = {dblp computer science bibliography, https://dblp.org}
}
@article{PIRET20124662,
title = "The orthogonal gradients method: A radial basis functions method for solving partial differential equations on arbitrary surfaces",
journal = "Journal of Computational Physics",
volume = "231",
number = "14",
pages = "4662 - 4675",
year = "2012",
issn = "0021-9991",
doi = "https://doi.org/10.1016/j.jcp.2012.03.007",
url = "http://www.sciencedirect.com/science/article/pii/S0021999112001477",
author = "C\'{e}cile Piret",
keywords = "Radial basis functions, RBF, Closest point method, Implicit surfaces, Level set method, Orthogonal gradients method, OGr method",
abstract = "Much work has been done on reconstructing arbitrary surfaces using the radial basis function (RBF) method, but one can hardly find any work done on the use of RBFs to solve partial differential equations (PDEs) on arbitrary surfaces. In this paper, we investigate methods to solve PDEs on arbitrary stationary surfaces embedded in R3 using the RBF method. We present three RBF-based methods that easily discretize surface differential operators. We take advantage of the meshfree character of RBFs, which give us a high accuracy and the flexibility to represent the most complex geometries in any dimension. Two out of the three methods, which we call the orthogonal gradients (OGr) methods are the result of our work and are hereby presented for the first time."
}

@TECHREPORT{IMM2005-03600,
    author       = "M. J. Powell",
    title        = "Five Lectures on Radial Basis Functions",
    year         = "2005",
    number       = "",
    series       = "IMM-Techncial Report-2005-03",
    institution  = "Informatics and Mathematical Modelling, Technical University of Denmark, {DTU}",
    address      = "Richard Petersens Plads, Building 321, {DK-}2800 Kgs. Lyngby",
    type         = "",
    url          = "http://www2.compute.dtu.dk/pubdb/pubs/3600-full.html",
    abstract     = "Professor Mike J. D. Powell spent three weeks at {IMM} in November - December 2004. During the visit he gave five lectures on radial basis functions. These notes are a TeXified version of his hand-outs, made by Hans Bruun Nielsen, {IMM}."
}

@article{Fuselier2013AHK,
  title={A High-Order Kernel Method for Diffusion and Reaction-Diffusion Equations on Surfaces},
  author={Edward J. Fuselier and Grady B. Wright},
  journal={Journal of Scientific Computing},
  year={2013},
  volume={56},
  pages={535-565}
}
@article{Yvonnet2004,
author = {Yvonnet, J. and Ryckelynck, D. and Lorong, P. and Chinesta, F.},
title = {A new extension of the natural element method for non-convex and discontinuous problems: the constrained natural element method (C-NEM)},
journal = {International Journal for Numerical Methods in Engineering},
volume = {60},
number = {8},
pages = {1451-1474},
keywords = {C-NEM, non-convex domains, visibility criterion, constrained Voronoi diagram},
doi = {10.1002/nme.1016},
url = {https://onlinelibrary.wiley.com/doi/abs/10.1002/nme.1016},
eprint = {https://onlinelibrary.wiley.com/doi/pdf/10.1002/nme.1016},
abstract = {Abstract In this paper a new extension of the mesh-free natural element method (NEM) is presented. In this approach, coined as constrained natural element method (C-NEM), a visibility criterion is introduced to select natural neighbours in the computation of the shape functions. The computation of these shape functions is based on a modified, constrained Voronoi diagram. With this technique, some difficulties inherent to this method in non-convex domains are avoided and the analysis of problems involving cracks or discontinuities are now easily performed. As the NEM satisfies the Kronecker delta property, the imposition of essential boundary conditions is trivial, unlike other mesh-free methods. The C-NEM technique provides a description of integration cells that allows the use of the stabilized conforming nodal integration (SCNI) scheme instead of Gauss integration to enhance computational efficiency and accuracy. Two numerical examples in elastostatics are reported to evaluate the potential of the proposed technique in highly non-convex geometries, like a crack, through which the solution becomes discontinuous. Copyright © 2004 John Wiley \& Sons, Ltd.},
year = {2004}
}
@Article{Atluri1998,
author={Atluri, S. N.
and Zhu, T.},
title={A new Meshless Local Petrov-Galerkin (MLPG) approach in computational mechanics},
journal={Computational Mechanics},
year={1998},
month={Aug},
day={01},
volume={22},
number={2},
pages={117-127},
abstract={A local symmetric weak form (LSWF) for linear potential problems is developed, and a truly meshless method, based on the LSWF and the moving least squares approximation, is presented for solving potential problems with high accuracy. The essential boundary conditions in the present formulation are imposed by a penalty method. The present method does not need a ``finite element mesh'', either for purposes of interpolation of the solution variables, or for the integration of the ``energy''. All integrals can be easily evaluated over regularly shaped domains (in general, spheres in three-dimensional problems) and their boundaries. No post-smoothing technique is required for computing the derivatives of the unknown variable, since the original solution, using the moving least squares approximation, is already smooth enough. Several numerical examples are presented in the paper. In the example problems dealing with Laplace {\&} Poisson's equations, high rates of convergence with mesh refinement for the Sobolev norms ||{\textperiodcentered}||0 and ||{\textperiodcentered}||1 have been found, and the values of the unknown variable and its derivatives are quite accurate. In essence, the present meshless method based on the LSWF is found to be a simple, efficient, and attractive method with a great potential in engineering applications.},
issn={1432-0924},
doi={10.1007/s004660050346},
url={https://doi.org/10.1007/s004660050346}
}

@article{ONATE1996,
author = {Onate, E. and Idelsohn, S. and Zienkiwicz, O. C. and TAYLOR, R. L.},
title = {A FINITE POINT METHOD IN COMPUTATIONAL MECHANICS. APPLICATIONS TO CONVECTIVE TRANSPORT AND FLUID FLOW},
journal = {International Journal for Numerical Methods in Engineering},
volume = {39},
number = {22},
pages = {3839-3866},
keywords = {least squares, finite point method, mesh free techniques},
doi = {10.1002/(SICI)1097-0207(19961130)39:22<3839::AID-NME27>3.0.CO;2-R},
url = {https://onlinelibrary.wiley.com/doi/abs/10.1002/%28SICI%291097-0207%2819961130%2939%3A22%3C3839%3A%3AAID-NME27%3E3.0.CO%3B2-R},
eprint = {https://onlinelibrary.wiley.com/doi/pdf/10.1002/%28SICI%291097-0207%2819961130%2939%3A22%3C3839%3A%3AAID-NME27%3E3.0.CO%3B2-R},
abstract = {Abstract The paper presents a fully meshless procedure fo solving partial differential equations. The approach termed generically the ‘finite point method’ is based on a weighted least square interpolation of point data and point collocation for evaluating the approximation integrals. Some examples showing the accuracy of the method for solution of adjoint and non-self adjoint equations typical of convective-diffusive transport and also to the analysis of compressible fluid mechanics problem are presented.},
year = {1996}
}

@article{LISZKA198083,
title = "The finite difference method at arbitrary irregular grids and its application in applied mechanics",
journal = "Computers \& Structures",
volume = "11",
number = "1",
pages = "83 - 95",
year = "1980",
note = "Special Issue-Computational Methods in Nonlinear Mechanics",
issn = "0045-7949",
doi = "https://doi.org/10.1016/0045-7949(80)90149-2",
url = "http://www.sciencedirect.com/science/article/pii/0045794980901492",
author = "T. Liszka and J. Orkisz",
abstract = "Presented modification of the FDM enables local condensation of the mesh and easy discretisation of the boundary conditions in the case of an arbitrary shape of the domain. As a result an essential reduction of the required number of nodal points may usually be achieved. Thus the FDM can be competitive in some fields to the finite element method. Several troubles arising from the use of an irregulr mesh have been solved, the most important being elimination of singular or ill-conditioned stars, and a successful way of automatic discretization of boundary conditions was proposed. As a result of this complete automatization of the FDM has been reached. Many problems connected with the theory of the method proposed and its computer implementation have been discussed. FIDAM code—a system of computer programs designed for the solution of two-dimensional, linear and nonlinear, elliptic problems and three-dimensional parabolic problems is presented. Problems to be solved should be layed down in a local formulation as a set of partial derivative equations of the second order with boundary conditions not exceeding the same order. Various particular problems of applied mechanics and physics such as: torsion of bars, plane elasticity problems, deflections of plates and membranes, fluid now and temperature distribution have been solved using the FIDAM procedures. For the solution of nonlinear problems various iterative methods have been adopted—with special attention payed to the self-correcting method where the Newton-Raphson and incremental procedures are particular cases. The present version of the method allows fully automatic calculations to be carried out as in advanced programs of the finite element method and may be preferred in non-linear, optimization and time-dependent problems."
}

@article{Liu1995,
author = {Liu, Wing Kam and Jun, Sukky and Zhang, Yi Fei},
title = {Reproducing kernel particle methods},
journal = {International Journal for Numerical Methods in Fluids},
volume = {20},
number = {8‐9},
pages = {1081-1106},
keywords = {multiple scale decomposition, correction function, multi-resolution analysis, reproducing kernel function, wavelet, mesh- (or grid-) free particle methods},
doi = {10.1002/fld.1650200824},
url = {https://onlinelibrary.wiley.com/doi/abs/10.1002/fld.1650200824},
eprint = {https://onlinelibrary.wiley.com/doi/pdf/10.1002/fld.1650200824},
abstract = {Abstract A new continuous reproducing kernel interpolation function which explores the attractive features of the flexible time-frequency and space-wave number localization of a window function is developed. This method is motivated by the theory of wavelets and also has the desirable attributes of the recently proposed smooth particle hydrodynamics (SPH) methods, moving least squares methods (MLSM), diffuse element methods (DEM) and element-free Galerkin methods (EFGM). The proposed method maintains the advantages of the free Lagrange or SPH methods; however, because of the addition of a correction function, it gives much more accurate results. Therefore it is called the reproducing kernel particle method (RKPM). In computer implementation RKPM is shown to be more efficient than DEM and EFGM. Moreover, if the window function is C∞, the solution and its derivatives are also C∞ in the entire domain. Theoretical analysis and numerical experiments on the 1D diffusion equation reveal the stability conditions and the effect of the dilation parameter on the unusually high convergence rates of the proposed method. Two-dimensional examples of advection-diffusion equations and compressible Euler equations are also presented together with 2D multiple-scale decompositions.},
year = {1995}
}

@Article{Atluri2000,
author={Atluri, S. N.
and Zhu, T.-L.},
title={The meshless local Petrov-Galerkin (MLPG) approach for solving problems in elasto-statics},
journal={Computational Mechanics},
year={2000},
month={Mar},
day={01},
volume={25},
number={2},
pages={169-179},
abstract={{\enspace}The meshless local Petrov-Galerkin (MLPG) approach is an effective method for solving boundary value problems, using a local symmetric weak form and shape functions from the moving least squares approximation. In the present paper, the MLPG method for solving problems in elasto-statics is developed and numerically implemented. The present method is a truly meshless method, as it does not need a ``finite element mesh'', either for purposes of interpolation of the solution variables, or for the integration of the energy. All integrals in the formulation can be easily evaluated over regularly shaped domains (in general, spheres in three-dimensional problems) and their boundaries. The essential boundary conditions in the present formulation are imposed by a penalty method, as the essential boundary conditions can not be enforced directly when the non-interpolative moving least squares approximation is used. Several numerical examples are presented to illustrate the implementation and performance of the present MLPG method. The numerical examples show that the present MLPG approach does not exhibit any volumetric locking for nearly incompressible materials, and that high rates of convergence with mesh refinement for the displacement and energy norms are achievable. No post-processing procedure is required to compute the strain and stress, since the original solution from the present method, using the moving least squares approximation, is already smooth enough.},
issn={1432-0924},
doi={10.1007/s004660050467},
url={https://doi.org/10.1007/s004660050467}
}

@article{Gingold1977,
    author = {Gingold, R. A. and Monaghan, J. J.},
    title = "{Smoothed particle hydrodynamics: theory and application to non-spherical stars}",
    journal = {Monthly Notices of the Royal Astronomical Society},
    volume = {181},
    number = {3},
    pages = {375-389},
    year = {1977},
    month = {12},
    abstract = "{A new hydrodynamic code applicable to a space of an arbitrary number of dimensions is discussed and applied to a variety of polytropic stellar models. The principal feature of the method is the use of statistical techniques to recover analytical expressions for the physical variables from a known distribution of fluid elements. The equations of motion take the form of Newtonian equations for particles. Starting with a non-axisymmetric distribution of approximately 80 particles in three dimensions, the method is found to reproduce the structure of uniformly rotating and magnetic polytropes to within a few per cent. The method may be easily extended to deal with more complicated physical models.}",
    issn = {0035-8711},
    doi = {10.1093/mnras/181.3.375},
    url = {https://doi.org/10.1093/mnras/181.3.375},
    eprint = {https://academic.oup.com/mnras/article-pdf/181/3/375/3104055/mnras181-0375.pdf},
}

@Article{Zhu1998,
author={Zhu, T.
and Zhang, J.-D.
and Atluri, S. N.},
title={A local boundary integral equation (LBIE) method in computational mechanics, and a meshless discretization approach},
journal={Computational Mechanics},
year={1998},
month={Apr},
day={01},
volume={21},
number={3},
pages={223-235},
abstract={The Galerkin finite element method (GFEM) owes its popularity to the local nature of nodal basis functions, i.e., the nodal basis function, when viewed globally, is non-zero only over a patch of elements connecting the node in question to its immediately neighboring nodes. The boundary element method (BEM), on the other hand, reduces the dimensionality of the problem by one, through involving the trial functions and their derivatives, only in the integrals over the global boundary of the domain; whereas, the GFEM involves the integration of the ``energy'' corresponding to the trial function over a patch of elements immediately surrounding the node. The GFEM leads to banded, sparse and symmetric matrices; the BEM based on the global boundary integral equation (GBIE) leads to full and unsymmetrical matrices. Because of the seemingly insurmountable difficulties associated with the automatic generation of element-meshes in GFEM, especially for 3-D problems, there has been a considerable interest in element free Galerkin methods (EFGM) in recent literature. However, the EFGMs still involve domain integrals over shadow elements and lead to difficulties in enforcing essential boundary conditions and in treating nonlinear problems.},
issn={1432-0924},
doi={10.1007/s004660050297},
url={https://doi.org/10.1007/s004660050297}
}

@article{BELYTSCHKO2000385,
title = "Element-free galerkin methods for dynamic fracture in concrete",
journal = "Computer Methods in Applied Mechanics and Engineering",
volume = "187",
number = "3",
pages = "385 - 399",
year = "2000",
issn = "0045-7825",
doi = "https://doi.org/10.1016/S0045-7825(00)80002-X",
url = "http://www.sciencedirect.com/science/article/pii/S004578250080002X",
author = "T. Belytschko and D. Organ and C. Gerlach",
abstract = "Mixed-mode dynamic crack propagation in concrete is studied using the element-free Galerkin (EFG) method. The EFG methodology allows for arbitrary crack growth in terms of direction and speed. A fracture process zone (FPZ) model is used for fracture in concrete. Comparisons are made to mode I and mixed-mode experiments on concrete three-point-bend specimens. Good agreement is acheived with the experimental results."
}
@article{BELYTSCHKO19963,
title = "Meshless methods: An overview and recent developments",
journal = "Computer Methods in Applied Mechanics and Engineering",
volume = "139",
number = "1",
pages = "3 - 47",
year = "1996",
issn = "0045-7825",
doi = "https://doi.org/10.1016/S0045-7825(96)01078-X",
url = "http://www.sciencedirect.com/science/article/pii/S004578259601078X",
author = "T. Belytschko and Y. Krongauz and D. Organ and M. Fleming and P. Krysl",
abstract = "Meshless approximations based on moving least-squares, kernels, and partitions of unity are examined. It is shown that the three methods are in most cases identical except for the important fact that partitions of unity enable p-adaptivity to be achieved. Methods for constructing discontinuous approximations and approximations with discontinuous derivatives are also described. Next, several issues in implementation are reviewed: discretization (collocation and Galerkin), quadrature in Galerkin and fast ways of constructing consistent moving least-square approximations. The paper concludes with some sample calculations."
}

@article{BABUSKA1997,
author = {Babuska, I. and Melenk, J. M.},
title = {THE PARTITION OF UNITY METHOD},
journal = {International Journal for Numerical Methods in Engineering},
volume = {40},
number = {4},
pages = {727-758},
keywords = {finite element method, meshless finite element method, finite element methods for highly oscillatory solutions},
doi = {10.1002/(SICI)1097-0207(19970228)40:4<727::AID-NME86>3.0.CO;2-N},
url = {https://onlinelibrary.wiley.com/doi/abs/10.1002/%28SICI%291097-0207%2819970228%2940%3A4%3C727%3A%3AAID-NME86%3E3.0.CO%3B2-N},
eprint = {https://onlinelibrary.wiley.com/doi/pdf/10.1002/%28SICI%291097-0207%2819970228%2940%3A4%3C727%3A%3AAID-NME86%3E3.0.CO%3B2-N},
abstract = {Abstract A new finite element method is presented that features the ability to include in the finite element space knowledge about the partial differential equation being solved. This new method can therefore be more efficient than the usual finite element methods. An additional feature of the partition-of-unity method is that finite element spaces of any desired regularity can be constructed very easily. This paper includes a convergence proof of this method and illustrates its efficiency by an application to the Helmholtz equation for high wave numbers. The basic estimates for a posteriori error estimation for this new method are also proved. © 1997 by John Wiley \& Sons, Ltd.},
year = {1997}
}

@article{mirzadeh2011,
  title={A second-order discretization of the nonlinear {Poisson--Boltzmann} equation over irregular geometries using non-graded adaptive {Cartesian} grids},
  author={M. Mirzadeh and M. Theillard and F. Gibou},
  journal={J. Comput. Phys.},
  volume={230},
  pages={2125--2140},
  year={2011},
  publisher={Elsevier}
}

@Article{Theillard2013,
author={Theillard, Maxime
and Rycroft, Chris H.
and Gibou, Fr{\'e}d{\'e}ric},
title={A Multigrid Method on Non-Graded Adaptive Octree and Quadtree Cartesian Grids},
journal={Journal of Scientific Computing},
year={2013},
month={Apr},
day={01},
volume={55},
number={1},
pages={1-15},
abstract={In order to develop efficient numerical methods for solving elliptic and parabolic problems where Dirichlet boundary conditions are imposed on irregular domains, Chen et al. (J. Sci. Comput. 31(1):19--60, 2007) presented a methodology that produces second-order accurate solutions with second-order gradients on non-graded quadtree and octree data structures. These data structures significantly reduce the number of computational nodes while still allowing for the resolution of small length scales. In this paper, we present a multigrid solver for this framework and present numerical results in two and three spatial dimensions that demonstrate that the computational time scales linearly with the number of nodes, producing a very efficient solver for elliptic and parabolic problems with multiple length scales.},
issn={1573-7691},
doi={10.1007/s10915-012-9619-2},
url={https://doi.org/10.1007/s10915-012-9619-2}
}
@article{BELLOTTI2019266,
title = "A coupled level-set and reference map method for interface representation with applications to two-phase flows simulation",
journal = "Journal of Computational Physics",
volume = "392",
pages = "266 - 290",
year = "2019",
issn = "0021-9991",
doi = "https://doi.org/10.1016/j.jcp.2019.05.003",
url = "http://www.sciencedirect.com/science/article/pii/S0021999119303286",
author = "Thomas Bellotti and Maxime Theillard",
keywords = "Interface problems, Level set, Reference map, Reinitialization, Navier-Stokes, Two-phase flows",
abstract = "We present a novel hybrid methodology combining the reference map theory with the level-set method for tracking moving interfaces. Instead of directly advecting the level-set function, we track the reference map, which maps the deformed state into the original one. We then reconstruct the deformed level-set function from this mapping and the original geometry. Because of the better smoothness of the reference map and the reduced impact of the reinitialization, this new approach grants higher precision. This results in significant improvements regarding interface location precision and mass conservation, especially in situations involving small deformations. Moreover, the implementation of the method is straightforward since it is based on the same numerical techniques used for the level-set method. Our coupled method is extensively validated in the case of externally generated velocity fields and incorporated into our previously introduced two-phase incompressible Navier-Stokes flow solver."
}

@article{mirzadeh_theillard_helgadöttir_boy_gibou_2013, title={An Adaptive, Finite Difference Solver for the Nonlinear Poisson-Boltzmann Equation with Applications to Biomolecular Computations}, volume={13}, DOI={10.4208/cicp.290711.181011s}, number={1}, journal={Communications in Computational Physics}, publisher={Cambridge University Press}, author={Mirzadeh, Mohammad and Theillard, Maxime and Helgadöttir, Asdís and Boy, David and Gibou, Frédéric}, year={2013}, pages={150–173}}
@article{THEILLARD2013430,
title = "A second-order sharp numerical method for solving the linear elasticity equations on irregular domains and adaptive grids – Application to shape optimization",
journal = "Journal of Computational Physics",
volume = "233",
pages = "430 - 448",
year = "2013",
issn = "0021-9991",
doi = "https://doi.org/10.1016/j.jcp.2012.09.002",
url = "http://www.sciencedirect.com/science/article/pii/S0021999112005207",
author = "Maxime Theillard and Landry Fokoua Djodom and Jean-Léopold Vié and Frédéric Gibou",
keywords = "Second-order discretization, Linear elasticity, Hybrid finite volume/finite difference, Level-set, Irregular domains, Octree data structure, Quadtree/octree data structure, Non-graded adaptive grid",
abstract = "We present a numerical method for solving the equations of linear elasticity on irregular domains in two and three spatial dimensions. We combine a finite volume and a finite difference approaches to derive discretizations that produce second-order accurate solutions in the L∞-norm. Our discretization is ‘sharp’ in the sense that the physical boundary conditions (mixed Dirichlet/Neumann-type) are imposed at the interface and the solution is computed inside the irregular domain only, without the need of smearing the solution across the interface. The irregular domain is represented implicitly using a level-set function so that this approach is applicable to free moving boundary problems; we provide a simple example of shape optimization to illustrate this capability. In addition, we provide an extension of our method to the case of adaptive meshes in both two and three spatial dimensions: we use non-graded quadtree (2D) and octree (3D) data structures to represent the grid that is automatically refined near the irregular domain’s boundary. This extension to quadtree/octree grids produces second-order accurate solutions albeit non-symmetric linear systems, due to the node-based sampling nature of the approach. However, the linear system can be solved with simple linear solvers; in this work we use the BICGSTAB algorithm."
}

@article{CLERETDELANGAVANT2017271,
title = "Level-set simulations of soluble surfactant driven flows",
journal = "Journal of Computational Physics",
volume = "348",
pages = "271 - 297",
year = "2017",
issn = "0021-9991",
doi = "https://doi.org/10.1016/j.jcp.2017.07.003",
url = "http://www.sciencedirect.com/science/article/pii/S0021999117305053",
author = "Charles {Cleret de Langavant} and Arthur Guittet and Maxime Theillard and Fernando Temprano-Coleto and Frédéric Gibou",
keywords = "Navier–Stokes, Incompressible, Soluble surfactants, Surfactant driven flows, Marangoni forces, Quad/Octrees, Adaptive mesh refinement, Stable projection method",
abstract = "We present an approach to simulate the diffusion, advection and adsorption–desorption of a material quantity defined on an interface in two and three spatial dimensions. We use a level-set approach to capture the interface motion and a Quad/Octree data structure to efficiently solve the equations describing the underlying physics. Coupling with a Navier–Stokes solver enables the study of the effect of soluble surfactants that locally modify the parameters of surface tension on different types of flows. The method is tested on several benchmarks and applied to three typical examples of flows in the presence of surfactant: a bubble in a shear flow, the well-known phenomenon of tears of wine, and the Landau–Levich coating problem."
}

@article{THEILLARD201991,
title = "Sharp numerical simulation of incompressible two-phase flows",
journal = "Journal of Computational Physics",
volume = "391",
pages = "91 - 118",
year = "2019",
issn = "0021-9991",
doi = "https://doi.org/10.1016/j.jcp.2019.04.024",
url = "http://www.sciencedirect.com/science/article/pii/S0021999119302578",
author = "Maxime Theillard and Frédéric Gibou and David Saintillan",
keywords = "Two-phase flows, Incompressible Navier-Stokes, Projection method, Sharp jump conditions, Level-set, MAC grid",
abstract = "We present a numerical method for simulating incompressible immiscible fluids, in two and three spatial dimensions. It is constructed as a modified pressure correction projection method on adaptive non-graded Oc/Quadtree Cartesian grids, using the level-set framework to capture the moving interface between the two fluids. The sharp treatment of the interface position, of the fluid parameter discontinuities, and of the interfacial jump conditions ensures convergence in the L∞-norm. Using a novel construction for the pressure guess, we are able to alleviate the standard time step restriction incurred by capillary forces. The solver is validated numerically and employed to simulate the dynamics of physically relevant problems such as rising bubbles and viscous droplets in electric fields."
}

@article{DU20033933,
title = "Voronoi-based finite volume methods, optimal Voronoi meshes, and PDEs on the sphere",
journal = "Computer Methods in Applied Mechanics and Engineering",
volume = "192",
number = "35",
pages = "3933 - 3957",
year = "2003",
issn = "0045-7825",
doi = "https://doi.org/10.1016/S0045-7825(03)00394-3",
url = "http://www.sciencedirect.com/science/article/pii/S0045782503003943",
author = "Qiang Du and Max D. Gunzburger and Lili Ju",
abstract = "We first develop and analyze a finite volume scheme for the discretization of partial differential equations (PDEs) on the sphere; the scheme uses Voronoi tessellations of the sphere. For a model convection–diffusion problem, the finite volume scheme is shown to produce first-order accurate approximations with respect to a mesh-dependent discrete first-derivative norm. Then, we introduce the notion of constrained centroidal Voronoi tessellations (CCVTs) of the sphere; these are special Voronoi tessellation of the sphere for which the generators of the Voronoi cells are also the constrained centers of mass, with respect to a prescribed density function, of the cells. After discussing an algorithm for determining CCVT meshes on the sphere, we discuss and illustrate several desirable properties possessed by these meshes. In particular, it is shown that CCVT meshes define very high-quality uniform and non-uniform meshes on the sphere. Finally, we discuss, through some computational experiments, the performance of the CCVT meshes used in conjunction with the finite volume scheme for the solution of simple model PDEs on the sphere. The experiments show, for example, that the CCVT based finite volume approximations are second-order accurate if errors are measured in discrete L2 norms."
}

@article{Guittet2015ASP,
  title={A stable projection method for the incompressible Navier-Stokes equations on arbitrary geometries and adaptive Quad/Octrees},
  author={Arthur Guittet and Maxime Theillard and Fr{\'e}d{\'e}ric Gibou},
  journal={J. Comput. Phys.},
  year={2015},
  volume={292},
  pages={215-238}
}

@article{Sukumar1998,
author = {Sukumar, N. and Moran, B. and Belytschko, T.},
title = {The natural element method in solid mechanics},
journal = {International Journal for Numerical Methods in Engineering},
volume = {43},
number = {5},
pages = {839-887},
keywords = {natural neighbour interpolation, natural element method, 1st- and 2nd-order Voronoi diagrams, Delaunay triangle, elastostatics},
doi = {10.1002/(SICI)1097-0207(19981115)43:5<839::AID-NME423>3.0.CO;2-R},
url = {https://onlinelibrary.wiley.com/doi/abs/10.1002/%28SICI%291097-0207%2819981115%2943%3A5%3C839%3A%3AAID-NME423%3E3.0.CO%3B2-R},
eprint = {https://onlinelibrary.wiley.com/doi/pdf/10.1002/%28SICI%291097-0207%2819981115%2943%3A5%3C839%3A%3AAID-NME423%3E3.0.CO%3B2-R},
abstract = {Abstract The application of the Natural Element Method (NEM) to boundary value problems in two-dimensional small displacement elastostatics is presented. The discrete model of the domain Ω consists of a set of distinct nodes N, and a polygonal description of the boundary ∂Ω. In the Natural Element Method, the trial and test functions are constructed using natural neighbour interpolants. These interpolants are based on the Voronoi tessellation of the set of nodes N. The interpolants are smooth (C∞) everywhere, except at the nodes where they are C0. In one-dimension, NEM is identical to linear finite elements. The NEM interpolant is strictly linear between adjacent nodes on the boundary of the convex hull, which facilitates imposition of essential boundary conditions. A methodology to model material discontinuities and non-convex bodies (cracks) using NEM is also described. A standard displacement-based Galerkin procedure is used to obtain the discrete system of linear equations. Application of NEM to various problems in solid mechanics, which include, the patch test, gradient problems, bimaterial interface, and a static crack problem are presented. Excellent agreement with exact (analytical) solutions is obtained, which exemplifies the accuracy and robustness of NEM and suggests its potential application in the context of other classes of problems—crack growth, plates, and large deformations to name a few. © 1998 John Wiley \& Sons, Ltd.},
year = {1998}
}

@book{Gonzalez2009,
author="Gonz{\'a}lez, David
and Cueto, El{\'i}as
and Doblar{\'e}, Manuel",
editor="Ferreira, A. J. M.
and Kansa, E. J.
and Fasshauer, G. E.
and Leit{\~a}o, V. M. A.",
title="Towards an Isogeometric Meshless Natural Element Method",
bookTitle="Progress on Meshless Methods",
year="2009",
publisher="Springer Netherlands",
address="Dordrecht",
pages="237--257",
abstract="The problem of generalizing the Natural Element Method in terms of higher-order consistency and continuity is addressed here with several possible solutions. In this work we review some of them. First, a study of the possible benefits of enriching the interpolation using the Partition of Unity paradigm is considered. Different enrichments were tested leading to different reproducing properties. Another possible solution to the problem is the use of iterated Voronoi diagram interpolants, due to G. Farin. This last solution is done by means of the de Boor's algorithm, the same employed to obtain B-splines by linear combinations of linear interpolants in one dimension. We propose another form of B-spline-like interpolants that employs the de Boor's algorithm, but with a simpler structure. In order to obtain a smooth interpolant, a review of the Hiyoshi-Sugihara interpolant is also made. In this case, although high-order smoothness can be achieved, the consistency remains linear. By employing the proposed algorithms, however, this consistency can be improved to the desired degree.",
isbn="978-1-4020-8821-6",
doi="10.1007/978-1-4020-8821-6_14",
url="https://doi.org/10.1007/978-1-4020-8821-6_14"
}

@article{Duarte1996,
author = {Duarte, C. Armando and Oden, J. Tinsley},
title = {H-p clouds—an h-p meshless method},
journal = {Numerical Methods for Partial Differential Equations},
volume = {12},
number = {6},
pages = {673-705},
doi = {10.1002/(SICI)1098-2426(199611)12:6<673::AID-NUM3>3.0.CO;2-P},
url = {https://onlinelibrary.wiley.com/doi/abs/10.1002/%28SICI%291098-2426%28199611%2912%3A6%3C673%3A%3AAID-NUM3%3E3.0.CO%3B2-P},
eprint = {https://onlinelibrary.wiley.com/doi/pdf/10.1002/%28SICI%291098-2426%28199611%2912%3A6%3C673%3A%3AAID-NUM3%3E3.0.CO%3B2-P},
abstract = {Abstract A new methodology to build discrete models of boundary-value problems is presented. The h-p cloud method is applicable to arbitrary domains and employs only a scattered set of nodes to build approximate solutions to BVPs. This new method uses radial basis functions of varying size of supports and with polynomial-reproducing properties of arbitrary order. The approximating properties of the h-p cloud functions are investigated in this article and several theorems concerning these properties are presented. Moving least squares interpolants are used to build a partition of unity on the domain of interest. These functions are then used to construct, at a very low cost, trial and test functions for Galerkin approximations. The method exhibits a very high rate of convergence and has a greater flexibility than traditional h-p finite element methods. Several numerical experiments in 1-D and 2-D are also presented. © 1996 John Wiley \& Sons, Inc.},
year = {1996}
}

@article{MAIDUY2002133,
title = "Mesh-free radial basis function network methods with domain decomposition for approximation of functions and numerical solution of Poisson's equations",
journal = "Engineering Analysis with Boundary Elements",
volume = "26",
number = "2",
pages = "133 - 156",
year = "2002",
issn = "0955-7997",
doi = "https://doi.org/10.1016/S0955-7997(01)00092-3",
url = "http://www.sciencedirect.com/science/article/pii/S0955799701000923",
author = "Nam Mai-Duy and Thanh Tran-Cong",
keywords = "Mesh-free radial basis function network methods, Domain decomposition, Boundary integral equations, Function approximation, Poisson's equation",
abstract = "This paper presents the combination of new mesh-free radial basis function network (RBFN) methods and domain decomposition (DD) technique for approximating functions and solving Poisson's equations. The RBFN method allows numerical approximation of functions and solution of partial differential equations (PDEs) without the need for a traditional ‘finite element’-type (FE) mesh while the combined RBFN–DD approach facilitates coarse-grained parallelisation of large problems. Effect of RBFN parameters on the quality of approximation of function and its derivatives is investigated and compared with the case of single domain. In solving Poisson's equations, an iterative procedure is employed to update unknown boundary conditions at interfaces. At each iteration, the interface boundary conditions are first estimated by using boundary integral equations (BIEs) and subdomain problems are then solved by using the RBFN method. Volume integrals in standard integral equation representation (IE), which usually require volume discretisation, are completely eliminated in order to preserve the mesh-free nature of RBFN methods. The numerical examples show that RBFN methods in conjunction with DD technique achieve not only a reduction of memory requirement but also a high accuracy of the solution."
}

@article{Wright2010,
author = {Wright, G. B. and Flyer, N. and Yuen, D. A.},
title = {A hybrid radial basis function–pseudospectral method for thermal convection in a 3-D spherical shell},
journal = {Geochemistry, Geophysics, Geosystems},
volume = {11},
number = {7},
pages = {},
keywords = {mantle convection, radial basis function, pseudospectral, isoviscous, spherical geometry},
doi = {10.1029/2009GC002985},
url = {https://agupubs.onlinelibrary.wiley.com/doi/abs/10.1029/2009GC002985},
eprint = {https://agupubs.onlinelibrary.wiley.com/doi/pdf/10.1029/2009GC002985},
abstract = {A novel hybrid spectral method that combines radial basis function (RBF) and Chebyshev pseudospectral methods in a “2 + 1” approach is presented for numerically simulating thermal convection in a 3-D spherical shell. This is the first study to apply RBFs to a full 3-D physical model in spherical geometry. In addition to being spectrally accurate, RBFs are not defined in terms of any surface-based coordinate system such as spherical coordinates. As a result, when used in the lateral directions, as in this study, they completely circumvent the pole issue with the further advantage that nodes can be “scattered” over the surface of a sphere. In the radial direction, Chebyshev polynomials are used, which are also spectrally accurate and provide the necessary clustering near the boundaries to resolve boundary layers. Applications of this new hybrid methodology are given to the problem of convection in the Earth's mantle, which is modeled by a Boussinesq fluid at infinite Prandtl number. To see whether this numerical technique warrants further investigation, the study limits itself to an isoviscous mantle. Benchmark comparisons are presented with other currently used mantle convection codes for Rayleigh number (Ra) 7 × 103 and 105. Results from a Ra = 106 simulation are also given. The algorithmic simplicity of the code (mostly due to RBFs) allows it to be written in less than 400 lines of MATLAB and run on a single workstation. We find that our method is very competitive with those currently used in the literature.},
year = {2010}
}

@article{Mairhuber1956,
 ISSN = {00029939, 10886826},
 URL = {http://www.jstor.org/stable/2033359},
 author = {John C. Mairhuber},
 journal = {Proceedings of the American Mathematical Society},
 number = {4},
 pages = {609--615},
 publisher = {American Mathematical Society},
 title = {On Haar's Theorem Concerning Chebychev Approximation Problems Having Unique Solutions},
 volume = {7},
 year = {1956}
}

@article{curtis1959,
author = "Curtis, Philip C.",
fjournal = "Pacific Journal of Mathematics",
journal = "Pacific J. Math.",
number = "4",
pages = "1013--1027",
publisher = "Pacific Journal of Mathematics, A Non-profit Corporation",
title = "$n$-parameter families and best approximation.",
url = "https://projecteuclid.org:443/euclid.pjm/1103038880",
volume = "9",
year = "1959"
}

@article{Micchelli1986,
	Abstract = {Among other things, we prove that multiquadric surface interpolation is always solvable, thereby settling a conjecture of R. Franke.},
	Author = {Micchelli, Charles A. },
	Da = {1986/12/01},
	Date-Added = {2020-07-14 00:06:33 +0000},
	Date-Modified = {2020-07-14 00:06:33 +0000},
	Doi = {10.1007/BF01893414},
	Id = {Micchelli1986},
	Isbn = {1432-0940},
	Journal = {Constructive Approximation},
	Number = {1},
	Pages = {11--22},
	Title = {Interpolation of scattered data: Distance matrices and conditionally positive definite functions},
	Ty = {JOUR},
	Url = {https://doi.org/10.1007/BF01893414},
	Volume = {2},
	Year = {1986},
	Bdsk-Url-1 = {https://doi.org/10.1007/BF01893414},
	Bdsk-Url-2 = {http://dx.doi.org/10.1007/BF01893414}}

@article{Hardy1971,
author = {Hardy, Rolland L.},
title = {Multiquadric equations of topography and other irregular surfaces},
journal = {Journal of Geophysical Research (1896-1977)},
volume = {76},
number = {8},
pages = {1905-1915},
keywords = {General: Mathematics, Geomorphology: Theoretical Studies},
doi = {10.1029/JB076i008p01905},
url = {https://agupubs.onlinelibrary.wiley.com/doi/abs/10.1029/JB076i008p01905},
eprint = {https://agupubs.onlinelibrary.wiley.com/doi/pdf/10.1029/JB076i008p01905},
abstract = {A new analytical method of representing irregular surfaces that involves the summation of equations of quadric surfaces having unknown coefficients is described. The quadric surfaces are located at significant points throughout the region to be mapped. Procedures are given for solving multiquadric equations of topography that are based on coordinate data. Contoured multiquadric surfaces are compared with topography and other irregular surfaces from which the multiquadric equation was derived.},
year = {1971}
}

@article{fornberg_flyer_2015, title={Solving PDEs with radial basis functions}, volume={24}, DOI={10.1017/S0962492914000130}, journal={Acta Numerica}, publisher={Cambridge University Press}, author={Fornberg, Bengt and Flyer, Natasha}, year={2015}, pages={215–258}}

@article{Sarler2004,
author = {Sarler, Bozidar and Perko, Janez and Chen, Ching-Shyang},
year = {2004},
month = {03},
pages = {187-212},
title = {Radial basis function collocation method solution of natural convection in porous media},
volume = {14},
journal = {International Journal of Numerical Methods for Heat \& Fluid Flow},
doi = {10.1108/09615530410513809}
}
@article{10.1145/355984.355989,
author = {Paige, Christopher C. and Saunders, Michael A.},
title = {LSQR: An Algorithm for Sparse Linear Equations and Sparse Least Squares},
year = {1982},
issue_date = {March 1982},
publisher = {Association for Computing Machinery},
address = {New York, NY, USA},
volume = {8},
number = {1},
issn = {0098-3500},
url = {https://doi.org/10.1145/355984.355989},
doi = {10.1145/355984.355989},
journal = {ACM Trans. Math. Softw.},
month = mar,
pages = {43–71},
numpages = {29}
}

@article{Saunders1995,
	Abstract = {We examine two iterative methods for solving rectangular systems of linear equations: LSQR for over-determined systemsAx ≈b, and Craig's method for under-determined systemsAx = b. By including regularization, we extend Craig's method to incompatible systems, and observe that it solves the same damped least-squares problems as LSQR. The methods may therefore be compared on rectangular systems of arbitrary shape.},
	Author = {Saunders, Michael A. },
	Da = {1995/12/01},
	Date-Added = {2020-07-13 20:02:15 +0000},
	Date-Modified = {2020-07-13 20:02:15 +0000},
	Doi = {10.1007/BF01739829},
	Id = {Saunders1995},
	Isbn = {1572-9125},
	Journal = {BIT Numerical Mathematics},
	Number = {4},
	Pages = {588--604},
	Title = {Solution of sparse rectangular systems using LSQR and CRAIG},
	Ty = {JOUR},
	Url = {https://doi.org/10.1007/BF01739829},
	Volume = {35},
	Year = {1995},
	Bdsk-Url-1 = {https://doi.org/10.1007/BF01739829},
	Bdsk-Url-2 = {http://dx.doi.org/10.1007/BF01739829}}

@article{SOLEIMANI20101411,
title = "Local RBF-DQ method for two-dimensional transient heat conduction problems",
journal = "International Communications in Heat and Mass Transfer",
volume = "37",
number = "9",
pages = "1411 - 1418",
year = "2010",
issn = "0735-1933",
doi = "https://doi.org/10.1016/j.icheatmasstransfer.2010.06.033",
url = "http://www.sciencedirect.com/science/article/pii/S0735193310001806",
author = "Soheil Soleimani and M. Jalaal and H. Bararnia and E. Ghasemi and D.D. Ganji and F. Mohammadi",
keywords = "Conduction, Radial basis function, Mesh-free method, Irregular geometry",
abstract = "The meshless local radial basis function-based differential quadrature (RBF-DQ) method is applied on two-dimensional heat conduction for different irregular geometries. This method is the combination of differential quadrature approximation of derivatives and function approximation of radial basis function. Four different geometries with regular and irregular boundaries are considered, and numerical results are compared with those gained by finite element (FE) solution achieved by COMSOL commercial code. Outcomes prove that current technique is in very good agreement with FEM and this fact that RBF-DQ method is an accurate and flexible method in solution of heat conduction problems."
}

@article{SHU20052001,
title = "An upwind local RBF-DQ method for simulation of inviscid compressible flows",
journal = "Computer Methods in Applied Mechanics and Engineering",
volume = "194",
number = "18",
pages = "2001 - 2017",
year = "2005",
issn = "0045-7825",
doi = "https://doi.org/10.1016/j.cma.2004.07.008",
url = "http://www.sciencedirect.com/science/article/pii/S0045782504003640",
author = "C. Shu and H. Ding and H.Q. Chen and T.G. Wang",
keywords = "Radial basis function, Meshless method, Upwind scheme, Inviscid compressible flow",
abstract = "In this paper, an upwind local radial basis function-based differential quadrature (RBF-DQ) scheme is presented for simulation of inviscid compressible flows with shock wave. RBF-DQ is a naturally mesh-free method. The scheme consists of two parts. The first part is to use the local RBF-DQ method to discretize the Euler equation in conservative, differential form on a set of scattered nodes. The second part is to apply the upwind method to evaluate the flux at the mid-point between the reference knot and its supporting knots. The proposed scheme is validated by its application to simulate the supersonic flow in a symmetric, convergent channel and the shock tube problem. The obtained numerical results agree very well with the theoretical data."
}

@Article{cmes.2003.004.085,
AUTHOR = {N. Mai-Duy, T. Tran-Cong},
TITLE = {Indirect RBFN Method with Thin Plate Splines for Numerical Solution of Differential Equations},
JOURNAL = {Computer Modeling in Engineering  Sciences},
VOLUME = {4},
YEAR = {2003},
NUMBER = {1},
PAGES = {85--102},
URL = {http://www.techscience.com/CMES/v4n1/24803},
ISSN = {1526-1506},
ABSTRACT = {This paper reports a mesh-free Indirect Radial Basis Function Network method (IRBFN) using Thin Plate Splines (TPSs) for numerical solution of Differential Equations (DEs) in rectangular and curvilinear coordinates. The adjustable parameters required by the method are the number of centres, their positions and possibly the order of the TPS. The first and second order TPSs which are widely applied in numerical schemes for numerical solution of DEs are employed in this study. The advantage of the TPS over the multiquadric basis function is that the former, with a given order, does not contain the adjustable shape parameter (i.e. the RBF's width) and hence TPS-based RBFN methods require less parametric study. The direct TPS-RBFN method is also considered in some cases for the purpose of comparison with the indirect TPS-RBFN method. The TPS-IRBFN method is verified successfully with a series of problems including linear elliptic PDEs, nonlinear elliptic PDEs, parabolic PDEs and Navier-Stokes equations in rectangular and curvilinear coordinates. Numerical results obtained show that the method achieves the norm of the relative error of the solution of$O(10^{-6})$ for the case of 1D second order DEs using a density of$51$, of$O(10^{-7})$ for the case of 2D elliptic PDEs using a density of$20\times 20$ and a Reynolds number$Re=200$ for the case of Jeffery-Hamel flow with a density of$43\times 12$.},
DOI = {10.3970/cmes.2003.004.085}
}

@INPROCEEDINGS{Fasshauer97solvingpartial,
    author = {Gregory E. Fasshauer},
    title = {Solving Partial Differential Equations by Collocation with Radial Basis Functions},
    booktitle = {In: Surface Fitting and Multiresolution Methods A. Le M'ehaut'e, C. Rabut and L.L. Schumaker (eds.), Vanderbilt},
    year = {1997},
    pages = {131--138},
    publisher = {University Press}
}

@article{SARLER20061269,
title = "Meshfree explicit local radial basis function collocation method for diffusion problems",
journal = "Computers  Mathematics with Applications",
volume = "51",
number = "8",
pages = "1269 - 1282",
year = "2006",
note = "Radial Basis Functions and Related Multivariate Meshfree Approximation Methods: Theory and Applications",
issn = "0898-1221",
doi = "https://doi.org/10.1016/j.camwa.2006.04.013",
url = "http://www.sciencedirect.com/science/article/pii/S0898122106000836",
author = "B. Sarler and R. Vertnik",
abstract = "This paper formulates a simple explicit local version of the classical meshless radial basis function collocation (Kansa) method. The formulation copes with the diffusion equation, applicable in the solution of a broad spectrum of scientific and engineering problems. The method is structured on multiquadrics radial basis functions. Instead of global, the collocation is made locally over a set of overlapping domains of influence and the time-stepping is performed in an explicit way. Only small systems of linear equations with the dimension of the number of nodes included in the domain of influence have to be solved for each node. The computational effort thus grows roughly linearly with the number of the nodes. The developed approach thus overcomes the principal large-scale problem bottleneck of the original Kansa method. Two test cases are elaborated. The first is the boundary value problem (NAFEMS test) associated with the steady temperature field with simultaneous involvement of the Dirichlet, Neumann and Robin boundary conditions on a rectangle. The second is the initial value problem, associated with the Dirichlet jump problem on a square. The accuracy of the method is assessed in terms of the average and maximum errors with respect to the density of nodes, number of nodes in the domain of influence, multiquadrics free parameter, and timestep length on uniform and nonuniform node arrangements. The developed meshless method outperforms the classical finite difference method in terms of accuracy in all situations except immediately after the Dirichlet jump where the approximation properties appear similar."
}

@article{Kovacevic,
author = {Kovacevic, I and Poredos, A and Sarler, Bozidar},
year = {2003},
month = {12},
pages = {575-599},
title = {Solving the Stefan problem with the radial basis function collocation method},
volume = {44},
journal = {Numerical Heat Transfer: Part B: Fundamentals},
doi = {10.1080/716100496}
}

@article{doi:10.1002/fld.165,
author = {Mai-Duy, Nam and Tran-Cong, Thanh},
title = {Numerical solution of Navier–Stokes equations using multiquadric radial basis function networks},
journal = {International Journal for Numerical Methods in Fluids},
volume = {37},
number = {1},
pages = {65-86},
keywords = {mesh-free method, Navier–Stokes equations, radial basis function networks, streamfunction–vorticity formulation},
doi = {10.1002/fld.165},
url = {https://onlinelibrary.wiley.com/doi/abs/10.1002/fld.165},
eprint = {https://onlinelibrary.wiley.com/doi/pdf/10.1002/fld.165},
abstract = {Abstract A numerical method based on radial basis function networks (RBFNs) for solving steady incompressible viscous flow problems (including Boussinesq materials) is presented in this paper. The method uses a ‘universal approximator’ based on neural network methodology to represent the solutions. The method is easy to implement and does not require any kind of ‘finite element-type’ discretization of the domain and its boundary. Instead, two sets of random points distributed throughout the domain and on the boundary are required. The first set defines the centres of the RBFNs and the second defines the collocation points. The two sets of points can be different; however, experience shows that if the two sets are the same better results are obtained. In this work the two sets are identical and hence commonly referred to as the set of centres. Planar Poiseuille, driven cavity and natural convection flows are simulated to verify the method. The numerical solutions obtained using only relatively low densities of centres are in good agreement with analytical and benchmark solutions available in the literature. With uniformly distributed centres, the method achieves Reynolds number Re = 100 000 for the Poiseuille flow (assuming that laminar flow can be maintained) using the density of \$11\times 11\$, Re = 400 for the driven cavity flow with a density of \$33\times 33\$ and Rayleigh number Ra = 1 000 000 for the natural convection flow with a density of \$27\times 27\$. Copyright © 2001 John Wiley \& Sons, Ltd.},
year = {2001}
}

@Article{cmes.2005.007.185,
AUTHOR = {B. Sarler},
TITLE = {A Radial Basis Function Collocation Approach in Computational Fluid Dynamics},
JOURNAL = {Computer Modeling in Engineering \& Sciences},
VOLUME = {7},
YEAR = {2005},
NUMBER = {2},
PAGES = {185--194},
URL = {http://www.techscience.com/CMES/v7n2/29727},
ISSN = {1526-1506},
ABSTRACT = {This paper explores the application of the mesh-free radial basis function collocation method for solution of heat transfer and fluid flow problems. The solution procedure is represented for a Poisson reformulated general transport equation in terms of a-symmetric, symmetric and modified (double consideration of the boundary nodes) collocation approaches. In continuation, specifics of a primitive variable solution procedure for the coupled mass, momentum, and energy transport representing the natural convection in an incompressible Newtonian Bussinesq fluid are elaborated. A comparison of different collocation strategies is performed based on the two dimensional De Vahl Davis steady natural convection benchmark with Prandtl number Pr = 0.71, and Rayleigh numbers Ra = 10<sup>3</sup>, 10<sup>4</sup>, 10<sup>5</sup>, 10<sup>6</sup>. Multiquadrics radial basis functions are used. The three methods are assessed in terms of streamfunction extreme, cavity Nusselt number, and mid-plane velocity components. Best performance is achieved with the modified approach.},
DOI = {10.3970/cmes.2005.007.185}
}

@article{KANSA1990147,
title = "Multiquadrics—A scattered data approximation scheme with applications to computational fluid-dynamics—II solutions to parabolic, hyperbolic and elliptic partial differential equations",
journal = "Computers  Mathematics with Applications",
volume = "19",
number = "8",
pages = "147 - 161",
year = "1990",
issn = "0898-1221",
doi = "https://doi.org/10.1016/0898-1221(90)90271-K",
url = "http://www.sciencedirect.com/science/article/pii/089812219090271K",
author = "E.J. Kansa",
abstract = "This paper is the second in a series of investigations into the benefits of multiquadrics (MQ). MQ is a true scattered data, multidimensional spatial approximation scheme. In the previous paper, we saw that MQ was an extremely accurate approximation scheme for interpolation and partial derivative estimates for a variety of two-dimensional functions over both gridded and scattered data. The theory of Madych and Nelson shows for the space of all conditionally positive definite functions to which MQ belongs, a semi-norm exists which is minimized by such functions. In this paper, MQ is used as the spatial approximation scheme for parabolic, hyperbolic and the elliptic Poisson's equation. We show that MQ is not only exceptionally accurate, but is more efficient than finite difference schemes which require many more operations to achieve the same degree of accuracy."
}

@article{lee_local_2003,
	title = {Local multiquadric approximation for solving boundary value problems},
	volume = {30},
	issn = {1432-0924},
	url = {https://doi.org/10.1007/s00466-003-0416-5},
	doi = {10.1007/s00466-003-0416-5},
	abstract = {This paper presents a truly meshless approximation strategy for solving partial differential equations based on the local multiquadric (LMQ) and the local inverse multiquadric (LIMQ) approximations. It is different from the traditional global multiquadric (GMQ) approximation in such a way that it is a pure local procedure. In constructing the approximation function, the only geometrical data needed is the local configuration of nodes fallen within its influence domain. Besides this distinct characteristic of localization, in the context of meshless-typed approximation strategies, other major advantages of the present strategy include: (i) the existence of the shape functions is guaranteed provided that all the nodal points within an influence domain are distinct; (ii) the constructed shape functions strictly satisfy the Kronecker delta condition; (iii) the approximation is stable and insensitive to the free parameter embedded in the formulation and; (iv) the computational cost is modest and the matrix operations require only inversion of matrices of small size which is equal to the number of nodes inside the influence domain. Based on the present LMQ and LIMQ approximations, a collocation procedure is developed for solutions of 1D and 2D boundary value problems. Numerical results indicate that the present LMQ and LIMQ approximations are more stable than their global counterparts. In addition, it demonstrates that both approximation strategies are highly efficient and able to yield accurate solutions regardless of the chosen value for the free parameter.},
	number = {5},
	journal = {Computational Mechanics},
	author = {Lee, C. K. and Liu, X. and Fan, S. C.},
	month = apr,
	year = {2003},
	pages = {396--409}
}

@Article{Theillard:2017,
author ="Theillard, Maxime and Alonso-Matilla, Roberto and Saintillan, David",
title  ="Geometric control of active collective motion",
journal  ="Soft Matter",
year  ="2017",
volume  ="13",
issue  ="2",
pages  ="363-375",
publisher  ="The Royal Society of Chemistry",
doi  ="10.1039/C6SM01955B",
url  ="http://dx.doi.org/10.1039/C6SM01955B",
abstract  ="Recent experimental studies have shown that confinement can profoundly affect self-organization in semi-dilute active suspensions{,} leading to striking features such as the formation of steady and spontaneous vortices in circular domains and the emergence of unidirectional pumping motions in periodic racetrack geometries. Motivated by these findings{,} we analyze the two-dimensional dynamics in confined suspensions of active self-propelled swimmers using a mean-field kinetic theory where conservation equations for the particle configurations are coupled to the forced Navier–Stokes equations for the self-generated fluid flow. In circular domains{,} a systematic exploration of the parameter space casts light on three distinct states: equilibrium with no flow{,} stable vortex{,} and chaotic motion{,} and the transitions between these are explained and predicted quantitatively using a linearized theory. In periodic racetracks{,} similar transitions from equilibrium to net pumping to traveling waves to chaos are observed in agreement with experimental observations and are also explained theoretically. Our results underscore the subtle effects of geometry on the morphology and dynamics of emerging patterns in active suspensions and pave the way for the control of active collective motion in microfluidic devices."}

@Article{Theillard:2018,
author ="Theillard, Maxime and Saintillan, David",
title  ="Computational mean field modelling of confined active suspensions",
journal  ="in preparation",
year  ="2018",
}

@article{LUCY1977,
       author = {Lucy, L. B.},
        title = "{A numerical approach to the testing of the fission hypothesis.}",
     keywords = {Binary Stars, Hypotheses, Nuclear Fission, Numerical Analysis, Protostars, Stellar Evolution, Astronomical Models, Difference Equations, Gas Dynamics, Monte Carlo Method, Numerical Stability, Particle Motion, Stellar Rotation, Astrophysics},
         year = {1977},
       volume = {82},
        pages = {1013-1024},
          doi = {10.1086/112164},
      adsnote = {Provided by the SAO/NASA Astrophysics Data System}
}

@article{Vivek2020,
	abstract = {The present article outlines the historical development of Meshless method. The high computational cost associated with conventional numerical methods and their lack of accuracy, particularly in the area of astrophysics is the motivation and an initiation towards the invention of latest meshless methods. The review starts with Smoothed particle hydrodynamics meshless method which was invented by Gingold, Monghan in 1977 for astrophysics applications. Furthermore, shortcomings of the SPH method, it's solution, and recent development has been reviewed for different applications. After that advent features of Element free Galerkin method with mean least square method and Lagrangian multiplier have been critically reviewed for different cases. Later, the radial basis function method is explained along with its extended version to solve a partial differential equation. Many types of RBF methods have been compared. Also, several consistency and accuracy problems are reviewed for Reproducing Kernel Particle Method. Recent advancement like generalization of RKPM method, hybrid of RKPM with Lagrangian multiplier has been briefed for a different approach. Additionally, feasibility, flexibility and several applications of Meshless local Petrov Galerkin have been reviewed. Lastly recent development of MLPG with Laplace transform to solve PDE have been studied. Afterwards, a different attribute of point collocation method for inherent multiresolution capability and error estimation has been reviewed for different applications. Briefly, the article provides a condensed overview of several meshless methods, its recent advancement and their applications.},
	author = {Vivek G. Patel and Nikunj V. Rachchh},
	date-modified = {2024-03-05 19:23:54 -0800},
	doi = {https://doi.org/10.1016/j.matpr.2020.02.328},
	issn = {2214-7853},
	journal = {Materials Today: Proceedings},
	keywords = {Meshless methods, Element free galerkin, Radial basis function, Reproducing kernel particle method, Mean least square},
	note = {10th International Conference of Materials Processing and Characterization},
	pages = {1598-1603},
	read = {1},
	title = {Meshless method -- Review on recent developments},
	url = {https://www.sciencedirect.com/science/article/pii/S221478532031083X},
	volume = {26},
	year = {2020},
	bdsk-url-1 = {https://www.sciencedirect.com/science/article/pii/S221478532031083X},
	bdsk-url-2 = {https://doi.org/10.1016/j.matpr.2020.02.328}}

@article{Belyt1994,
author = {Belytschko, T. and Lu, Y. Y. and Gu, L.},
title = {Element-free Galerkin methods},
journal = {International Journal for Numerical Methods in Engineering},
volume = {37},
number = {2},
pages = {229-256},
doi = {https://doi.org/10.1002/nme.1620370205},
url = {https://onlinelibrary.wiley.com/doi/abs/10.1002/nme.1620370205},
eprint = {https://onlinelibrary.wiley.com/doi/pdf/10.1002/nme.1620370205},
abstract = {Abstract An element-free Galerkin method which is applicable to arbitrary shapes but requires only nodal data is applied to elasticity and heat conduction problems. In this method, moving least-squares interpolants are used to construct the trial and test functions for the variational principle (weak form); the dependent variable and its gradient are continuous in the entire domain. In contrast to an earlier formulation by Nayroles and coworkers, certain key differences are introduced in the implementation to increase its accuracy. The numerical examples in this paper show that with these modifications, the method does not exhibit any volumetric locking, the rate of convergence can exceed that of finite elements significantly and a high resolution of localized steep gradients can be achieved. The moving least-squares interpolants and the choices of the weight function are also discussed in this paper.},
year = {1994}
}

@article{JOLDES2019152,
	abstract = {The ability to predict patient-specific soft tissue deformations is key for computer-integrated surgery systems and the core enabling technology for a new era of personalized medicine. Element-Free Galerkin (EFG) methods are better suited for solving soft tissue deformation problems than the finite element method (FEM) due to their capability of handling large deformation while also eliminating the necessity of creating a complex predefined mesh. Nevertheless, meshless methods based on EFG formulation, exhibit three major limitations: (i) meshless shape functions using higher order basis cannot always be computed for arbitrarily distributed nodes (irregular node placement is crucial for facilitating automated discretization of complex geometries); (ii) imposition of the Essential Boundary Conditions (EBC) is not straightforward; and, (iii) numerical (Gauss) integration in space is not exact as meshless shape functions are not polynomial. This paper presents a suite of Meshless Total Lagrangian Explicit Dynamics (MTLED) algorithms incorporating a Modified Moving Least Squares (MMLS) method for interpolating scattered data both for visualization and for numerical computations of soft tissue deformation, a novel way of imposing EBC for explicit time integration, and an adaptive numerical integration procedure within the Meshless Total Lagrangian Explicit Dynamics algorithm. The appropriateness and effectiveness of the proposed methods is demonstrated using comparisons with the established non-linear procedures from commercial finite element software ABAQUS and experiments with very large deformations. To demonstrate the translational benefits of MTLED we also present a realistic brain-shift computation.},
	author = {Grand Joldes and George Bourantas and Benjamin Zwick and Habib Chowdhury and Adam Wittek and Sudip Agrawal and Konstantinos Mountris and Damon Hyde and Simon K. Warfield and Karol Miller},
	doi = {https://doi.org/10.1016/j.media.2019.06.004},
	issn = {1361-8415},
	journal = {Medical Image Analysis},
	keywords = {Surgical simulation, Soft tissues, Meshless Total Lagrangian Explicit Dynamics, Nonlinear computational mechanics},
	pages = {152-171},
	title = {Suite of meshless algorithms for accurate computation of soft tissue deformation for surgical simulation},
	url = {https://www.sciencedirect.com/science/article/pii/S1361841518305723},
	volume = {56},
	year = {2019},
	bdsk-url-1 = {https://www.sciencedirect.com/science/article/pii/S1361841518305723},
	bdsk-url-2 = {https://doi.org/10.1016/j.media.2019.06.004}}

@article{HE2020225,
	abstract = {In this paper, we develop a meshless method by using radial basis functions (RBFs) to solve time domain Maxwell's equations resulting from simulating wave propagation in electromagnetic wave splitter and rotator devices. To simulate wave propagation in these devices, we introduce a perfectly matched layer (PML) to reduce the unbounded physical domain problem to a bounded domain problem. Using PML leads to a multi-physics problem with different governing equations in different subdomains, which makes the simulation quite challenging. By following the idea of leap-frog FDTD scheme, we develop the leap-frog RBF meshless method to solve these complicated coupled modeling equations. Extensive numerical results by using multiquadric RBF and Gaussian RBF are presented to demonstrate the effectiveness of our meshless method.},
	author = {Bin He and Jichun Li and Meng Chen and Yunqing Huang},
	doi = {https://doi.org/10.1016/j.enganabound.2020.06.010},
	issn = {0955-7997},
	journal = {Engineering Analysis with Boundary Elements},
	keywords = {Maxwell's equations, Radial basis functions, Meshless method, Perfectly matched layer, Metamaterials},
	pages = {225-242},
	title = {A leap-frog meshless method with radial basis functions for simulating electromagnetic wave splitter and rotator},
	url = {https://www.sciencedirect.com/science/article/pii/S0955799720301661},
	volume = {118},
	year = {2020},
	bdsk-url-1 = {https://www.sciencedirect.com/science/article/pii/S0955799720301661},
	bdsk-url-2 = {https://doi.org/10.1016/j.enganabound.2020.06.010}}

@article{LU1995131,
	abstract = {Element-free Galerkin method (EFG) is extended to dynamic problems. EFG method, which is based on moving least square interpolants (MLS), requires only nodal data; no element connectivity is needed. This makes the method particularly attractive for moving dynamic crack problems, since remeshing can be avoided. In contrast to the earlier formulation for static problems by authors, the weak form of kinematic boundary conditions for dynamic problems is introduced in the implementation to enforce the kinematic boundary conditions. With this formulation, the stiffness matrix is symmetric and positive semi-definite, and hence the consistency, conergence and stability analyses of time integration remain the same as those in finite element method. Numerical examples are presented to illustrate the performance of this method. The relationship between the element-free Galerkin method and the smooth particle hydrodynamics (SPH) method is also discussed in this paper. Results are presented for some one-dimensional problems and two-dimensional problems with static and moving cracks.},
	author = {Y.Y. Lu and T. Belytschko and M. Tabbara},
	doi = {https://doi.org/10.1016/0045-7825(95)00804-A},
	issn = {0045-7825},
	journal = {Computer Methods in Applied Mechanics and Engineering},
	number = {1},
	pages = {131-153},
	title = {Element-free Galerkin method for wave propagation and dynamic fracture},
	url = {https://www.sciencedirect.com/science/article/pii/004578259500804A},
	volume = {126},
	year = {1995},
	bdsk-url-1 = {https://www.sciencedirect.com/science/article/pii/004578259500804A},
	bdsk-url-2 = {https://doi.org/10.1016/0045-7825(95)00804-A}}

@article{MAZHAR202114,
	abstract = {This paper presents a review of recent progress made towards the applications of the meshfree particle methods (MPMs) for solving coupled fluid-structure interaction (FSI) problems. Meshfree methods are categorized based on their mathematical formulation and treatment of computational data points. The advantages and limitations of these methods, particularly related to FSI applications, have been identified. A detailed account of salient work related to the FSI problems involving complex geometries, viscous flows, and large structural deformations has been presented and the benchmark solutions are identified for future research. Compared to their mesh-based counterparts, MPMs are found better suited in negotiating moving boundaries and complex geometries, features that are the hallmark of FSI problems. However, the biggest challenge to their wider acceptability is their implementation and programming complexity, higher computational cost, and lack of commercial software packages. So far, meshfree methods have mostly been limited to applications, where conventional methods show limited performance. Owing to its promising growth potential, partitioned FSI is the prime emphasis of this paper. Various aspects of partitioned FSI have been identified and classified for meshfree FSI problems, which include problem formulation strategies, domains discretization approaches, solver coupling methodology, interface treatment, benchmark problems, computational load, and availability of commercial software. Furthermore, various challenges involved in employing MPMs for FSI have also been identified and discussed along with the state-of-the-art techniques used in meshfree methods and FSI applications, and a future way forward has been proposed. In essence, this paper is an effort to identify and classify key aspects of MPM applications for FSI and suggest potential avenues to explore the full potential of MPM capabilities for the solution of coupled problems.},
	author = {Farrukh Mazhar and Ali Javed and Jing Tang Xing and Aamer Shahzad and Mohtashim Mansoor and Adnan Maqsood and Syed Irtiza Ali Shah and Kamran Asim},
	doi = {https://doi.org/10.1016/j.enganabound.2020.11.005},
	issn = {0955-7997},
	journal = {Engineering Analysis with Boundary Elements},
	keywords = {Meshfree particle methods, Meshless methods, Fluid-solid interaction, Fluid-structure interface, FSI coupling, Immersed and body-conforming meshes},
	pages = {14-40},
	title = {On the meshfree particle methods for fluid-structure interaction problems},
	url = {https://www.sciencedirect.com/science/article/pii/S0955799720302885},
	volume = {124},
	year = {2021},
	bdsk-url-1 = {https://www.sciencedirect.com/science/article/pii/S0955799720302885},
	bdsk-url-2 = {https://doi.org/10.1016/j.enganabound.2020.11.005}}

@article{Kumar2021,
	abstract = {In this work, authors are actively showcasing the finding of recent papers using meshless method in different structure problems which are either 1-D, 2-D, 3-D problems or problems related to composite material. Also presents the comparison of research going in that field of most common and widely used method, i.e. most widely used method of finite element (FEM). Over the past few decades meshless methods where implemented into many application areas which are ranging from long established problems related to astral physics to that of solid mechanics problem also in engineering and mathematical models, vibration and fluid mechanics investigation and optimization of the numerical results of the equation of partial differential problems. From critical reviews of different meshless methods, authors have suggested the applicability of different meshless method to improve computational efficiency.},
	author = {Manish Kumar and Aditya Kumar Jha and Yash Bhagoria and Pankaj Gupta},
	doi = {10.1088/1757-899X/1116/1/012119},
	journal = {IOP Conference Series: Materials Science and Engineering},
	month = {apr},
	number = {1},
	pages = {012119},
	publisher = {IOP Publishing},
	title = {A review to explore different meshless methods in various Structural problems},
	url = {https://dx.doi.org/10.1088/1757-899X/1116/1/012119},
	volume = {1116},
	year = {2021},
	bdsk-url-1 = {https://dx.doi.org/10.1088/1757-899X/1116/1/012119}}

@article{Antonio2023,
	abstract = {Mesh-based and particle methods were conceived as two different discretization strategies to solve partial differential equations. In the last two decades computational methods have diversified and a myriad of hybrid formulations that combine elements of these two approaches have been developed to solve Computational fluid dynamics problems. In this work we present a review about the meshless-FV family of methods, an analysis is carried out showing that the MLS-SPH-ALE method can be considered as a general formulation from which a set of particle-based methods can be recovered. Moreover, we show the relations between the MLS-SPH-ALE method and the finite volume method. The MLS-SPH-ALE method is a versatile particle-based method that was developed to circumvent the consistency issues of particle methods caused by the use of the kernel approximation. The MLS-SPH-ALE method is developed from the differential equation in ALE form using the partition unity property which is automatically fulfilled by the Moving Least Squares approximation.},
	author = {Eir{\'\i}s, Antonio and Ram{\'\i}rez, Luis and Couceiro, Iv{\'a}n and Fern{\'a}ndez-Fidalgo, Javier and Par{\'\i}s, Jos{\'e} and Nogueira, Xes{\'u}s},
	date = {2023/11/01},
	date-added = {2024-03-06 16:36:42 -0800},
	date-modified = {2024-03-06 16:36:42 -0800},
	doi = {10.1007/s11831-023-09965-2},
	id = {Eir{\'\i}s2023},
	isbn = {1886-1784},
	journal = {Archives of Computational Methods in Engineering},
	number = {8},
	pages = {4959--4981},
	title = {MLS-SPH-ALE: A Review of Meshless-FV Methods and a Unifying Formulation for Particle Discretizations},
	url = {https://doi.org/10.1007/s11831-023-09965-2},
	volume = {30},
	year = {2023},
	bdsk-url-1 = {https://doi.org/10.1007/s11831-023-09965-2}}

@article{KUMAR2019858,
	abstract = {We introduce an improved meshfree approximation scheme which is based on the local maximum-entropy strategy as a compromise between shape function locality and entropy in an information-theoretical sense. The improved version is specifically designed for severe, finite deformation and offers significantly enhanced stability as opposed to the original formulation. This is achieved by (i) formulating the quasistatic mechanical boundary value problem in a suitable updated-Lagrangian setting, (ii) introducing anisotropy in the shape function support to accommodate directional variations in nodal spacing with increasing deformation and eliminate tensile instability, (iii) spatially bounding and evolving shape function support to restrict the domain of influence and increase efficiency, (iv) truncating shape functions at interfaces in order to stably represent multi-component systems like composites or polycrystals. The new scheme is applied to benchmark problems of severe elastic and elastoplastic deformation that demonstrate its performance both in terms of accuracy (as compared to exact solutions and, where applicable, finite element simulations) and efficiency. Importantly, the presented formulation overcomes the classical tensile instability found in most meshfree interpolation schemes, as shown for stable simulations of, e.g., the inhomogeneous extension of a hyperelastic block up to 100% or the torsion of a hyperelastic cube by 200$\,^{\circ}$ --- both in an updated Lagrangian setting and without the need for remeshing.},
	author = {Siddhant Kumar and Kostas Danas and Dennis M. Kochmann},
	doi = {https://doi.org/10.1016/j.cma.2018.10.030},
	issn = {0045-7825},
	journal = {Computer Methods in Applied Mechanics and Engineering},
	keywords = {Meshfree/particle methods, Updated Lagrangian, Maximum entropy, Tensile instability, Finite deformation},
	pages = {858-886},
	title = {Enhanced local maximum-entropy approximation for stable meshfree simulations},
	url = {https://www.sciencedirect.com/science/article/pii/S0045782518305346},
	volume = {344},
	year = {2019},
	bdsk-url-1 = {https://www.sciencedirect.com/science/article/pii/S0045782518305346},
	bdsk-url-2 = {https://doi.org/10.1016/j.cma.2018.10.030}}

@article{DUAN200866,
	abstract = {In this paper, we prove some inverse inequalities for the finite dimensional subspace interpolated by radial basis functions. The properties of these inequalities are worse than those in finite element spaces. These bad properties result in a huge condition number for both the collocation method and Galerkin method using radial basis functions. As applications, theoretical estimates of the condition numbers of a meshless method using radial basis functions are given.},
	author = {Yong Duan},
	doi = {https://doi.org/10.1016/j.camwa.2007.03.011},
	issn = {0898-1221},
	journal = {Computers \& Mathematics with Applications},
	keywords = {Radial basis function, Inverse inequalities, Collocation method, Galerkin method, Meshless},
	number = {1},
	pages = {66-75},
	title = {A note on the meshless method using radial basis functions},
	url = {https://www.sciencedirect.com/science/article/pii/S089812210700332X},
	volume = {55},
	year = {2008},
	bdsk-url-1 = {https://www.sciencedirect.com/science/article/pii/S089812210700332X},
	bdsk-url-2 = {https://doi.org/10.1016/j.camwa.2007.03.011}}

@article{KANSA2000123,
	abstract = {Madych and Nelson [1] proved multiquadric (MQ) mesh-independent radial basis functions (RBFs) enjoy exponential convergence. The primary disadvantage of the MQ scheme is that it is global, hence, the coefficient matrices obtained from this discretization scheme are full. Full matrices tend to become progressively more ill-conditioned as the rank increases. In this paper, we explore several techniques, each of which improves the conditioning of the coefficient matrix and the solution accuracy. The methods that were investigated are 1.(1) replacement of global solvers by block partitioning, LU decomposition schemes,2.(2) matrix preconditioners,3.(3) variable MQ shape parameters based upon the local radius of curvature of the function being solved,4.(4) a truncated MQ basis function having a finite, rather than a full band-width,5.(5) multizone methods for large simulation problems, and6.(6) knot adaptivity that minimizes the total number of knots required in a simulation problem. The hybrid combination of these methods contribute to very accurate solutions. Even though FEM gives rise to sparse coefficient matrices, these matrices in practice can become very ill-conditioned. We recommend using what has been learned from the FEM practitioners and combining their methods with what has been learned in RBF simulations to form a flexible, hybrid approach to solve complex multidimensional problems.},
	author = {E.J. Kansa and Y.C. Hon},
	doi = {https://doi.org/10.1016/S0898-1221(00)00071-7},
	issn = {0898-1221},
	journal = {Computers \& Mathematics with Applications},
	keywords = {Multiquadric radial basis functions applied to PDEs, Two-dimensional Poisson elliptic partial differential equations, Domain decomposition methods for systems of linear equations, Truncated multiquadric radial basis functions, Multizone decomposition methods},
	number = {7},
	pages = {123-137},
	title = {Circumventing the ill-conditioning problem with multiquadric radial basis functions: Applications to elliptic partial differential equations},
	url = {https://www.sciencedirect.com/science/article/pii/S0898122100000717},
	volume = {39},
	year = {2000},
	bdsk-url-1 = {https://www.sciencedirect.com/science/article/pii/S0898122100000717},
	bdsk-url-2 = {https://doi.org/10.1016/S0898-1221(00)00071-7}}

@article{Franke1982ScatteredDI,
  title={Scattered data interpolation: tests of some methods},
  author={Richard Franke},
  journal={Mathematics of Computation},
  year={1982},
  volume={38},
  pages={181-200},
  url={https://api.semanticscholar.org/CorpusID:8290519}
}

@article{Chung1981,
	author = {K. C. Chung},
	doi = {10.1080/01495728108961797},
	eprint = {https://doi.org/10.1080/01495728108961797},
	journal = {Numerical Heat Transfer},
	number = {3},
	pages = {345-357},
	publisher = {Taylor \& Francis},
	title = {A GENERALIZED FINITE-DIFFERENCE METHOD FOR HEAT TRANSFER PROBLEMS OF IRREGULAR GEOMETRIES},
	url = {https://doi.org/10.1080/01495728108961797},
	volume = {4},
	year = {1981},
	bdsk-url-1 = {https://doi.org/10.1080/01495728108961797}}

\end{document}